\documentclass[11pt,reqno]{amsart} 
\usepackage[utf8]{inputenc}
\usepackage{amsmath,amssymb,amsthm,mathtools,microtype}
\usepackage[letterpaper,margin=1in]{geometry} 
\usepackage{setspace}
\usepackage{needspace}
\usepackage[shortlabels]{enumitem}
\usepackage{algorithm}
\usepackage[foot]{amsaddr}
\makeatletter
\renewcommand{\@setemails}{\begingroup\raggedright\mbox{\emailaddrname:\space}{\ttfamily\emails}.\par\endgroup}
\makeatother
\usepackage{xcolor}
\usepackage{tikz}
\usepackage{flafter}
\definecolor{darkblue}{rgb}{0,0,0.5}
\usepackage[colorlinks=true,linkcolor=darkblue,citecolor=darkblue,urlcolor=darkblue]{hyperref}
\usepackage[capitalize,noabbrev]{cleveref}

\usepackage[backend=biber, style=alphabetic, sorting=nyt, maxnames=100, backref=true, giveninits=true, sortcites=true]{biblatex} 
\hypersetup{pdftitle={Fractional colorings of K\_r-free graphs with local demands},pdfauthor={Ewan Davies, Abhishek Dhawan, Abhishek Methuku, Huy Nguyen, Jaehyeon Seo}}
\hypersetup{bookmarksnumbered=true, bookmarksopen=true, bookmarksopenlevel=1}

\makeatletter
\renewcommand\section{\@startsection{section}{1}%
  \z@{1.2\linespacing\@plus\linespacing}{.5\linespacing}
  {\normalfont\scshape\centering}}%
\makeatother

\DeclareFieldFormat{eprint:arxiv}{\textsc{arXiv}: \href{https://arxiv.org/abs/#1}{\nolinkurl{#1}}}
\DeclareFieldFormat{eprint}{\textsc{arXiv}: \href{https://arxiv.org/abs/#1}{\nolinkurl{#1}}}
\AtEveryBibitem{\clearfield{issn}\clearfield{isbn}\clearlist{language}}

\renewbibmacro*{doi+eprint+url}{%
    \iftoggle{bbx:eprint}{%
		\usebibmacro{eprint}%
	}{}
    \iffieldundef{eprint}{%
        \printfield{doi}%
    }{%
    }
	\newunit\newblock%
	\iffieldundef{doi}{%
    \iffieldundef{eprint}{%
		\usebibmacro{url+urldate}%
	}{}}{}%
}

\makeatletter
\newenvironment{breakablealgorithm}
  {
     \refstepcounter{algorithm}
     \hrule height.8pt depth0pt \kern2pt
     \renewcommand{\caption}[2][\relax]{
       {\raggedright\textbf{\ALG@name~\thealgorithm} ##2\par}%
       \ifx\relax##1\relax 
         \addcontentsline{loa}{algorithm}{\protect\numberline{\thealgorithm}##2}%
       \else 
         \addcontentsline{loa}{algorithm}{\protect\numberline{\thealgorithm}##1}%
       \fi
       \kern2pt\hrule\kern2pt
       \smallskip
     }
  }{
      \smallskip
      \kern2pt\hrule\relax
  }
\makeatother

\setlist{nosep,leftmargin=2em}
\setlist[enumerate]{font=\upshape}
\numberwithin{equation}{section}
\newtheorem{theorem}{Theorem}[section]
\newtheorem{lemma}[theorem]{Lemma}
\newtheorem{corollary}[theorem]{Corollary}
\newtheorem{conjecture}[theorem]{Conjecture}
\newtheorem{question}[theorem]{Question}
\theoremstyle{definition}
\newtheorem{definition}[theorem]{Definition}
\theoremstyle{remark}

\crefname{conjecture}{conjecture}{conjectures}
\crefname{question}{question}{questions}

\AddToHook{env/lemma/begin}{\crefalias{theorem}{lemma}}
\AddToHook{env/proposition/begin}{\crefalias{theorem}{proposition}}
\AddToHook{env/corollary/begin}{\crefalias{theorem}{corollary}}
\AddToHook{env/conjecture/begin}{\crefalias{theorem}{conjecture}}
\AddToHook{env/claim/begin}{\crefalias{theorem}{claim}}
\AddToHook{env/question/begin}{\crefalias{theorem}{question}}

\newcommand{\step}[1]{\par\smallskip\noindent\textit{#1}}

\newcommand{\bone}{\mathbf{1}}
\newcommand{\ind}{\mathbf{1}}
\newcommand{\cI}{\mathcal{I}}
\newcommand{\defeq}{\coloneqq}
\newcommand{\eqdef}{\eqqcolon}
\newcommand{\ratiobound}[1]{R_{#1}}
\newcommand{\rtconstant}{\kappa_{\mathrm{rt}}}
\newcommand{\rtdegree}{\nu_{\mathrm{rt}}}
\newcommand{\rtfactor}{\mathfrak{B}_{n,r}}
\newcommand{\rtcontrol}{\mathfrak{C}}
\newcommand{\rtpower}{\mathfrak{p}}
\DeclareMathOperator{\E}{\mathbb{E}}
\DeclareMathOperator{\Prob}{\mathbb{P}}
\DeclareMathOperator{\TV}{TV}
\DeclareMathOperator{\KL}{KL}
\DeclarePairedDelimiter{\abs}{\lvert}{\rvert}
\DeclarePairedDelimiter{\norm}{\lVert}{\rVert}

\allowdisplaybreaks[1]

\title[Fractional colorings of $K_r$-free graphs with local demands]{Fractional colorings of $K_r$-free graphs with local demands}
\author{Ewan Davies$^\star$}
\address{$^\star$Department of Computer Science, Colorado State University}
\email{research@ewandavies.org}
\author{Abhishek Dhawan$^\dagger$}
\address{$^\dagger$Department of Mathematics, University of Illinois Urbana--Champaign}
\email{adhawan2@illinois.edu}
\author{Abhishek Methuku$^\dagger$}
\email{methuku@illinois.edu}
\author{Huy Nguyen$^\dagger$}
\email{huyn4@illinois.edu}
\author{Jaehyeon Seo$^\ddagger$}
\address{$^\ddagger$Department of Mathematics, Yonsei University}
\email{jaehyeonseo@yonsei.ac.kr}
\date{}

\begin{document}

\begin{abstract}
We study fractional coloring with local demands, a framework introduced by Kelly and Postle in which one seeks a probability distribution on the independent sets of a graph such that each vertex \(v\) is included with probability at least a prescribed target \(\phi(v)\).
For \(K_r\)-free graphs, Kelly and Postle conjectured that such a distribution exists provided \(\phi(v)\le c_r\frac{\log d(v)}{d(v)}\) for a fixed \(c_r>0\), where \(d(v)\) is the degree of \(v\), and they established this up to a factor of \((\log\log d(v))^2\).
In this paper, we prove their conjecture, showing that
\[
    \phi(v)\le\frac{c}{r}\frac{\log(d(v)+2)}{d(v)+1}
\]
suffices for an absolute constant \(c>0\).
As a corollary, we obtain a lower bound on the independence number of \(K_r\)-free graphs in terms of their degree sequence.
This strengthens the bound that OpenAI proved in resolving the Ajtai--Erd\H{o}s--Koml\'os--Szemer\'edi conjecture.
Moreover, our argument yields a deterministic algorithm that produces the desired probability distribution in time \(2^{O(r^2)}n^{O(1)}\), where \(n\) is the number of vertices and the implied constants are absolute.
\end{abstract}

\maketitle

\pagenumbering{gobble}
\tableofcontents
\newpage
\pagenumbering{arabic}

\section{Introduction}\label{section: intro}

\subsection{Background and main results}\label{subsection: background}

Let $G$ be a graph.
A \emph{proper coloring} of $G$ assigns colors to its vertices so that no two adjacent vertices receive the same color.
The \emph{chromatic number} $\chi(G)$ is the minimum number of colors in a proper coloring of $G$.
An \emph{independent set} is a set of pairwise nonadjacent vertices; we write $\cI(G)$ for the family of independent sets of $G$, and the \emph{independence number} $\alpha(G)$ is the maximum size of a set in $\cI(G)$.
These parameters are among the most central and classical objects of study in graph theory.

We write $d_G(v)$ for the degree of a vertex $v$.
For any $n$-vertex graph $G$ with maximum degree $\Delta$, elementary greedy arguments yield the well-known bounds \[ \alpha(G)\ge\frac{n}{\Delta+1} \qquad\text{and}\qquad \chi(G)\le\Delta+1. \]
A central line of inquiry in extremal graph theory asks to what extent local sparsity constraints allow substantial improvements over these greedy baselines.
The most basic such constraint is the absence of small cliques.
Let $K_r$ denote the complete graph on $r$ vertices; a graph is \emph{$K_r$-free} if it contains no copy of $K_r$, and \emph{triangle-free} if it is $K_3$-free.
Throughout, $\log$ denotes the natural logarithm, and $h(d)\sim k(d)$ means that $h(d)/k(d)\to1$ as $d\to\infty$.

For independent sets, this line of research began with the groundbreaking work of Ajtai, Koml\'os, and Szemer\'edi~\cite{ajtai1980note, ajtai1981dense}, who proved that every $n$-vertex triangle-free graph with average degree $d$ satisfies $\alpha(G)=\Omega\left(\frac{n\log d}{d}\right)$.
Shearer~\cite{shearer1983note} subsequently showed that every such graph satisfies \begin{equation}\label{eq: shearer function} \alpha(G)\ge n f(d), \qquad f(d)\defeq\frac{d\log d-d+1}{(d-1)^2}\quad(d>0,\ d\ne1), \end{equation} where $f(0)\defeq1$ and $f(1)\defeq1/2$.
We call $f$ \emph{Shearer's function}.
Since $f(d)\sim\log d/d$ as $d\to\infty$, Shearer's bound sharpens the leading constant to $1-o(1)$, which remains the state of the art.

For graph coloring, Kim~\cite{kim1995brooks} established that every graph with maximum degree $\Delta$ and girth at least $5$ satisfies $\chi(G)\le(1+o(1))\frac{\Delta}{\log\Delta}$.
Independently, Johansson~\cite{Joh_triangle} showed that triangle-free graphs satisfy $\chi(G)=O\left(\frac{\Delta}{\log\Delta}\right)$.\footnote{Although widely cited as a DIMACS Technical Report, Johansson's original manuscript is no longer easily accessible. A complete exposition by Molloy and Reed~\cite{MolloyReed} proves the explicit bound $\chi(G)\le160\frac{\Delta}{\log\Delta}$ for sufficiently large $\Delta$.}
Molloy~\cite{molloy2019list} later matched Kim's leading constant for triangle-free graphs, bringing the upper bound to $(1+o(1))\frac{\Delta}{\log\Delta}$; alternative proofs and generalizations have since been developed by several authors~\cite{bernshteyn2019johansson, hurley2021first, DKPS, bonamy2022bounding, martinsson2021simplified}.
More broadly, coloring under local sparsity conditions---where vertex neighborhoods span few edges---has grown into a vibrant subfield of graph theory~\cite{AKSConjecture, vu2002general, anderson2024coloring, anderson2025coloring, bradshaw2025toward, KangKelly2023nibble}.
The algorithmic aspects of these probabilistic methods are also of interest.
Many of them yield constructions that run in polynomial time when the maximum degree is fixed, and running times that are also polynomial in the degree have been obtained in some settings~\cite{shearer1983note,DKPS20}.
Pettie and Su~\cite{PS15} studied the distributed setting, refining Johansson's method and improving the leading constant for triangle-free graphs.

For $K_r$-free graphs with fixed $r\ge4$, the corresponding bounds remained less well understood for many years.
Shearer~\cite{shearer1995independence} established, using methods somewhat different from his analysis of the triangle-free case, that every $n$-vertex $K_r$-free graph with average degree $d$ satisfies \begin{equation}\label{eq: shearer} \alpha(G)=\Omega\biggl(\frac{n\log d}{rd\log\log d}\biggr). \end{equation}
Extending this to coloring, Johansson~\cite{J96-Kr} proved\footnote{As with his triangle-free paper, Johansson's manuscript was circulated as a DIMACS Technical Report and is no longer easily accessible; Bansal, Gupta, and Guruganesh~\cite[Section~9]{bansal2015lovasz} provide a simplified exposition of the argument.} that every $K_r$-free graph with maximum degree $\Delta$ satisfies \begin{equation}\label{eq: johansson} \chi(G)=O\biggl(\frac{r\Delta\log\log\Delta}{\log\Delta}\biggr). \end{equation}
See~\cite{molloy2019list, bernshteyn2019johansson, dhawan2025bounds} for simplified proofs and improved leading constants, and~\cite{bonamy2022bounding, DKPS} for asymptotic improvements when $r=\omega\bigl(\frac{\log\Delta}{\log\log\Delta}\bigr)$.
Ajtai, Erd\H{o}s, Koml\'os, and Szemer\'edi~\cite{ajtai1981turan} conjectured that the $\log\log$ factor in~\eqref{eq: shearer} can be removed.

\begin{conjecture}\label{conj:aeks}
    For every fixed $r\ge3$, there is a constant $c_r>0$ such that every $n$-vertex $K_r$-free graph $G$ with average degree $d\ge2$ satisfies \[ \alpha(G)\ge c_r \frac{n\log d}{d}. \]
\end{conjecture}

The analogous improvement of Johansson's bound~\eqref{eq: johansson} was conjectured by Alon, Krivelevich, and Sudakov~\cite{AKSConjecture}.
These conjectures have driven more than four decades of research~\cite{alon1996independence, alon1999coloring, davies2020graph, anderson2025coloring, AndersonBernshteynDhawan, dhawan2025independent, molloy2019list, bernshteyn2019johansson}, in part because of their connection to the off-diagonal Ramsey numbers $R(r,t)$.
Recently, \Cref{conj:aeks} was established in~\cite{OAI2026logarithmic}, and Kang~\cite{kang2026ode} refined the argument to obtain $\alpha(G)=\Omega\left(\frac{n\log d}{rd}\right)$ with an absolute implied constant.

In this work, we focus on the \emph{fractional chromatic number} $\chi_f(G)$, the optimal value of the linear programming relaxation of the standard integer program for $\chi(G)$.
It is often fruitful to view fractional colorings through a probabilistic lens: for a nonempty graph $G$, $\chi_f(G)$ is the least $q\ge1$ such that there is a probability distribution on $\cI(G)$ whose sampled set $I$ satisfies $\Prob(v\in I)\ge1/q$ for every $v\in V(G)$.
Many upper bounds on $\chi_f(G)$ are proved by randomized constructions of independent sets that exploit this characterization.
Since such a distribution has expected size at least $n/q$, we have $n/\alpha(G)\le\chi_f(G)\le\chi(G)$.

There are two major lines of inquiry concerning the fractional chromatic number.
First, one can often establish stronger upper bounds for $\chi_f(G)$ in regimes where analogous bounds for $\chi(G)$ are either false or currently out of reach.
For instance, Reed's Conjecture~\cite{reed1998omega} is known to hold for the fractional chromatic number (see~\cite[Ch.~21]{MolloyReed}), and Kilakos and Reed~\cite{kilakos1993fractionally} proved the fractional relaxation of the Total Coloring Conjecture~\cite{vizing1968some, behzad1965graphs}.
In fact, in the foundational paper introducing fractional coloring, which predates the proof of the Four Color Theorem~\cite{appel1976every}, Hilton, Rado, and Scott~\cite{hilton19735} proved that every planar graph has fractional chromatic number strictly less than $5$.
More recently, Martinsson and Steiner~\cite{martinsson2025random} showed that $\chi_f(G)=O(d/\log d)$ for $d$-degenerate triangle-free graphs $G$; since there exist $d$-degenerate triangle-free graphs with chromatic number $d+1$~\cite{kostochka1999properties, alon1999coloring}, this exhibits a sharp separation in the fractional setting.
An alternative proof and several extensions of this result were obtained in~\cite{dhawan2026fractional}.

Second, extending classical bounds on the independence number to the fractional chromatic number has attracted significant attention.
Dvo\v{r}\'ak, Sereni, and Volec~\cite{dvovrak2014subcubic} established that triangle-free subcubic graphs have fractional chromatic number at most $14/5$, generalizing a classic theorem of Staton~\cite{staton1979some} and resolving a conjecture of Heckman and Thomas~\cite{heckman2001new}.
Heckman and Thomas~\cite{heckman2006independent} also proved that triangle-free subcubic planar graphs satisfy $\alpha(G)\ge3n/8$, and conjectured that this bound extends to the fractional chromatic number.

A stronger form of the bounds above asks for an estimate in terms of the degree sequence rather than the average or maximum degree.
For triangle-free graphs, Shearer~\cite{shearer1991note} proved such a bound: \[ \alpha(G)\ge\sum_{v\in V(G)} f(d_G(v)), \] where $f$ is Shearer's function defined in~\eqref{eq: shearer function}.
Motivated by degree-sequence bounds of this kind, Kelly and Postle~\cite{kelly2024fractional} introduced the theory of \emph{fractional coloring with local demands}, in which each vertex $v$ requires a prescribed probability of belonging to the sampled independent set.
Formally, we adopt the following definition (cf.~\cite[Definition~1.2 and Proposition~1.4]{kelly2024fractional}).

\begin{definition}[$\phi$-colorings]\label[definition]{def:demand}
A \emph{demand function} for a graph $G$ is a function $\phi\colon V(G)\to[0,1]$.
A \emph{$\phi$-coloring} of $G$ is a probability distribution on $\cI(G)$ such that a set $I$ sampled from it satisfies $\Prob(v\in I)\ge\phi(v)$ for every $v\in V(G)$.
\end{definition}

The constant demand function $\phi\equiv1/q$ recovers the probabilistic description of $\chi_f(G)$ above.
Local demands provide the flexibility needed to sample larger independent sets in graphs with nonuniform degrees.
For instance, a vertex of minimum degree is a better candidate for inclusion than one of maximum degree, as including it removes fewer vertices from consideration.
Kelly and Postle conjectured a local-demand version of Shearer's bound for triangle-free graphs and made progress toward it; Martinsson and Steiner subsequently proved this conjecture.

\begin{theorem}[{Martinsson--Steiner~\cite{martinsson2025random}}]\label{theo:martinsson_steiner_local_demands}
    Let $G$ be a triangle-free graph without isolated vertices and let $\phi$ be a demand function for $G$. Then $G$ admits a $\phi$-coloring provided that \[ \phi(v)\le(1-o(1))\frac{\log d_G(v)}{d_G(v)}\qquad\text{for all }v\in V(G). \]
\end{theorem}

\Cref{theo:martinsson_steiner_local_demands} is proved via an intricate weighted modification of Shearer-style induction.
Davies~\cite{davies2026random} extended this argument to graphs whose neighborhoods have bounded maximum average degree or bounded (fractional) chromatic number.

For $K_r$-free graphs with $r\ge4$, Kelly and Postle proposed the following strengthening of the \hyperref[conj:aeks]{Ajtai--Erd\H{o}s--Koml\'os--Szemer\'edi Conjecture} and its degree-sequence form.

\begin{conjecture}[{Kelly--Postle~\cite[Conjecture~2.4]{kelly2024fractional}}]\label{conj:kelly-postle}
For every integer $r\ge4$, there is a constant $a_r>0$ such that every $K_r$-free graph $G$ admits a $\phi$-coloring for \[ \phi(v)\defeq a_r\frac{\log(d_G(v)+2)}{d_G(v)+1}\qquad(v\in V(G)). \]
\end{conjecture}

Taking the expected size of the sampled independent set shows that \Cref{conj:kelly-postle} implies a degree-sequence bound on $\alpha(G)$, which in turn implies \Cref{conj:aeks}.
Toward this conjecture, Kelly and Postle~\cite[Theorem~2.5]{kelly2024fractional} proved a weaker bound with an additional factor of $(\log\log d_G(v))^2$ in the denominator, for vertices of degree sufficiently large in terms of $r$.

Our main result confirms \Cref{conj:kelly-postle}, with $a_r$ inversely proportional to $r$.

\begin{theorem}\label{thm:main}
There is an absolute constant $\kappa>0$ such that the following holds for every integer $r\ge4$.
Every $K_r$-free graph $G$ admits a $\phi$-coloring for every demand function $\phi$ satisfying \[ \phi(v)\le\frac{\kappa}{r}\frac{\log(d_G(v)+2)}{d_G(v)+1}\qquad(v\in V(G)). \]
\end{theorem}

\Cref{thm:main} holds at every degree, including degree zero, and requires no minimum-degree assumption.
Taking the expected size of the sampled set, we immediately obtain a degree-sequence strengthening of \Cref{conj:aeks}.

\begin{corollary}\label{cor:global}
Let $\kappa$ be as in \Cref{thm:main}, let $r\ge4$, and let $G$ be a $K_r$-free graph. Then \[ \alpha(G)\ge\frac{\kappa}{r}\sum_{v\in V(G)}\frac{\log(d_G(v)+2)}{d_G(v)+1}. \]
\end{corollary}

This bound generalizes Shearer's degree-sequence bound~\cite{shearer1991note} to $K_r$-free graphs, and its dependence on $r$ matches that of Kang~\cite{kang2026ode}.
We deduce \Cref{thm:main} from a weighted bound, described next, and give a deterministic construction of the required distribution in time \(2^{O(r^2)}n^{O(1)}\), with absolute implied constants.\footnote{Running time bounds stated in terms of the number \(n\) of vertices assume \(n\ge1\). 
}

\subsection{Weighted independent sets}\label{subsection: weighted independent sets intro}
A \emph{vertex weight} is a real-valued function on the vertex set.
It is positive or nonnegative when each of its values is positive or nonnegative, respectively.
For a vertex weight $a$ and a set $A\subseteq V(G)$, write $a(A)\defeq\sum_{v\in A}a(v)$.
For a positive vertex weight $s$, define the \emph{weighted degree} \[ d_{G,s}(v)\defeq\frac{s(N_G(v))}{s(v)}, \] where $N_G(v)$ is the set of neighbors of $v$.
Taking $s\equiv\bone$ gives the ordinary degree $d_G(v)=|N_G(v)|$.
For a nonnegative vertex weight $u$, put \[ \alpha_u(G)\defeq\max\{u(I):I\in\cI(G)\}, \qquad \Phi_{G,u}(s)\defeq\sum_{v\in V(G)}u(v)f(d_{G,s}(v)), \] where $f$ is Shearer's function.
The two vertex weights play separate roles: $s$ determines the weighted degrees in the desired bound, and $u$ specifies the weight of the independent set we seek.
We call $u$ the \emph{objective weight} to distinguish these roles.

\begin{theorem}\label{thm:existence}
For every integer $r\ge4$, there is a rational constant $c_r\in(0,1]$ such that the following holds.
Let $G$ be a $K_r$-free graph, let $s$ be a positive real vertex weight, and let $u$ be a nonnegative real vertex weight.
Then \[ \alpha_u(G)\ge c_r\Phi_{G,u}(s). \]
Moreover, $c_r$ can be chosen so that $c_r\ge\kappa_1/r$ for an absolute constant $\kappa_1>0$.
\end{theorem}

Since the bound in \Cref{thm:existence} holds for every nonnegative $u$ with $s$ fixed, it yields a distribution meeting local demands.
Indeed, by linear programming duality, a graph $G$ admits a $\phi$-coloring if and only if $\alpha_u(G)\ge\sum_{v}u(v)\phi(v)$ for every nonnegative vertex weight $u$ (see \Cref{lem:demands}).
Applying this criterion with $\phi(v)\defeq c_rf(d_{G,s}(v))$ gives the following.

\begin{corollary}\label{cor:local}
Let $r\ge4$, let $G$ be a $K_r$-free graph, and let $s$ be a positive real vertex weight.
Then $G$ admits a $\phi$-coloring for \[ \phi(v)\defeq c_r f(d_{G,s}(v)) \qquad(v\in V(G)). \]
\end{corollary}

Since $f(d)(d+1)/\log(d+2)$ is bounded below by a positive absolute constant for $d\ge0$, taking $s\equiv\bone$ in \Cref{cor:local} yields \Cref{thm:main}.
Convexity and monotonicity of $f$ also give the average-degree consequence \begin{equation}\label{eq:global-consequences} \alpha(G)\ge c_r\sum_{v\in V(G)} f(d_G(v))\ge c_r|V(G)|f(d) \end{equation} for every nonempty $K_r$-free graph of average degree $d$.
This gives an alternative proof of \Cref{conj:aeks} with $c_r=\Omega(1/r)$, matching the dependence on $r$ obtained by Kang~\cite{kang2026ode}.

\subsection{Algorithmic consequences}
Existence alone does not yield an efficient construction, since a graph may have exponentially many independent sets.
For rational input weights, we give a deterministic algorithm that finds an independent set meeting the weighted bound up to a factor of two, and we then use it to construct a distribution meeting local demands.
In particular, the improvement over the greedy bound by a factor of order \(\log d\) is achieved by a deterministic algorithm whose running time is fixed-parameter tractable in \(r\).
The input length \(L\) counts the bits needed to encode the graph by its adjacency matrix and all rational vertex weights, with each numerator and denominator written in binary; in particular, \(L\ge\max\{n,1\}\).
Running times count bit operations, with all arithmetic performed exactly on rational numbers in binary.
The algorithms below take \(r\ge4\) as part of the input, and there are absolute constants \(\rtconstant,\rtdegree>0\) for which their running times satisfy the stated bounds; in particular, the polynomial degree is independent of \(r\).
We do not optimize these constants.

\begin{theorem}\label{thm:algorithm}
There is a deterministic algorithm that, given an integer \(r\ge4\), a \(K_r\)-free graph \(G\), and rational vertex weights \(s>0\) and \(u\ge0\), returns an independent set \(I\) satisfying
\[
    u(I)\ge\frac{c_r}{2}\Phi_{G,u}(s)
\]
in time at most \(2^{\rtconstant r^2}L^{\rtdegree}\).
\end{theorem}

To construct a local-demand distribution, we give a multiplicative-weights procedure with rational updates, following the fractional covering framework~\cite[Section~3.3.3]{arora2012multiplicative}.
For an \(n\)-vertex graph with \(n\ge1\), it makes \(O(rn\log(2n))\) calls to \Cref{thm:algorithm} with varying objective weights; see \Cref{subsec:construct-distribution}.

\begin{corollary}\label{cor:algorithm-local}
There is a deterministic algorithm that, given an integer \(r\ge4\) and a \(K_r\)-free graph \(G\) on \(n\) vertices, returns a list \(I_1,\ldots,I_m\) of independent sets satisfying
\[
    \frac{\abs{\{j\in\{1,\ldots,m\}:v\in I_j\}}}{m}
    \ge\frac{c_r}{8}f(d_G(v))\qquad(v\in V(G)).
\]
The list may contain repeated sets.
For \(n\ge1\), the list has length \(m=O(rn\log(2n))\), with an absolute implied constant, and the running time is at most \(2^{\rtconstant r^2}n^{\rtdegree}\).

In particular, for an absolute constant $\kappa_{\rm alg}>0$, choosing an index $j$ uniformly from $\{1,\ldots,m\}$ gives \[ \Prob(v\in I_j)\ge\frac{\kappa_{\rm alg}}{r} \frac{\log(d_G(v)+2)}{d_G(v)+1}\qquad(v\in V(G)). \]
\end{corollary}

\subsection{\texorpdfstring{Comparison with the recent proofs of \Cref{conj:aeks}}{Comparison with the recent proofs of the AEKS conjecture}}
OpenAI~\cite{OAI2026logarithmic} proved \Cref{conj:aeks} by bounding the total weight of triangles at a maximizer of an auxiliary \emph{potential function}, adapted from Davies~\cite[Section~3.3]{davies2026random}.
Kang~\cite{kang2026ode} refined this argument to obtain the bound $\Omega(n\log d/(rd))$ with an absolute constant.
We adapt their method for controlling triangles, but our setting differs from theirs in two key respects.

First, the demands are local.
Rather than a single large independent set, we need, for every nonnegative objective weight $u$, an independent set meeting the target $\sum_v u(v)\phi(v)$.
We therefore extend the weighted induction of Martinsson and Steiner~\cite{martinsson2025random} and Davies~\cite{davies2026random} with a new potential that incorporates the objective weights $u$ and Shearer's function $f$.
The triangle bound we need is correspondingly weighted by $u$, and the linear dependence on $r$ comes from Kang's deterministic neighborhood restriction~\cite{kang2026ode}.

Second, we seek a deterministic algorithm rather than an existence proof.
This requires removing the randomness in OpenAI's argument~\cite{OAI2026logarithmic}, which we do using Kang's neighborhood restriction~\cite{kang2026ode} and the method of conditional expectations.
It also requires avoiding the potentially intractable task of globally maximizing the potential; instead, we make explicit changes that increase the potential.
We discuss both differences in more detail in \Cref{sec:overview}.

\subsection{Notation}\label{subsection: notation}
Throughout the rest of the paper, we use the following notation.
All graphs are finite and simple.
For a graph $G$, we write $V(G)$ and $E(G)$ for its vertex and edge sets, $N_G[v]\defeq N_G(v)\cup\{v\}$, $N_G(A)\defeq\bigcup_{v\in A}N_G(v)$, and $N_G[A]\defeq A\cup N_G(A)$.
The set $N_G[v]$ is the \emph{closed neighborhood} of $v$.
For $A\subseteq V(G)$, we let $G[A]$ denote the subgraph induced by $A$, and we write $G-A\defeq G[V(G)\setminus A]$ and $G_x\defeq G-N_G[x]$.
We omit the graph subscript when the graph is clear from context.

We view vertex weights as vectors indexed by $V(G)$, writing their values as $a(v)$ and writing auxiliary parameters with subscripts, such as $\theta_v$.
We use the notation $a(A)\defeq\sum_{y\in A}a(y)$ for any real-valued function $a$ on a finite set $S$ and any $A\subseteq S$, and for subscripted parameters we write $a(A)\defeq\sum_{y\in A}a_y$.
Inequalities between vertex weights are coordinatewise: $a\le b$ means $a(v)\le b(v)$ for every vertex $v$.
When passing to an induced subgraph, we restrict each vertex weight to its vertex set without changing its name.
Unqualified vertex sums refer to the graph currently under consideration.
We let $\bone$ denote the vertex weight equal to one at every vertex.

For a real vector $z$ indexed by $V(G)$, $\|z\|_1$ and $\|z\|_\infty$ denote the sum and the maximum of the absolute values of its coordinates.
We use $\ind_{\mathcal A}$ for the indicator of an event $\mathcal A$.

\section{Proof overview}\label{sec:overview}
We first explain the argument in its simplest setting: finding one large independent set in a regular graph.
This setting shows where Shearer's function comes from and why triangles are the main obstacle (\Cref{subsec:overview-regular,subsec:overview-triangle-cost}).
We then introduce the objective weights needed for local demands, together with a potential that replaces Shearer's function (\Cref{subsec:overview-induction,subsec:overview-potential}), and compute exactly how the potential changes when a vertex is selected (\Cref{subsec:overview-deletion}).
The remaining obstacle is a bound on the weighted triangle count, which we obtain from the maximality of the potential (\Cref{subsec:overview-increase,subsec:overview-triangle-lemma}).
Finally, we explain how to turn the argument into a deterministic algorithm (\Cref{subsec:overview-algorithm,subsec:overview-derandomization}).

The triangle bound adapts the idea of OpenAI~\cite{OAI2026logarithmic}, refined by Kang~\cite{kang2026ode}, that if the total weight of triangles is too large, the vertex weights can be changed to increase the potential.
Our argument departs from theirs in two respects.
First, OpenAI's potential depends only on auxiliary vertex weights, whereas ours incorporates an independently prescribed objective weight $u$ through Shearer's function; this is what yields the bound in \Cref{thm:existence} for every nonnegative $u$, as local demands require.
Second, the existence proof relies on a global maximizer of the potential, whose computation may be intractable.
We therefore track the effect of selecting a vertex at arbitrary auxiliary weights, make progress through explicit weight updates and vertex selections, and make the procedure deterministic by computing the relevant conditional expectations exactly.

\subsection{Regular triangle-free graphs}\label{subsec:overview-regular}
For Shearer's function $f$ from~\eqref{eq: shearer function} and a graph $H$, define \[ \Psi(H)\defeq\Phi_{H,\bone}(\bone)=\sum_{v\in V(H)}f(d_H(v)). \]
We sketch a proof by induction on $|V(H)|$ that every triangle-free graph $H$ satisfies $\alpha(H)\ge\Psi(H)$; this is Shearer's degree-sequence bound~\cite{shearer1991note}.
Since deleting vertices destroys regularity, the induction hypothesis must apply to all triangle-free graphs.
For brevity, we sketch the induction step only when the current graph $G$ is $d$-regular with $n$ vertices and $d\ge2$, where the calculation is most transparent.
In this case $\Psi(G)=nf(d)$, which is asymptotically $n\log d/d$.
The standard greedy argument accounts only for the deleted vertices and gives $n/(d+1)$, which is tight for disjoint unions of copies of $K_{d+1}$.
For triangle-free graphs, the extra factor of $\log d$ comes from also tracking how the degrees of the surviving vertices decrease; when triangles are present, additional terms arise that must be controlled (\Cref{subsec:overview-triangle-cost}).
Although the sketch is simple, it is a substantial challenge to design a method that generalizes cleanly from regular to irregular graphs, from the unweighted to the weighted setting, and from triangle-free to $K_r$-free graphs.

\begin{enumerate}[label=\arabic*.,ref=\arabic*]
\item\label{step:overview-regular-select} \textbf{Select a vertex.}
For $x\in V(G)$, let $G_x\defeq G-N[x]$.
By induction, $G_x$ has an independent set of size at least $\Psi(G_x)$, and adding $x$ to it shows that $\alpha(G)\ge1+\Psi(G_x)$.
Thus the induction step succeeds as soon as some vertex $x$ has nonnegative \emph{gain}: \begin{equation}\label{eq:overview-regular-gain} 1+\Psi(G_x)-\Psi(G)\ge0. \end{equation}

\item\label{step:overview-regular-deletion} \textbf{Compute the deletion.}
We call the computation of $\Psi(G)-\Psi(G_x)$ the \emph{deletion calculation}.
A vertex $v$ \emph{survives} if $v\in V(G_x)$.
\begin{enumerate}[label=(\alph*),ref=(\alph*)]
\item\label{step:overview-deleted} \emph{Deleted vertices.}
The $d+1$ terms indexed by $N[x]$ are removed, and each equals $f(d)$.
\item\label{step:overview-surviving} \emph{Surviving vertices.}
If a surviving vertex loses $k$ neighbors, its term changes from $f(d)$ to $f(d-k)$.
Since $f$ is convex, \[ f(d-k)-f(d)\ge-kf'(d)\qquad(0\le k\le d). \]
Thus, each lost neighbor increases the term of a surviving vertex by at least $-f'(d)>0$.
\item\label{step:overview-lost-neighbors} \emph{Counting lost neighbors.}
The total number of neighbors lost by surviving vertices equals the number of edges from $N(x)$ to $V(G)\setminus N[x]$.
Since $G$ is triangle-free, each of the $d$ neighbors of $x$ has $d-1$ further neighbors, all outside $N[x]$, so this number is exactly $d(d-1)$.
This is the only place where we use that $G$ is triangle-free.
\end{enumerate}

\item\label{step:overview-regular-combine} \textbf{Combine.}
Subtracting the deleted terms in Step~\ref{step:overview-regular-deletion}\ref{step:overview-deleted} and adding the increases from Step~\ref{step:overview-regular-deletion}\ref{step:overview-surviving}--\ref{step:overview-lost-neighbors} gives \[ 1+\Psi(G_x)-\Psi(G) \ge1-(d+1)f(d)-d(d-1)f'(d)=0, \] where the equality is the differential identity \begin{equation}\label{eq:overview-ode} d(d-1)f'(d)+(d+1)f(d)=1 \end{equation} satisfied by Shearer's function.
Thus~\eqref{eq:overview-regular-gain} holds for every vertex $x$.
\end{enumerate}

In the general argument, which also involves vertex weights, these three steps change in two ways.
First, the (weighted) degrees vary, so in Step~\ref{step:overview-regular-combine} we average the gain over $x$ rather than bounding it for each $x$ separately.
Second, the inequality in Step~\ref{step:overview-regular-deletion}\ref{step:overview-surviving} is replaced by a delicate exact deletion identity involving the vertex weights.

\subsection{The cost of triangles}\label{subsec:overview-triangle-cost}
Keep the assumption that $G$ is $d$-regular, but now allow triangles.
Let $t_x$ be the number of triangles containing $x$, equivalently the number of edges with both endpoints in $N(x)$.
Steps~\ref{step:overview-regular-select} and~\ref{step:overview-regular-deletion}\ref{step:overview-deleted}--\ref{step:overview-surviving} of \Cref{subsec:overview-regular} are unchanged, and only the count in Step~\ref{step:overview-regular-deletion}\ref{step:overview-lost-neighbors} differs.
The quantity $d(d-1)$ counts ordered pairs $(y,z)$ with $y\in N(x)$ and $z\in N(y)\setminus\{x\}$, and it counts each of the $t_x$ edges inside $N(x)$ twice.
Thus the surviving vertices lose only $d(d-1)-2t_x$ neighbors in total, and~\eqref{eq:overview-ode} gives \begin{equation}\label{eq:overview-triangle-cost} 1+\Psi(G_x)-\Psi(G) \ge1-(d+1)f(d)-\bigl(d(d-1)-2t_x\bigr)f'(d) =2t_xf'(d). \end{equation}
Since $f'(d)<0$, this lower bound is negative whenever $x$ lies in a triangle: each edge inside $N(x)$ costs two lost neighbors, each worth at least $-f'(d)$.

Excluding a larger clique does not by itself make this cost small.
For example, the complete tripartite graph $K_{m,m,m}$ is $K_4$-free and $2m$-regular, but every vertex lies in $m^2$ triangles.
Since $f'(d)\sim-\log d/d^2$, the lower bound $2m^2f'(2m)$ tends to $-\infty$ as $m\to\infty$, so the estimate does not guarantee a nonnegative gain.
To offset this cost, we introduce auxiliary vertex weights; maximizing a potential over these weights will bound the weighted triangle count (\Cref{subsec:overview-increase}).

\subsection{The statement proved by induction}\label{subsec:overview-induction}
From now on, fix an integer $r\ge4$ and a $K_r$-free graph $G$.
By the linear programming duality criterion from \Cref{subsection: weighted independent sets intro}, $G$ admits a $\phi$-coloring for a demand function $\phi$ if and only if $\alpha_u(G)\ge\sum_vu(v)\phi(v)$ for every nonnegative objective weight $u$ (see \Cref{lem:demands}).
It suffices to consider positive $u$: for each fixed $s$, both $\alpha_u(G)$ and $\Phi_{G,u}(s)$ are continuous in $u$, so an inequality between them that holds for $u+\varepsilon\bone$ for every $\varepsilon>0$ also holds for $u$.
We therefore fix a positive vertex weight $u$ on $G$ and prove \begin{equation}\label{eq:overview-target} \alpha_u(G)\ge c_r\sup_{s>0}\Phi_{G,u}(s), \end{equation} where $s$ ranges over all positive vertex weights and $c_r>0$ depends only on $r$; we choose $c_r$ at the end of \Cref{subsec:overview-increase}.
In particular, \eqref{eq:overview-target} gives the bound in \Cref{thm:existence} for every prescribed $s$, and taking $s\equiv\bone$ gives the bound in terms of ordinary degrees.

We prove~\eqref{eq:overview-target} by induction on $|V(G)|$, following the deletion argument of \Cref{subsec:overview-regular,subsec:overview-triangle-cost}.
The supremum over $s$ plays two roles.
First, when we apply the induction hypothesis to $G_x$, the weight on $G_x$ is inherited from $G$ and is not under our control, so the hypothesis must hold for every positive weight, just as the hypothesis in \Cref{subsec:overview-regular} had to hold for every degree sequence.
Second, on $G$ itself we may use a maximizing weight, and its optimality controls the triangles of \Cref{subsec:overview-triangle-cost}: an excessive weighted triangle count would yield a better weight.

Both roles are difficult to handle directly with $\Phi_{G,u}$, because $d_{G,s}(v)$ depends on the weights of $v$ and of all its neighbors.
Consequently, deleting $N[x]$ changes weighted degrees up to distance two from $x$, which requires convexity estimates as in Step~\ref{step:overview-regular-deletion}\ref{step:overview-surviving}, and the optimality conditions for $s$ involve entire neighborhoods.
We therefore replace $\Phi_{G,u}$ by a function with the same supremum whose vertex terms depend only on their own weights.
For now, we assume that \(G\) has no isolated vertices and admits a positive maximizing weight; we remove these assumptions at the end of \Cref{subsec:overview-increase}.

\subsection{The potential}\label{subsec:overview-potential}
The idea is to treat the weighted degree of each vertex as a free variable.
Let $w$ be an auxiliary positive vertex weight.
For each vertex $v$, setting $t=d_{G,w}(v)$ gives \[ u(v)f(t)+t\,w(v)^2=u(v)f(d_{G,w}(v))+w(v)w(N(v)). \]
Summing over $v$, the last terms add up to $2\sum_{xy\in E(G)}w(x)w(y)$, since each edge is counted from both endpoints.
Thus $\Phi_{G,u}(w)$ is obtained by evaluating the left side at $t=d_{G,w}(v)$, summing over $v$, and subtracting twice the edge sum.
In this expression, the neighbors of $v$ enter the vertex term only through the particular value $t=d_{G,w}(v)=w(N(v))/w(v)$.
If we instead minimize over all $t>0$, the result $\inf_{t>0}\{u(v)f(t)+t\,w(v)^2\}$ is a function of $u(v)$ and $w(v)$ alone.
The weights of the neighbors then appear only in the edge sum.
Accordingly, for $a,z>0$ we define \[ U_a(z)\defeq\inf_{t>0}\{af(t)+tz^2\},\qquad U_a(0)\defeq0, \] and we define the \emph{potential} \begin{equation}\label{eq:overview-potential} \mathcal U_{G,u}(w)\defeq\sum_vU_{u(v)}(w(v)) -2\sum_{xy\in E(G)}w(x)w(y). \end{equation}
Each vertex term depends only on the weight at that vertex, and the weights interact only through the edge sum, which is bilinear.
Throughout, $u$ and $s$ are fixed by the input, while $w$ is auxiliary: the existence proof takes $w$ to be a maximizer, and the algorithm modifies it step by step.

The infimum defining $U_a(z)$ is attained at the unique $t>0$ satisfying $-af'(t)=z^2$.
Choosing $t=d_{G,w}(v)$ in each vertex term instead shows that $\mathcal U_{G,u}(w)\le\Phi_{G,u}(w)$.
For the reverse comparison, given any positive weight $s$, set \begin{equation}\label{eq:overview-q-def} q(v)\defeq\sqrt{-u(v)f'(d_{G,s}(v))}. \end{equation}
Since $G$ has no isolated vertices, these values are finite and positive, and \Cref{lem:lift} shows that \begin{equation}\label{eq:overview-initialization} \mathcal U_{G,u}(q)\ge\Phi_{G,u}(s). \end{equation}
Hence $\Phi_{G,u}$ and $\mathcal U_{G,u}$ have the same supremum over positive vertex weights.
Thus replacing the entropy terms used by OpenAI~\cite[Section~2]{OAI2026logarithmic} with the scalar terms $U_{u(v)}$ loses nothing in the weighted local target.
If $s$ attains the supremum of $\Phi_{G,u}$, then $q$ attains that of $\mathcal U_{G,u}$, and we write $w\defeq q$ for this maximizer.
Since $\mathcal U_{G,u}$ extends continuously to nonnegative weights, $w$ also maximizes it over nonnegative vertex weights; this will matter in \Cref{subsec:overview-increase}, where we set some weights to zero.
Moreover, each $w(v)$ can be varied in either direction with the other weights fixed, so every partial derivative vanishes: \begin{equation}\label{eq:overview-potential-derivative} \frac{\partial\mathcal U_{G,u}}{\partial w(v)}(w)=U'_{u(v)}(w(v))-2w(N(v))=0 \qquad(v\in V(G)), \end{equation} where $U'_a$ denotes the derivative of $U_a$.
Unlike the optimality conditions for $\Phi_{G,u}$, each of these conditions involves only $w(v)$ and the weight of the neighborhood of $v$.

The induction step now reduces to a statement about the potential.
At the maximizer, $\mathcal U_{G,u}(w)$ equals the supremum in~\eqref{eq:overview-target}, so it suffices to find an independent set of weight at least $c_r\,\mathcal U_{G,u}(w)$.
For each $x$, the induction hypothesis for $G_x$, evaluated at the restriction of $w$ to $V(G_x)$, gives
\[
    \alpha_u(G_x)\ge c_r\,\Phi_{G_x,u}(w)\ge c_r\,\mathcal U_{G_x,u}(w).
\]
Adjoining $x$ to such an independent set, it therefore suffices to find a vertex $x$ with \begin{equation}\label{eq:overview-goal} u(x)+c_r\,\mathcal U_{G_x,u}(w)\ge c_r\,\mathcal U_{G,u}(w). \end{equation}
Selecting $x$ deletes the vertex terms indexed by $N[x]$ and the edge terms meeting $N[x]$, and leaves every other term of $\mathcal U_{G,u}$ unchanged.
The next subsection computes a weighted average of these changes exactly.

\subsection{An exact deletion calculation}\label{subsec:overview-deletion}
We follow the three steps of \Cref{subsec:overview-regular}, with $\Psi$ replaced by the potential and with an arbitrary positive auxiliary weight $w$; the algorithm will need the calculation at weights that are not maximizers.
In the deletion calculation below, the terms indexed by a vertex $v\in N[x]$ combine into $U_{u(v)}(w(v))-2w(v)w(N(v))$.
To understand this net contribution, we split the vertex term using its tangent line at $w(v)$: \[ U_{u(v)}(w(v))-2w(v)w(N(v))=\xi_v+w(v)\zeta_v, \] where \[ \xi_v\defeq U_{u(v)}(w(v))-w(v)U'_{u(v)}(w(v)),\qquad \zeta_v\defeq U'_{u(v)}(w(v))-2w(N(v)). \]
The two parts play different roles.
The number $\xi_v$ is the value at zero of the tangent line to $U_{u(v)}$ at $w(v)$; since $U_{u(v)}$ is strictly concave with $U_{u(v)}(0)=0$, we have $\xi_v>0$.
The number $\zeta_v$ is the partial derivative of the potential with respect to $w(v)$, as in~\eqref{eq:overview-potential-derivative}, so it vanishes at the maximizer of \Cref{subsec:overview-potential}.
Thus, at a maximizer, each deleted vertex contributes exactly $\xi_v$, and in general $\zeta_v$ measures how far $w$ is from satisfying the optimality conditions.

\begin{enumerate}[label=(Del\arabic*),ref=(Del\arabic*)]
\item\label{step:overview-potential-select} \textbf{Select a vertex.}
Define the gain from selecting $x$ as the difference between the two sides of~\eqref{eq:overview-goal}: \[ g_x\defeq u(x)+c_r\,\mathcal{U}_{G_x,u}(w)-c_r\,\mathcal{U}_{G,u}(w). \]
In \Cref{subsec:overview-regular}, we bounded the gain separately for every vertex.
Here we instead bound a weighted average of the gains: we show that $\sum_xw(x)g_x\ge0$, and since every $w(x)$ is positive, some $g_x$ is then nonnegative.
This requires $c_r$ to be sufficiently small in terms of $r$; we determine how small below.

\item\label{step:overview-potential-deletion} \textbf{Compute the deletion.}
The deletion calculation is now the identity \begin{equation}\label{eq:overview-deletion-calculation} \mathcal U_{G,u}(w)-\mathcal U_{G_x,u}(w) =\underbrace{\sum_{v\in N[x]}U_{u(v)}(w(v))}_{\text{\ref{step:overview-potential-deleted}}} -\underbrace{2\sum_{v\in N[x]}w(v)w(N(v))}_{\text{\ref{step:overview-potential-surviving}}} +\underbrace{2\sum_{yz\in E(G[N[x]])}w(y)w(z)}_{\text{\ref{step:overview-potential-correction}}}. \end{equation}
\begin{enumerate}[label=(\alph*),ref=(\alph*)]
\item\label{step:overview-potential-deleted} \emph{Deleted vertices.}
The vertex terms indexed by $N[x]$ are removed, as before.
\item\label{step:overview-potential-surviving} \emph{Deleted edges.}
Unlike in Step~\ref{step:overview-regular-deletion}\ref{step:overview-surviving}, the vertex terms of surviving vertices do not change, since each depends only on the weight of its own vertex.
Instead, the gain from deleting $N[x]$ appears in the edge terms: deleting an edge $yz$ meeting $N[x]$ removes its negative contribution $-2w(y)w(z)$ to the potential.
The sum in~\ref{step:overview-potential-surviving} counts such edges from their endpoints in $N[x]$.
\item\label{step:overview-potential-correction} \emph{Correcting the count.}
An edge with both endpoints in $N[x]$ is counted twice in~\ref{step:overview-potential-surviving}, and the last sum corrects this.
Such an edge either contains $x$ or lies inside $N(x)$; in the latter case it forms a triangle with $x$, as for the edges counted by $t_x$ in \Cref{subsec:overview-triangle-cost}.
\end{enumerate}

\item\label{step:overview-potential-combine} \textbf{Combine.}
Average~\eqref{eq:overview-deletion-calculation} over $x$ with coefficients $w(x)$.
The terms in~\ref{step:overview-potential-deleted} and~\ref{step:overview-potential-surviving} indexed by $v$ appear exactly when $x\in N[v]$, so they receive total coefficient $w(N[v])$, the weighted analogue of $|N[v]|=d+1$.
An edge lies inside $N[x]$ when $x$ is one of its endpoints or a common neighbor of both endpoints.
In the second case, the edge and $x$ form a triangle, so each triangle is counted three times, once for each of its edges, and the factor of two in~\ref{step:overview-potential-correction} brings this to six.
Put \[ S\defeq\sum_vu(v)w(v),\qquad T_G(w)\defeq\sum_{\{x,y,z\}\text{ a triangle}}w(x)w(y)w(z), \] where each triangle is counted once.
Applying the differential identity~\eqref{eq:overview-ode} at the minimizers defining the vertex terms $U_{u(v)}(w(v))$ gives the following exact formula, proved in \Cref{lem:deletion}: \begin{equation}\label{eq:overview-deletion} \sum_xw(x)g_x=(1-c_r)S-6c_rT_G(w) -c_r\sum_v\zeta_v\bigl(w(v)w(N(v))-\xi_v/2\bigr). \end{equation}
This is the weighted analogue of~\eqref{eq:overview-triangle-cost}.
If $c_r=1$ and every $\zeta_v$ vanishes, only the triangle term $-6T_G(w)$ remains, just as only $2t_xf'(d)$ remained there.
Taking $c_r<1$ leaves the positive term $(1-c_r)S$ to absorb the triangle cost.

At a maximizing weight, every $\zeta_v$ is zero, so the last sum vanishes.
It therefore suffices to show that $T_G(w)\le C_rS$ at a maximizing weight, for some $C_r$ depending only on $r$: if $c_r\le(1+6C_r)^{-1}$, then \[ \sum_xw(x)g_x\ge\bigl(1-c_r(1+6C_r)\bigr)S\ge0, \] and the induction step follows.
For the algorithm, the last sum in~\eqref{eq:overview-deletion} is the exact error incurred when the $\zeta_v$ do not vanish.
\end{enumerate}

\subsection{Bounding triangles or increasing the potential}\label{subsec:overview-increase}
By~\eqref{eq:overview-deletion}, the induction step reduces to showing that $T_G(w)\le C_rS$ at a maximizing weight $w$.
We prove this using maximality: if $T_G(w)$ were too large, we would construct a nonnegative vertex weight with a larger potential.
\Cref{lem:triangle-test,lem:improvement} carry this out in the following steps.

\begin{enumerate}[label=(T\arabic*),ref=(T\arabic*)]
\item\label{step:overview-b-reduce} \textbf{Reduce to a dichotomy.}
The triangles containing a vertex $v$ correspond to the edges inside $N(v)$, so $T_G(w)$ is small when every neighborhood has small weight.
To quantify this, fix positive parameters $\theta_v$, to be chosen in~\ref{step:overview-b-change}, and think of $\theta_v$ as an allowed size for $w(N(v))$.
Put \[ D_\theta\defeq\sum_v\theta_vw(v)w(N(v)) =\sum_{xy\in E(G)}w(x)w(y)(\theta_x+\theta_y). \]
If \(w(N(v))\le R\theta_v\) for every \(v\) and some \(R>0\), then the edges inside \(N(v)\), counted in both directions, have total weight at most \(w(N(v))^2\le R\theta_vw(N(v))\).
Multiplying by \(w(v)\) and summing over \(v\) counts each triangle six times, so \(6T_G(w)\le RD_\theta\).

Neighborhoods need not be this small, but \Cref{lem:triangle-test} shows that a large triangle weight forces a dense pair of vertex sets.
For vertex sets $A,B$, write \begin{equation}\label{eq:overview-edge-weight} e_w(A,B)\defeq\sum_{x\in A}\sum_{y\in B\cap N(x)}w(x)w(y), \end{equation} where an edge with both endpoints in $A\cap B$ is counted in both directions.
\Cref{lem:triangle-test}, outlined in \Cref{subsec:overview-triangle-lemma}, applies to arbitrary positive $\theta_v$ and gives one of the following alternatives:
\begin{enumerate}[label=(\roman*),ref=(\roman*)]
\item\label{case:overview-triangle-bound} $T_G(w)\le C_\triangle D_\theta$ for an absolute constant $C_\triangle\ge1$;
\item\label{case:overview-triangle-witness} there are vertex sets $A,B$ such that
\begin{itemize}
\item the edge weight $e_w(A,B)$ is large compared with $\sum_{v\in A}\theta_vw(v)$, and
\item $w(A)$, $w(B)$, and $\sum_{v\in B}\theta_vw(v)$ are bounded above.
\end{itemize}
\end{enumerate}
\Cref{lem:triangle-test} gives the precise inequalities.
In~\ref{step:overview-b-heavy}--\ref{step:overview-b-change}, we show that alternative~\ref{case:overview-triangle-witness} contradicts the maximality of $w$, so alternative~\ref{case:overview-triangle-bound} must hold.

\item\label{step:overview-b-heavy} \textbf{Find heavy neighborhoods.}
Suppose that alternative~\ref{case:overview-triangle-witness} holds.
\Cref{lem:improvement} finds disjoint subsets $A'\subseteq A$ and $B'\subseteq B$, and a nonempty subset $A''\subseteq A'$, such that \[ w(N(v)\cap B')>12p\theta_v+\delta \qquad(v\in A''), \] where $p\ge1$ is an integer supplied by \Cref{lem:triangle-test} and $\delta>0$.
Thus every vertex of $A''$ has a heavy neighborhood in $B'$, which will produce an edge gain when we change the weights in~\ref{step:overview-b-change}.

\item\label{step:overview-b-sparse} \textbf{Choose a set with few internal edges.}
Increasing the weights on $A''$ also increases the cost of the edges inside $A''$, so we first pass to a subset with few internal edges.
For a parameter $Q\ge1$, start with $C\defeq A''$.
While some $v\in C$ satisfies $w(N(v)\cap C)>w(C)/Q$, replace $C$ by $N(v)\cap C$.
The vertices chosen in successive restrictions form a clique, so after $r-2$ restrictions the current set is independent, since an edge in it would complete a copy of $K_r$.
Each restriction retains more than a $1/Q$ fraction of the weight, so the procedure stops after at most $r-2$ restrictions with a nonempty set $C$ satisfying \[ w(C)\ge Q^{-(r-2)}w(A''),\qquad w(N(v)\cap C)\le w(C)/Q\quad(v\in C). \]
This restriction, used in the proof of~\cite[Lemma~3.1]{kang2026ode} and stated as \Cref{lem:sparse-restriction}, is the only step that uses the assumption that $G$ is $K_r$-free.
In OpenAI~\cite[Section~2]{OAI2026logarithmic}, random rescaling controls the cost of the edges inside the set whose weights are increased; Kang's restriction does this deterministically.

\item\label{step:overview-b-change} \textbf{Change the weights and choose $\theta$.}
We obtain $w'$ from $w$ by multiplying the weights on $C$ by $1+\alpha$, setting the weights on $B'$ to zero, and leaving all other weights unchanged.
Since $w'$ is nonnegative and $w$ maximizes the potential over nonnegative weights, $\mathcal U_{G,u}(w')\le\mathcal U_{G,u}(w)$.
We contradict this inequality.
By~\eqref{eq:overview-potential-derivative}, the first-order terms cancel, so the gain \(2\alpha\,e_w(C,B')\) from the edges between \(C\) and \(B'\) must cover the cost of the edges inside \(C\) and the vertex losses due to concavity.
For our potential, these vertex losses are the errors in replacing each \(U_{u(v)}\) by its tangent line at \(w(v)\), and we bound them for arbitrary objective weights \(u\).

Since $C\subseteq A''$, the edge gain contains $24\alpha p\sum_{v\in C}\theta_vw(v)$, while the vertex losses on $C$ are at most $48r\alpha p\sum_{v\in C}\xi_v$.
We therefore choose \[ \theta_v\defeq \frac{2r\xi_v}{w(v)}, \] which makes these two quantities equal, so the edge gain covers the vertex losses on $C$.
The factor $r$ comes from the at most $r-2$ restrictions in~\ref{step:overview-b-sparse}.
The same choice of $\theta$ controls the vertex losses on $B'$, and the restriction in~\ref{step:overview-b-sparse} controls the cost of the edges inside $C$.
The proof of \Cref{lem:improvement} chooses $\alpha$ and $Q$ so that the remaining edge gain exceeds all remaining costs.
Thus $\mathcal U_{G,u}(w')>\mathcal U_{G,u}(w)$, contradicting maximality.

\item\label{step:overview-b-conclude} \textbf{Conclude.}
At a positive maximizer, alternative~\ref{case:overview-triangle-witness} is therefore impossible, so $T_G(w)\le C_\triangle D_\theta$.
The choice of $\theta$ also bounds $D_\theta$ in terms of $S$: by the optimality conditions at $w$, \[ D_\theta=2r\sum_v\xi_vw(N(v))\le2rS. \]
Thus $T_G(w)\le2C_\triangle rS$, which is the required triangle bound with $C_r\defeq2C_\triangle r$.
By \Cref{subsec:overview-deletion}, the induction step holds whenever $c_r\le(1+12C_\triangle r)^{-1}$.
The smaller value of $c_r$ specified in~\eqref{eq:constant} also absorbs the approximation errors in the algorithm, and it remains bounded above and below by absolute positive constants divided by $r$.
\end{enumerate}

The preceding argument assumes that \(G\) has no isolated vertices and that a positive weight attains the supremum.
If \(v\) is isolated, then \(\alpha_u(G)=u(v)+\alpha_u(G-v)\), and \(v\) contributes exactly \(u(v)\) to \(\Phi_{G,u}(s)\) for every \(s\); since \(c_r\le1\), \eqref{eq:overview-target} for \(G\) follows from \eqref{eq:overview-target} for \(G-v\) (see \Cref{sec:existence}).
If no positive weight attains the supremum, a standard limiting argument reduces the problem to smaller induced subgraphs with no edges between them, whose independent sets can be combined by induction; see \Cref{subsec:boundary-limit}.
This completes the proof of~\eqref{eq:overview-target}, and \Cref{lem:demands} then gives \Cref{cor:local}.

\subsection{Triangle control by vertex splitting}\label{subsec:overview-triangle-lemma}
We now outline the proof of \Cref{lem:triangle-test}, which gives either the triangle bound in alternative~\ref{case:overview-triangle-bound} or the sets in alternative~\ref{case:overview-triangle-witness}; the label SD stands for \emph{structural dichotomy}.
The idea is to modify $G$ by deleting edges and splitting vertices until every neighborhood is small relative to $\theta$, so that the counting argument in~\ref{step:overview-b-reduce} applies, while showing that these modifications destroy only a controlled amount of triangle weight.
Whenever the modifications cannot proceed with a controlled loss, the proof finds the sets in alternative~\ref{case:overview-triangle-witness}.
This argument is adapted from OpenAI~\cite[Sections~3--5]{OAI2026logarithmic}, with a quantitative refinement due to Kang~\cite[Sections~4--5]{kang2026ode}.
It applies to arbitrary positive $w$ and $\theta$, and it uses neither clique exclusion nor maximality.

\begin{enumerate}[label=(SD\arabic*),ref=(SD\arabic*)]
\item\label{step:overview-triangle-split} \textbf{Split vertices.}
Splitting a vertex means partitioning its incident edges into groups and replacing the vertex by one copy for each group.
Every copy inherits the original weight $w(v)$ and parameter $\theta_v$, and each original edge is assigned to one copy at each of its endpoints.
Since $D_\theta$ is a sum over edges (see~\ref{step:overview-b-reduce}), and each edge keeps the weights and parameters of its endpoints, splitting preserves $D_\theta$, while deleting edges can only decrease it.
Splitting does, however, destroy triangles: a triangle survives only if none of its edges is deleted and, at each of its vertices, its two edges are assigned to the same copy.

\item\label{step:overview-triangle-loss} \textbf{Control triangle losses.}
To split a vertex $x$, we assign each neighbor $v\in N(x)$ a label, and neighbors with the same label are assigned to the same copy of $x$.
A triangle through $x$ survives this split exactly when its other two vertices receive the same label, so we need labels that agree for most pairs of adjacent neighbors.
We take the label of $v$ to be drawn from the endpoint distribution of a random walk started at $v$.
The labels are chosen jointly using the shared rejection sampler of Kleinberg and Tardos~\cite[Section~3]{kleinberg2002approximation}, for which Angel and Spinka~\cite[Section~2.1]{angel2021pairwise} bound the probability that two labels differ in terms of the distance between the corresponding distributions.
These distances are compared using entropy increments, an idea with precedents in work of Erschler and Karlsson~\cite[Lemma~5.1]{erschler2010homomorphisms} and of Benjamini, Duminil-Copin, Kozma, and Yadin~\cite[Section~2, Theorem~8]{benjamini2015disorder}.
Further edge deletions then bound the neighborhood weight of each copy.
Either the expected triangle loss from these deletions and splits is controlled, or the proof finds the sets in alternative~\ref{case:overview-triangle-witness}.
In the existence proof the labels may be random; \Cref{subsec:overview-derandomization} explains how to choose them deterministically.

\item\label{step:overview-triangle-stop} \textbf{Stop and count.}
Write $D_{\rm initial},T_{\rm initial}$ and $D_{\rm final},T_{\rm final}$ for the values of $D_\theta$ and of the triangle weight before and after the modifications.
If alternative~\ref{case:overview-triangle-witness} does not occur, the modifications ensure \[ D_{\rm final}\le D_{\rm initial},\qquad T_{\rm final}\ge\tfrac{1}{2}T_{\rm initial}-D_{\rm initial}. \]
We stop when either \(T_{\rm final}\le2D_{\rm final}\) or every remaining vertex satisfies \(w(N(v))\le \ratiobound{*}\theta_v\) for an absolute constant \(\ratiobound{*}\).
In the latter case, the counting argument in~\ref{step:overview-b-reduce} gives \(6T_{\rm final}\le \ratiobound{*}D_{\rm final}\).
Either way, \(T_{\rm final}\le\max\{2,\ratiobound{*}/6\}\,D_{\rm final}\), and combining this with the two displayed inequalities gives
\[
    T_{\rm initial}\le2(T_{\rm final}+D_{\rm initial})
    \le(2\max\{2,\ratiobound{*}/6\}+2)D_{\rm initial}
    \le C_\triangle D_{\rm initial},
\]
for a suitable absolute constant \(C_\triangle\), chosen in \Cref{sec:triangle-walks}.
\end{enumerate}

\subsection{Finding a weighted independent set}\label{subsec:overview-algorithm}
We now turn from the existence proof to the algorithm of \Cref{thm:algorithm}.
The existence proof uses a global maximizer of the potential twice: to make every $\zeta_v$ vanish in~\eqref{eq:overview-deletion}, and to rule out alternative~\ref{case:overview-triangle-witness} in \Cref{subsec:overview-increase}.
To obtain an efficient construction, we replace this global optimization step by explicit changes to the auxiliary weights, each of which makes measurable progress.

The input consists of an integer \(r\ge4\), a \(K_r\)-free graph \(G\), rational objective weights \(u\ge0\), and rational degree weights \(s>0\).
After handling trivial cases, we normalize so that $\sum_vu(v)=1$.
Throughout, $I$ denotes the independent set selected so far and $H$ the remaining induced subgraph, with no edges between $I$ and $V(H)$.
We track the quantity \[ P\defeq u(I)+c_r\,\mathcal{U}_{H,u}(w). \]
Since $\mathcal U_{H,u}(w)\le\Phi_{H,u}(w)\le u(V(H))$ because $f\le1$, and since $c_r\le1$, we always have $P\le u(I)+u(V(H))\le1$.
We initialize $w$ as a rational approximation of the weight $q$ from~\eqref{eq:overview-q-def}, so that by~\eqref{eq:overview-initialization} the initial value of $P$ is at least $c_r\Phi_{G,u}(s)$, up to an approximation error.
All quantities below are computed in the current graph $H$ with the current weights $w$.
Each step does one of the following.

\begin{enumerate}[label=(Alg\arabic*),ref=(Alg\arabic*)]
\item\label{step:overview-c-range} \textbf{Control extreme weights.}
Discard vertices with negligible objective weight or auxiliary weight; the total loss in $P$ from these discards is small.
Select isolated vertices and vertices with sufficiently large auxiliary weight, adding each to $I$ and deleting its closed neighborhood from $H$.
These selections do not decrease $P$, and afterward all remaining weights lie in a controlled range.

\item\label{step:overview-c-derivative} \textbf{Improve one weight.}
If some partial derivative $\zeta_v$ is sufficiently far from zero, move $w(v)$ in the direction of the derivative: increase it if $\zeta_v>0$, and decrease it otherwise.
We take a step that increases the potential by a definite amount, stopping early if $w(v)$ reaches a threshold handled by~\ref{step:overview-c-range}.

\item\label{step:overview-c-select-or-change} \textbf{Select or reweight.}
Otherwise, $\abs{\zeta_v}\le\eta$ for every vertex $v$, where the tolerance $\eta$ is chosen so that both alternatives below increase $P$.
We choose rational upper approximations to the parameters $\theta_v$ from~\ref{step:overview-b-change} and apply the algorithmic version of the structural dichotomy, \Cref{lem:triangle-test-algorithm}, in $H$.

In alternative~\ref{case:overview-triangle-bound}, we have $T_H(w)\le C_\triangle D_\theta$.
Substituting this bound into the averaged deletion identity~\eqref{eq:overview-deletion}, the triangle term is controlled as in \Cref{subsec:overview-increase}.
Since every $\abs{\zeta_v}\le\eta$, the error term involving $\zeta_v$ cannot cancel the positive contribution, so $\sum_xw(x)g_x>0$.
Evaluating the gains to sufficient accuracy, we find a vertex $x$ with $g_x>0$; adding $x$ to $I$ and deleting $N[x]$ from $H$ increases $P$ by exactly $g_x$.

In alternative~\ref{case:overview-triangle-witness}, we obtain sets $A,B$ as in~\ref{step:overview-b-reduce}, and \Cref{lem:improvement} uses them to construct $w'$ by increasing the weights on a set $C$ and setting the weights on a disjoint set $B'$ to zero.
For our choice of $\eta$, the argument of~\ref{step:overview-b-change} shows that this reweighting increases the potential, and hence $P$.
Since vertices of weight zero contribute nothing to the potential, we may then delete $B'$ from $H$ without changing $P$.
\end{enumerate}

Every step either removes a vertex from \(H\) or increases \(P\) by a definite positive amount, even after rational rounding, and the initialization and discards incur only small losses.
The witness bound \(p=O((1+\log\log(2n))^3)\), together with uniform bounds on the numerical computations and the deterministic structural dichotomy, gives running time \(2^{\rtconstant r^2}L^{\rtdegree}\); see \Cref{subsec:numerics-complexity}.
When $H$ becomes empty, $P=u(I)$, so $u(I)$ is at least $c_r\Phi_{G,u}(s)$ up to an additive error.
Finally, we return the heavier of $I$ and the singleton consisting of a vertex of largest objective weight; this absorbs the additive error at the cost of the factor~$2$ in \Cref{thm:algorithm}.
See \Cref{sec:outer,app:numerics} for the thresholds and error bounds.

\Cref{cor:algorithm-local} follows from a multiplicative-weights construction in the fractional covering framework of Arora, Hazan, and Kale~\cite[Section~3.3.3]{arora2012multiplicative}.
We use their update rule but prove the guarantee directly from \Cref{thm:algorithm}, rather than invoking their covering theorem.
In each round, we call \Cref{thm:algorithm} with \(s\equiv\bone\) and the current positive rational objective weights, then update these weights according to the returned independent set.
After \(O(rn\log(2n))\) rounds, every vertex \(v\) appears with frequency at least \((c_r/8)f(d_G(v))\).
The input lengths in these calls are polynomial in \(n,r\), with absolute degree, giving total running time \(2^{\rtconstant r^2}n^{\rtdegree}\); see \Cref{subsec:construct-distribution,subsec:numerics-complexity}.

\subsection{Choosing the splits deterministically}\label{subsec:overview-derandomization}
The only random choices in the algorithm are the labels in~\ref{step:overview-triangle-loss}, which are drawn from walk-endpoint distributions using the shared rejection sampler.
OpenAI~\cite[Lemma~5.3]{OAI2026logarithmic} bounds the expected triangle weight lost in such a random split and deduces the existence of labels losing no more than this bound.
We instead fix the random choices one at a time by the method of conditional expectations, as Kleinberg and Tardos do for their labeling algorithm~\cite[Section~5]{kleinberg2002approximation}, each time choosing an outcome that does not decrease the conditional expectation of the surviving triangle weight.
The resulting split therefore retains at least the expected triangle weight of the random construction, which is what the inequalities in~\ref{step:overview-triangle-stop} require.
The key point is that these conditional expectations can be computed exactly: with all previous choices fixed, the expected surviving triangle weight is the sum over triangles of their weights times the probability that the two labels agree at each corner and that all three edges are kept.
\Cref{sec:derandomization} computes these probabilities in polynomial time and shows that polynomially many decisions with polynomially many candidate outcomes suffice.
\subsubsection*{Organization}
\Cref{sec:prelim} proves the weighted independent-set criterion and deduces \Cref{cor:local,thm:main} from \Cref{thm:existence}.
\Cref{sec:potential} develops the potential and the exact deletion identity.
\Cref{sec:triangle} states the structural dichotomy and uses potential increases to bound the triangle weight at a positive maximizer.
\Cref{sec:existence,sec:outer} prove the existence and algorithmic theorems, respectively, assuming the structural dichotomy.
\Cref{sec:triangle-walks} develops the neighborhood walks and entropy bounds, and \Cref{sec:triangle-proof} uses them to prove the structural dichotomy by vertex splitting.
\Cref{sec:derandomization} makes the construction deterministic using conditional expectations.
\Cref{sec:conclusion} contains concluding remarks.
The appendices establish the properties of \(f\) and \(U_a\), the bound on total variation distance in terms of relative entropy, and the rational implementation and uniform running time bounds used by the algorithms.


\section{Weighted independent sets and local demands}\label{sec:prelim}
We first prove a standard consequence of linear programming duality: weighted independent-set inequalities yield a single distribution meeting all local demands.
We give a direct geometric proof using the independence polytope.
\begin{lemma}\label{lem:demands}
Let $G$ be a finite graph and let $\phi\colon V(G)\to[0,1]$.
There is a $\phi$-coloring if and only if \[ \alpha_u(G)\ge\sum_vu(v)\phi(v) \quad\text{for every }u\colon V(G)\to[0,\infty). \]
It suffices to check this inequality for positive $u$.
\end{lemma}
\begin{proof}
For a random independent set \(I\) giving a \(\phi\)-coloring, the ``only if'' part follows from \(\alpha_u(G)\ge\E u(I)\ge\sum_vu(v)\phi(v)\).

Conversely, let \(\mathcal{D}\defeq\operatorname{conv}\{\ind_I:I\in\cI(G)\}\) be the \emph{independence polytope} of \(G\).
Its points are precisely the vertex inclusion probabilities obtained from distributions on \(\cI(G)\).
It is compact and convex, and deleting vertices from the sampled sets shows that it is downward closed.
Thus a \(\phi\)-coloring exists exactly when \(\phi\in\mathcal{D}\).

If \(\phi\notin\mathcal{D}\), choose a point \(p\in\mathcal{D}\) nearest to \(\phi\) in Euclidean distance.
Downward closedness gives \(p\le\phi\), since reducing any coordinate above \(\phi\) would decrease the distance.
Set \(u\defeq\phi-p\ge0\), which is nonzero.
For any \(z\in\mathcal{D}\), convexity and minimality of \(p\) imply that the squared distance from \(p+\varepsilon(z-p)\) to \(\phi\) has a nonnegative right derivative at \(\varepsilon=0\).
Hence \(\sum_vu(v)(z(v)-p(v))\le0\), and therefore \[ \alpha_u(G)=\max_{z\in\mathcal{D}}\sum_vu(v)z(v) \le\sum_vu(v)p(v) <\sum_vu(v)\phi(v), \] contrary to the assumed inequality.
Finally, continuity shows that positive weights suffice.
\end{proof}

We now deduce \Cref{cor:local,thm:main} from \Cref{thm:existence}.

\begin{proof}[Proof of \Cref{cor:local,thm:main}]
Since \(0<f\le1\) and \(c_r\le1\), the values \(c_rf(d_{G,s}(v))\) form a demand function.
\Cref{cor:local} follows from \Cref{lem:demands,thm:existence}.
To deduce \Cref{thm:main}, note that the continuous positive function \((d+1)f(d)/\log(d+2)\) tends to one at infinity, so its infimum \(\eta_0\) on \([0,\infty)\) is positive.
Taking \(s\equiv1\) and \(\kappa\defeq\kappa_1\eta_0\) gives the required bound, and the same distribution meets every smaller demand.
\end{proof}


\section{The potential and the deletion identity}\label{sec:potential}
In this section, we construct a potential \(\mathcal U\) with the same supremum as \(\Phi\) and derive an exact identity for the average gain from selecting a vertex and deleting its closed neighborhood (\Cref{lem:deletion}).
This identity reduces the induction step to bounding weighted triangles and, for the algorithm, an error term involving the partial derivatives of \(\mathcal U\).
Throughout this section we assume \(u>0\); vertices of zero objective weight are removed at the start of both the existence proof and the algorithm.

\subsection{The potential}
For auxiliary weights \(w>0\), write \[ \Lambda_v\defeq w(N(v)),\quad M_G(w)\defeq\sum_{xy\in E(G)}w(x)w(y), \] where each edge is counted once.
As in \Cref{subsec:overview-potential}, the identity \[ \Phi_{G,u}(w) =\sum_v\bigl(u(v)f(d_{G,w}(v))+d_{G,w}(v)w(v)^2\bigr)-2M_G(w) \] suggests minimizing each vertex term over its degree parameter.

For \(a,x>0\), define \[ U_a(x)\defeq\inf_{t>0}\{af(t)+tx^2\},\qquad U_a(0)\defeq0, \] and put \[ \mathcal U_{G,u}(w)\defeq\sum_v U_{u(v)}(w(v))-2M_G(w). \]
Each vertex term depends only on the weight at that vertex.
The dependence between different vertex weights occurs in the edge sum.
The extension \(U_a(0)=0\) is continuous, so \(\mathcal U\) is also defined on nonnegative vertex weights.
This lets us represent the deletion of a vertex by setting its weight to zero.
Unlike \(\Phi\), the potential \(\mathcal U\) need not be nonnegative and is not invariant under scaling \(w\).

We suppress the subscripts \(G,u\) on \(\mathcal U_{G,u}\) and \(\Phi_{G,u}\) when both are fixed, and sometimes suppress only \(u\) when comparing graphs.
For example, \(\mathcal U_{G-A}(w)\) means \(\mathcal U_{G-A,u|_{V\setminus A}}(w|_{V\setminus A})\).

The scalar facts used below are proved in \Cref{app:scalar}.
The function \(U_a\) is increasing and concave, with \(0<U_a(x)\le a\) for \(x>0\).
The infimum defining \(U_a(x)\) has a unique minimizer \(d_*(a,x)>0\), characterized by \begin{equation}\label{eq:scalar} -af'(d_*)=x^2,\qquad U_a'(x)=2d_*x. \end{equation}

The definition of \(U_a\) gives the comparison \(\mathcal U_{G,u}(w)\le\Phi_{G,u}(w)\), because an infimum is at most any value of the expression being minimized.
At a nonisolated vertex \(v\), choose \(t=d_{G,w}(v)=\Lambda_v/w(v)>0\) in the definition of \(U_{u(v)}(w(v))\) to obtain \[ U_{u(v)}(w(v)) \le u(v)f(d_{G,w}(v))+w(v)\Lambda_v. \]
For an isolated vertex, the same inequality follows by letting $t$ tend to zero from above in the defining expression, since $f(t)\to f(0)=1$ and $\Lambda_v=0$.
Since \(\sum_vw(v)\Lambda_v=2M_G(w)\), summing gives, for every positive \(w\), \begin{equation}\label{eq:potential-comparison} \mathcal U_{G,u}(w)\le\Phi_{G,u}(w). \end{equation}
This comparison holds for arbitrary graphs, including graphs with isolated vertices.
The next lemma gives an inequality in the opposite direction after changing the vertex weight.

\begin{lemma}\label{lem:lift}
Let \(G\) have no isolated vertices, let \(s\) be a positive vertex weight, and set \(q(v)\defeq\sqrt{-u(v)f'(d_{G,s}(v))}\) at every vertex.
Then \[ \mathcal U_{G,u}(q)\ge\Phi_{G,u}(s). \]
\end{lemma}
\begin{proof}
Every weighted degree is positive, so each \(q(v)\) is finite and positive.
By \eqref{eq:scalar}, the minimizer defining \(U_{u(v)}(q(v))\) is exactly \(d_{G,s}(v)\).
Consequently, grouping the vertex sum by edges gives \begin{align*} \mathcal U_{G,u}(q)-\Phi_{G,u}(s) &=\sum_v d_{G,s}(v)q(v)^2-2M_G(q)\\
    &=\sum_{xy\in E(G)} \left(\sqrt{\frac{s(y)}{s(x)}}q(x) -\sqrt{\frac{s(x)}{s(y)}}q(y)\right)^2\ge0. \qedhere \end{align*}
\end{proof}

The formula for \(q\) is not defined at isolated vertices, since \(-f'(d)\) tends to infinity as \(d\downarrow0\).
Such vertices are handled separately in the existence proof and selected immediately by the algorithm.

\begin{corollary}\label{cor:transfer}
Let \(G\) have no isolated vertices.
We have \[ \sup_{w>0}\mathcal U_{G,u}(w)=\sup_{s>0}\Phi_{G,u}(s). \]
If \(\Phi_{G,u}\) attains its supremum at a positive vertex weight, then \(\mathcal U_{G,u}\) has a positive global maximizer \(w\).
This vertex weight also maximizes \(\mathcal U_{G,u}\) over nonnegative vertex weights and satisfies \begin{equation}\label{eq:transfer-optimality} U'_{u(v)}(w(v))=2w(N(v))\qquad(v\in V(G)). \end{equation}
\end{corollary}
\begin{proof}
The inequality \eqref{eq:potential-comparison} and \Cref{lem:lift} give the equality of suprema.
If \(s\) maximizes \(\Phi\), the positive vertex weight \(q\) from \Cref{lem:lift} satisfies \[ \mathcal U(q)\ge\Phi(s)=\sup_{w>0}\mathcal U(w), \] so \(q\) maximizes \(\mathcal U\).
By continuity, \(q\) also maximizes \(\mathcal U\) over nonnegative vertex weights.
Since \(q>0\), the first-order optimality conditions give \eqref{eq:transfer-optimality} with \(w=q\).
\end{proof}

\subsection{The deletion identity}
Here \(G\) is an arbitrary graph, and we retain \(u>0\).
For \(a,x>0\), put \[ \xi_a(x)\defeq U_a(x)-xU_a'(x). \]
This is the tangent loss when \(x\) is set to zero, since \(U_a(0)=0\).
As proved in \Cref{app:scalar}, \(0<\xi_a(x)\le U_a(x)\le a\) and \begin{equation}\label{eq:scalar-ode} d_*\xi_a(x)\le a,\qquad (1+d_*)\xi_a(x)+2d_*x^2=a, \end{equation} where \(d_*\defeq d_*(a,x)\).
For the current weights put \[ S\defeq\sum_vu(v)w(v),\qquad T_G(w)\defeq\sum_{\{x,y,z\}\in\mathcal T(G)}w(x)w(y)w(z), \] where \(\mathcal T(G)\) is the set of unordered triangles, so each triangle is counted once, and put \[ \xi_v\defeq\xi_{u(v)}(w(v)),\qquad \zeta_v\defeq U_{u(v)}'(w(v))-2\Lambda_v. \]
The number $\zeta_v$ is the partial derivative of $\mathcal U$ with respect to $w(v)$.
Thus $\zeta_v=0$ means $U_{u(v)}'(w(v))=2\Lambda_v$.
Fix \(0<c\le1\).
Selecting a vertex \(x\) adds \(u(x)\) to the weight of the independent set and leaves the graph \(G_x\defeq G-N[x]\).
Define the gain from this selection by \[ g_x\defeq u(x)+c\,\mathcal{U}_{G_x,u}(w)-c\,\mathcal{U}_{G,u}(w). \]

In the existence proof, suppose that \(\Phi\) attains its supremum and take the maximizer \(w\) supplied by \Cref{cor:transfer}.
As in \Cref{subsec:overview-potential}, the inductive bound on \(G_x\), together with \eqref{eq:potential-comparison}, gives an independent set there of weight at least \(c\,\mathcal U_{G_x,u}(w)\).
Adjoining \(x\) therefore reaches the target \(c\,\mathcal U_{G,u}(w)\) whenever \(g_x\ge0\).
The following identity expresses the weighted sum of these gains for arbitrary \(w>0\).

\begin{lemma}\label{lem:deletion}
For every positive vertex weight \(w\), without any optimality assumption, \[ \sum_xw(x)g_x =(1-c)S-6cT_G(w) -c\sum_v \zeta_v\bigl(w(v)\Lambda_v-\xi_v/2\bigr). \]
\end{lemma}
\begin{proof}
We follow the outline in \Cref{subsec:overview-deletion}.
Deleting a set \(A\) removes its vertex terms and all edge terms meeting \(A\).
The sum \(\sum_{v\in A}w(v)\Lambda_v\) counts internal edges twice, so \begin{equation}\label{eq:deletion-cost} \begin{aligned} \mathcal U_G(w)-\mathcal U_{G-A}(w) &=\sum_{v\in A}U_{u(v)}(w(v))-2\sum_{v\in A}w(v)\Lambda_v+2M_{G[A]}(w)\\
    &=\sum_{v\in A}(\xi_v+w(v)\zeta_v)+2M_{G[A]}(w). \end{aligned} \end{equation}
The second equality uses \(U_{u(v)}(w(v))-2w(v)\Lambda_v=\xi_v+w(v)\zeta_v\).

We average this identity with \(A=N[x]\) and coefficients \(w(x)\).
A vertex \(v\) belongs to \(N[x]\) exactly when \(x\in N[v]\), so its coefficient is \(w(v)+\Lambda_v\).
An edge \(vy\) lies in \(G[N[x]]\) when \(x\) is one of its endpoints or a common neighbor.
Hence \begin{align} 2\sum_xw(x)M_{G[N[x]]}(w) &=2\sum_{vy\in E(G)}w(v)w(y)\bigl(w(v)+w(y)\bigr) +2\sum_{vy\in E(G)}w(v)w(y)w\bigl(N(v)\cap N(y)\bigr)\notag\\
    &=2\sum_vw(v)^2\Lambda_v +6\sum_{\{x,y,z\}\in\mathcal T(G)}w(x)w(y)w(z)\notag\\
    &=2\sum_vw(v)^2\Lambda_v+6T_G(w).\label{eq:deletion-edge-average} \end{align}

We now combine the averaged first sum in \eqref{eq:deletion-cost}, whose contribution at \(v\) is \((w(v)+\Lambda_v)(\xi_v+w(v)\zeta_v)\), with the endpoint contribution \(2w(v)^2\Lambda_v\) from \eqref{eq:deletion-edge-average}.
The scalar identity will reduce their sum to \(u(v)w(v)\), apart from terms involving \(\zeta_v\).
Let \(d_{*,v}\defeq d_*(u(v),w(v))\).
By \eqref{eq:scalar} and the definition of \(\zeta_v\), \[ \zeta_v=U'_{u(v)}(w(v))-2\Lambda_v =2d_{*,v}w(v)-2\Lambda_v, \] so \(w(v)d_{*,v}=\Lambda_v+\zeta_v/2\).
Substituting this into the identity in \eqref{eq:scalar-ode}, multiplied by \(w(v)\), gives \[ (w(v)+\Lambda_v)\xi_v+2w(v)^2\Lambda_v =u(v)w(v)-\zeta_v\bigl(w(v)^2+\xi_v/2\bigr). \]
Adding the remaining term \((w(v)+\Lambda_v)w(v)\zeta_v\) leaves \(u(v)w(v)+\zeta_v(w(v)\Lambda_v-\xi_v/2)\).
Thus the averaged deletion cost is \[ S+6T_G(w)+\sum_v\zeta_v(w(v)\Lambda_v-\xi_v/2). \]
The term \(6T_G(w)\) is the extra cost from triangles, and the last sum vanishes when all partial derivatives of \(\mathcal U\) are zero.
Subtracting \(c\) times this cost from \(\sum_xw(x)u(x)=S\) proves the lemma.
\end{proof}

When $\zeta_v=0$ for every vertex, the identity becomes $\sum_xw(x)g_x=(1-c)S-6cT_G(w)$.
In particular, if $G$ is nonempty and $6cT_G(w)\le(1-c)S$, then some $g_x\ge0$, since every $w(x)$ is positive.
Thus a bound $T_G(w)\le C_rS$ at a maximizer suffices when $c\le(1+6C_r)^{-1}$, as in \Cref{subsec:overview-deletion}.
For the algorithm, we instead assume $|\zeta_v|\le\eta$ at every vertex, for a tolerance $\eta>0$, and bound the absolute value of the last sum.


\section{Bounding triangles by increasing the potential}\label{sec:triangle}
Our goal is to bound the triangle weight at a positive maximizer of \(\mathcal U_{G,u}\), as outlined in \Cref{subsec:overview-increase}.
We use a structural dichotomy: either the triangle weight is already bounded, or there are two vertex sets with sufficiently large edge weight between them.
We state this dichotomy in \Cref{lem:triangle-test}; its proof, outlined in \Cref{subsec:overview-triangle-lemma}, is deferred to \Cref{sec:triangle-walks,sec:triangle-proof}, with the deterministic implementation in \Cref{sec:derandomization}.
The main argument of this section uses clique exclusion to turn the second alternative into a potential increase (\Cref{lem:improvement}).
At a maximizer, this increase is impossible, giving the triangle bound in \Cref{cor:stationary}.
We retain an explicit derivative error in the improvement estimate so that it also applies in the algorithm.

\subsection{A structural dichotomy}
For positive vertex weights \(w\), introduce positive vertex parameters \(\theta_v\) and define \[ b(A)\defeq\sum_{v\in A}\theta_vw(v),\qquad D_\theta\defeq\sum_v\theta_vw(v)\Lambda_v =\sum_{xy\in E(G)}w(x)w(y)(\theta_x+\theta_y), \] and \[ e_w(A,B)\defeq\sum_{x\in A}\sum_{y\in B\cap N(x)}w(x)w(y). \]
The quantity \(e_w(A,B)\) measures the weighted edges from \(A\) to \(B\).
The sets may overlap, in which case an edge with both endpoints in \(A\cap B\) is counted in both directions.
The parameters \(\theta_v\) are arbitrary in the next lemma; in \Cref{cor:stationary} we choose \(\theta_v\defeq2r\xi_v/w(v)\).

\begin{lemma}\label{lem:triangle-test}
There is an absolute rational constant \(C_\triangle\ge1\) with the following property.
For every graph \(G\) with positive real weights \(w,\theta\), one of the following holds:
\begin{enumerate}[label=(\roman*),ref=(\roman*)]
\item\label{case:triangle-test-bound} \(T_G(w)\le C_\triangle D_\theta\).
\item\label{case:triangle-test-witness} There are sets \(A,B\subseteq V(G)\), an integer \(p\ge1\), and \(\rho\in\{\theta_v:v\in V(G)\}\) such that \begin{gather} w(A),w(B)\le\rho2^p,\qquad b(B)\le\rho^22^{2p},\label{alg:size}\\ e_w(A,B)>24\bigl(pb(A)+\rho^22^{-p}\bigr).\label{alg:strong-witness} \end{gather}
\end{enumerate}
If \(G\) has no edges, including when it is empty, \ref{case:triangle-test-bound} holds.
\end{lemma}

For the algorithm, we also need a bound on \(p\) and an effective way to find the witnesses.
An edgeless graph satisfies \ref{case:triangle-test-bound} immediately, so the algorithmic statement below assumes that \(G\) has an edge.
The proof of the following lemma is completed in \Cref{sec:derandomization}, using the walks and reductions in \Cref{sec:triangle-walks,sec:triangle-proof}.
We call the algorithm in this lemma the \emph{dichotomy algorithm}.

\begin{lemma}\label{lem:triangle-test-algorithm}
Let \(G\) be a graph containing an edge, with positive rational weights \(w,\theta\), and let \(C_\triangle\) be the constant in \Cref{lem:triangle-test}.
Put \[ t_0 \defeq \max\{1,\lceil\log_2\max_v(\Lambda_v/\theta_v)\rceil\}. \]
A deterministic algorithm can either verify \ref{case:triangle-test-bound} or return \(A,B,p,\rho\) satisfying \ref{case:triangle-test-witness} with \(p\le30t_0^3\), in time polynomial in the binary input length, with an absolute polynomial degree.
\end{lemma}

In the second alternative, we will increase the weights on a subset of \(A\) and set the weights on a disjoint subset of \(B\) to zero.
The large edge weight between the sets will yield a potential increase when the partial derivatives of \(\mathcal U\) vanish, and an explicit lower bound on the change when they are small.
The bounds on \(w(B)\) and \(b(B)\) control the losses on the set whose weights are set to zero.
On the \(A\) side, the neighborhood restriction in the next subsection will control the internal edge cost of increasing the weights.

\subsection{Restricting to a set with few internal edges}
Increasing several vertex weights incurs a quadratic cost from edges between them.
We record the neighborhood restriction used in Kang's proof of \cite[Lemma 3.1]{kang2026ode}.
As outlined in \ref{step:overview-b-sparse}, it makes this internal edge weight small while retaining enough vertex weight.
This is the only ingredient that uses \(K_r\)-freeness.

\begin{lemma}\label{lem:sparse-restriction}
Let \(r\ge2\), let \(G\) be a \(K_r\)-free graph with positive real weights \(w\), let \(Y\subseteq V(G)\) be nonempty, and let \(Q\ge1\).
There is a nonempty set \(C\subseteq Y\) such that \[ w(C)\ge Q^{-(r-2)}w(Y),\qquad w(N(v)\cap C)\le w(C)/Q\quad(v\in C). \]
For rational \(w,Q\), such a set can be found by a deterministic algorithm whose running time is polynomial in the binary input length, with an absolute polynomial degree.
\end{lemma}
\begin{proof}
Start with \(C\defeq Y\).
If a vertex \(v\in C\) has \(w(N(v)\cap C)>w(C)/Q\), replace \(C\) by \(N(v)\cap C\).
Every restriction retains more than a \(1/Q\) fraction of the weight and leaves a nonempty set.
The chosen vertices form a clique, and the current set lies in their common neighborhood.
After \(r-2\) restrictions the current set is independent, since an edge in it would complete a \(K_r\).
Thus the procedure stops after at most \(r-2\) restrictions, proving both bounds.
Each restriction replaces \(C\) by a nonempty proper subset, since the chosen vertex is removed.
There are therefore at most \(\min\{r-2,\abs{Y}-1\}\) restrictions.
Each step uses only sums of the original rational weights and comparisons with \(w(C)/Q\), which proves the uniform polynomial running time assertion.
\end{proof}

\subsection{Increasing the potential using two disjoint vertex sets}
We now carry out the weight change in \ref{step:overview-b-change}.
To control the vertex terms when weights increase, we use the scalar tangent-loss estimate from \Cref{app:scalar}: \begin{equation}\label{eq:scalar-loss} U_a(x)+U_a'(x)(z-x)-U_a(z) \le4\xi_a(x)\frac{z-x}{x}\log(z/x)\qquad(a>0,\ 0<x\le z). \end{equation}
The next lemma assumes \(\abs{\partial\mathcal U_{G,u}/\partial w(v)}=\abs{\zeta_v}\le\eta\) at every vertex, where \(\eta\ge0\).
Taking \(\eta=0\) gives the version used at a maximizer.
Allowing a positive tolerance gives the version needed by the algorithm.
Fix an integer \(r\ge4\).
To bound the total change in the vertex weights, put \begin{equation}\label{eq:Jpw-definition} J_p(w)\defeq\bigl(1+2^{(12r-14)p}\bigr)w(V(G)). \end{equation}

\begin{lemma}\label{lem:improvement}
Let \(G\) be \(K_r\)-free, let \(u,w,\theta>0\), and fix \(\eta\ge0\).
Suppose \(\theta_vw(v)\ge 2r\xi_v\) for every \(v\), and \(\norm{\zeta}_\infty\le\eta\).
If \(A,B,p,\rho\) satisfy \eqref{alg:size} and \eqref{alg:strong-witness}, there is a nonnegative vertex weight \(w'\), zero on a nonempty set of vertices and at least \(w\) at every other vertex, such that \begin{gather*} \norm{w'-w}_1\le J_p(w),\\ \mathcal U(w')-\mathcal U(w)>2^{4p}p\,b(A)-\eta J_p(w). \end{gather*}
For rational \(w,\theta\), the vertex weight \(w'\) can be constructed deterministically from \(A,B,p,\rho\) in time polynomial in the binary input length and the numerical value of \(p\), with an absolute polynomial degree.
\end{lemma}
\begin{proof}
We first choose disjoint sets \(C\subseteq A\) and \(B'\subseteq B\) such that every vertex of \(C\) has enough neighbor weight in \(B'\), while \(C\) has few internal edges.
We then increase the weights on \(C\) and set those on \(B'\) to zero.
The edges between the sets provide the gain needed to absorb the scalar losses and internal edge costs.

\step{1. Make the sets disjoint.}
Assign each vertex independently to \(V_0\) or \(V_1\) with equal probability, and put \(A'\defeq A\cap V_0\), \(B'\defeq B\cap V_1\).
The expectation of \(e_w(A',B')-12p\,b(A')\) is \(e_w(A,B)/4-6p\,b(A)\).
Thus there is a partition for which \begin{equation}\label{eq:partition-gain} e_w(A',B')-12p\,b(A') \ge\frac{e_w(A,B)-24p\,b(A)}{4} >6\rho^22^{-p}. \end{equation}
The coefficient \(12\) makes the terms involving \(p\,b(A)\) cancel when applying \eqref{alg:strong-witness}.
Such a partition can be found by fixing the assignments one at a time, each time choosing a side that does not decrease the conditional expectation.
These expectations are sums of vertex and edge terms, so they are directly computable for rational data.
Since \(e_w(A',B')>0\), both \(A'\) and \(B'\) are nonempty.

\step{2. Choose the set whose weights will increase.}
We first retain only vertices with enough neighbor weight in \(B'\).
For \(v\in A'\), put \(d_{B'}(v)\defeq w(N(v)\cap B')\), and set \[ \delta\defeq\rho2^{-2p},\qquad A''\defeq\{v\in A':d_{B'}(v)>12p\theta_v+\delta\}. \]
The contribution of vertices in \(A'\setminus A''\) to \(e_w(A',B')\) is at most \(12p\,b(A'\setminus A'')+\delta w(A'\setminus A'')\), which is at most \(12p\,b(A')+\delta w(A')\).
Consequently, \eqref{eq:partition-gain} and \(w(A')\le\rho2^p\) give \[ e_w(A'',B') \ge e_w(A',B')-12p\,b(A')-\delta w(A') >5\rho^22^{-p}. \]
Since \(e_w(A'',B')\le w(A'')w(B')\le w(A'')\rho2^p\), we have \(w(A'')>5\rho2^{-2p}\).
In particular, \(A''\) is nonempty.

To control the internal edge cost, apply \Cref{lem:sparse-restriction} to \(A''\) with \(Q\defeq2^{12p}\), obtaining a nonempty set \(C\subseteq A''\) with \[ w(C)>\rho2^{-(12r-22)p},\qquad 2M_{G[C]}(w)\le w(C)^2/Q. \]

\step{3. Choose the weight change.}
Put \(\alpha\defeq\rho2^{8p}/w(C)\) and define \[ w'(v)\defeq\begin{cases} (1+\alpha)w(v),&v\in C,\\ 0,&v\in B',\\ w(v),&v\notin C\cup B'. \end{cases} \]
This change adds total weight \(\alpha w(C)=\rho2^{8p}\) to \(C\).
We have \(\alpha<2^{(12r-14)p}\), so \(\norm{w'-w}_1\le J_p(w)\) and \begin{equation}\label{eq:weight-change-log-bound} \log(1+\alpha)\le\bigl(1+(12r-14)p\bigr)\log2\le12rp. \end{equation}

\step{4. Compare the gain and the losses.}
Write \(\Delta\defeq w'-w\).
For each vertex \(v\), define the \emph{tangent loss} at \(v\) by \[ \ell_v\defeq U_{u(v)}(w(v))+U'_{u(v)}(w(v))\Delta_v-U_{u(v)}(w'(v)). \]
This is the error in approximating \(U_{u(v)}(w'(v))\) by the tangent line at \(w(v)\).
By concavity, \(\ell_v\ge0\), and \(\ell_v=0\) for \(v\notin C\cup B'\).
A direct calculation gives \begin{equation}\label{eq:potential-change} \mathcal U(w')-\mathcal U(w) =\sum_v\zeta_v\Delta_v-\sum_v\ell_v -2\sum_{xy\in E(G)}\Delta_x\Delta_y. \end{equation}
The linear term is at least \(-\eta\sum_v\abs{\Delta_v}\ge-\eta J_p(w)\).
For \(v\in C\), we have \(\Delta_v=\alpha w(v)\) and \(w'(v)=(1+\alpha)w(v)\).
Thus \eqref{eq:scalar-loss} and \eqref{eq:weight-change-log-bound} together with the hypothesis \(\theta_vw(v)\ge2r\xi_v\) give \begin{equation}\label{eq:increased-tangent-loss} \sum_{v\in C}\ell_v \le48r\alpha p\sum_{v\in C}\xi_v\le24\alpha p\,b(C). \end{equation}
For \(v\in B'\), we have \(w'(v)=0\) and \(U_{u(v)}(0)=0\), so \(\ell_v=\xi_v\).
Writing \(\xi(B')\defeq\sum_{v\in B'}\xi_v\), we obtain \begin{equation}\label{eq:deleted-tangent-loss} \sum_{v\in B'}\ell_v=\xi(B')\le b(B')/(2r). \end{equation}

In the quadratic edge term of \eqref{eq:potential-change}, edges from \(C\) to \(B'\) contribute a gain \(2\alpha\,e_w(C,B')\).
Edges inside \(C\) cost \(2\alpha^2M_{G[C]}(w)\), and edges inside \(B'\) cost \(2M_{G[B']}(w)\).
All other products \(\Delta_x\Delta_y\) vanish.
Every vertex of \(C\) lies in \(A''\), so \(e_w(C,B')>12p\,b(C)+\delta w(C)\).
The term \(24\alpha p\,b(C)\) in the cross-edge gain \(2\alpha\,e_w(C,B')\) covers the tangent losses on \(C\), leaving a surplus \(>2\alpha\delta w(C)\).
Apply \eqref{eq:increased-tangent-loss} and \eqref{eq:deleted-tangent-loss} together with the edge bounds above to \eqref{eq:potential-change}.
Combining these estimates with \(\sum_v\zeta_v\Delta_v\ge-\eta J_p(w)\), we obtain \[ \mathcal U(w')-\mathcal U(w) >-\eta J_p(w)-b(B')/(2r) +2\alpha\delta w(C)-2\alpha^2M_{G[C]}(w)-2M_{G[B']}(w). \]
The terms satisfy \[ \begin{aligned} 2\alpha\delta w(C)&=2\rho^2 2^{6p},\\
    2\alpha^2M_{G[C]}(w)&\le\alpha^2w(C)^2/Q=\rho^2 2^{4p},\\
    \frac{b(B')}{2r}+2M_{G[B']}(w) &\le\frac{b(B')}{2r}+w(B')^2\\
    &\le\Bigl(1+\frac{1}{2r}\Bigr)\rho^2 2^{2p}\le2\rho^2 2^{2p}, \end{aligned} \] where the second line uses \(2M_{G[C]}(w)\le w(C)^2/Q\), and the last line uses \eqref{alg:size} and \(r\ge4\).
Since \(p\ge1\), we have \(2\cdot2^{6p}-2^{4p}-2\cdot2^{2p}>2^{6p}\), and hence \[ \mathcal U(w')-\mathcal U(w)>\rho^2 2^{6p}-\eta J_p(w). \]
Finally, \eqref{alg:strong-witness} and \eqref{alg:size} give \[ pb(A)<e_w(A,B)\le w(A)w(B)\le\rho^22^{2p}. \]
It follows that \(2^{4p}p\,b(A)<\rho^2 2^{6p}\), which proves the stated estimate on the potential difference.

The partition uses at most \(\abs{V(G)}\) conditional-expectation decisions, and the neighborhood restriction uses at most \(\abs{V(G)}-1\) proper subset restrictions.
Every exponent of two used here has absolute value at most \(12rp\), so the numerator and denominator of the corresponding power of two together use at most \(12rp+2\) bits in binary.
The restriction leaves the weights unchanged, and the final rescaling uses sums of original weights followed by a bounded number of products and quotients per vertex.
The construction explicitly forms only the \(p\)-dependent powers \(2^{-2p}\), \(2^{12p}\), and \(2^{8p}\).
Powers involving \(r\), such as \(2^{-(12r-22)p}\), occur only in the estimates.
Thus both the operation count and the encoding lengths are bounded by polynomials of absolute degree in the binary input length and \(p\).
\end{proof}

The estimate \(\log(1+\alpha)\le12rp\) uses the fact that at most \(r-2\) neighborhood restrictions are made.
As explained in \ref{step:overview-b-change}, this leads to the factor \(2r\) in the hypothesis on \(\theta\), and hence to the dependence on \(r\) stated in \Cref{thm:existence}.

\subsection{The triangle bound at a maximizer}
We now specialize to a positive maximizer and choose \(\theta_v\defeq 2r\xi_v/w(v)\), completing \ref{step:overview-b-conclude}.
This choice both rules out alternative~\ref{case:triangle-test-witness} of \Cref{lem:triangle-test} and bounds \(D_\theta\) in terms of \(S\).

\begin{corollary}\label{cor:stationary}
Let \(G\) be \(K_r\)-free and \(u>0\).
If \(w>0\) is a global maximizer of \(\mathcal U_{G,u}\) over nonnegative vertex weights, then, with \(\theta_v\defeq 2r\xi_v/w(v)\), \[ T_G(w)\le 2C_\triangle rS. \]
\end{corollary}
\begin{proof}
Since \(w\) is a positive maximizer, all the partial derivatives \(\zeta_v\) are zero.
If sets \(A,B\) and parameters \(p,\rho\) satisfied \eqref{alg:size} and \eqref{alg:strong-witness}, \Cref{lem:improvement} with \(\eta=0\) would give \(w'\ge0\) with \(\mathcal U(w')>\mathcal U(w)\), contradicting maximality.
Hence \(T_G(w)\le C_\triangle D_\theta\).
Moreover, \(\Lambda_v=w(v)d_{*,v}\) because \(\zeta_v=0\), and \(d_{*,v}\xi_v\le u(v)\) by \eqref{eq:scalar-ode}.
Thus \[ D_\theta=2r\sum_v\xi_v\Lambda_v =2r\sum_vw(v)d_{*,v}\xi_v\le 2r\sum_vu(v)w(v)=2rS, \] and the conclusion follows.
\end{proof}


\section{The existence proof}\label{sec:existence}
Assuming the structural dichotomy in \Cref{lem:triangle-test}, we complete the proof of \Cref{thm:existence} by induction, following \Cref{subsec:overview-induction,subsec:overview-increase}.
We divide into two cases according to whether the limit of a normalized maximizing sequence has a zero coordinate or is positive.

For the rest of the paper, choose \begin{equation}\label{eq:constant} c\defeq c_r\defeq \frac{1}{1+96C_\triangle r}. \end{equation}
Then \(c_r\ge\kappa_1/r\) with \(\kappa_1\defeq1/(1+96C_\triangle)\), proving the final assertion of \Cref{thm:existence}.
The existence proof only requires \(c\le(1+12C_\triangle r)^{-1}\), as noted in \ref{step:overview-b-conclude}.
The smaller choice in \eqref{eq:constant} also bounds the derivative and rounding errors in the algorithm.

For fixed objective weights \(u\), put \[ F(G,u)\defeq \sup_{s>0}\Phi_{G,u}(s). \]
Since \(0<f\le1\), the supremum is finite.
We prove, for every \(K_r\)-free graph \(G\) and every nonnegative objective weight \(u\), that \(\alpha_u(G)\ge cF(G,u)\), by induction on \(\abs{V(G)}\); the empty case is immediate.
Deleting a vertex of zero objective weight leaves \(\alpha_u\) unchanged and cannot decrease \(F\).
An isolated vertex contributes \(u(v)\) to both \(\alpha_u\) and \(F\), so induction also handles this case since \(c\le1\).
Thus we may assume \(u>0\) and that \(G\) has no isolated vertices.

Choose positive vertex weights \(s^{(j)}\) with \(\Phi_{G,u}(s^{(j)})\to F(G,u)\).
Weighted degrees do not change under a common rescaling, so normalize \(s^{(j)}(V(G))=1\).
These vertex weights lie in the closed bounded set \(\{s\ge0:\sum_vs(v)=1\}\), so they have a convergent subsequence.
Write \(s^*\ge0\) for its limit; it also has total weight one.
Let \(A\defeq \{v:s^*(v)>0\}\neq\emptyset\).

\subsection{A boundary limit}\label{subsec:boundary-limit}
Suppose first that \(A\ne V(G)\).
We will bound \(F(G,u)\) by the sum of the corresponding suprema on two smaller induced subgraphs with no edges between them, so that induction applies.
Put \(B\defeq V(G)\setminus N_G[A]\).
There are three kinds of summands to consider.
By continuity, the total contribution from \(A\) tends to \(\Phi_{G[A],u}(s^*|_A)\).
For \(v\in N_G(A)\setminus A\), the denominator tends to zero while the numerator has a positive limit, so \(d_{G,s^{(j)}}(v)\to\infty\) and the summand tends to zero.
Finally, for vertices of \(B\), deleting the other vertices cannot increase their weighted degrees and hence cannot decrease their \(f\)-values, so \[ \sum_{v\in B}u(v)f(d_{G,s^{(j)}}(v)) \le\Phi_{G[B],u}(s^{(j)}|_B)\le F(G[B],u). \]
Taking the limit superior of the whole sum yields \[ F(G,u)\le F(G[A],u)+F(G[B],u). \]
Both induced graphs are smaller than \(G\), and there are no edges between \(A\) and \(B\).
By induction, the union of independent sets in these two graphs has weight at least \(cF(G,u)\), as required.

\subsection{A positive maximizer}\label{subsec:positive-maximizer}
It remains to consider \(A=V(G)\), when continuity shows that \(s^*\) attains \(F(G,u)\).
We will find a vertex \(x\) with \(g_x\ge0\), so that selecting it and applying induction to \(G_x\) reaches the target.
There are no isolated vertices, so \Cref{cor:transfer} supplies a positive global maximizer \(w\) of \(\mathcal U\), with \(\mathcal U_G(w)=F(G,u)\) and \(\zeta=0\).
\Cref{cor:stationary} gives \(T_G(w)\le 2C_\triangle rS\).
By \Cref{lem:deletion}, \[ \sum_xw(x)g_x=(1-c)S-6cT_G(w) \ge(1-c-12cC_\triangle r)S\ge0, \] where the last inequality follows from \eqref{eq:constant}.
All \(w(x)\) are positive, so at least one vertex \(x\) has \(g_x\ge0\).
For \(G_x=G-N[x]\), induction and \eqref{eq:potential-comparison} give \[ \alpha_u(G_x)\ge cF(G_x,u)\ge c\mathcal U_{G_x,u}(w). \]
An independent set in \(G_x\) can be extended by \(x\).
Consequently \[ \alpha_u(G)\ge u(x)+c\mathcal U_{G_x,u}(w) \ge c\mathcal U_{G,u}(w)=cF(G,u). \]
This completes the induction and proves \Cref{thm:existence}.


\section{The algorithmic proof}\label{sec:outer}
Assuming the dichotomy algorithm in \Cref{lem:triangle-test-algorithm}, we prove \Cref{thm:algorithm} by constructing an independent set \(I\) satisfying
\begin{equation}\label{eq:algorithm-bound}
    u(I)\ge\frac{c_r}{2}\Phi_{G,u}(s).
\end{equation}
We replace the maximizer in the existence proof with the three operations outlined in \Cref{subsec:overview-algorithm}.
We first control extreme weights and improve one coordinate whenever its partial derivative is far from zero.
When all partial derivatives are small, the dichotomy algorithm lets us select a vertex or reweight two sets and remove the vertices whose weights become zero.
We give the complete algorithm below, then justify its operations and bound its running time and total error.
Finally, we use the weighted independent-set algorithm to construct the local-demand distribution, as outlined at the end of \Cref{subsec:overview-algorithm}.

If \(u\equiv0\), return the empty independent set.
Otherwise normalize \(\sum_vu(v)=1\), which rescales both sides of the desired inequality \eqref{eq:algorithm-bound} by the same factor.
Let \(n\) be the number of vertices of the input graph.
The empty graph has already been handled; if \(n=1\), return its only vertex, using \(c_r\le1\) and \(f(0)=1\).
Henceforth \(n\ge2\); both \(n\) and the normalized objective weights \(u\) remain fixed throughout the loop.
Delete vertices with \(u(v)=0\); this cannot decrease the target sum.

The algorithm maintains an independent set \(I\) and a remaining induced graph \(H\defeq G[R]\), with \(R\cap N_G[I]=\emptyset\).
All quantities below are computed in the current graph \(H\); \(G\) always denotes the input graph.
Selecting a vertex means adding it to \(I\) and removing its closed neighborhood from \(H\), whereas discarding it means removing it without adding it to \(I\).
We track \begin{equation}\label{alg:tracked-potential} P\defeq u(I)+c\,\mathcal U_{H,u}(w)\le u(I)+\sum_{v\in R}u(v)\le1, \end{equation} where we used \(U_{u(v)}(w(v))\le u(v)\) and \(c\le1\).
Selecting \(x\) changes \(P\) by the gain \[ g_x = u(x)+c\,\mathcal U_{H-N_H[x],u}(w) -c\,\mathcal U_{H,u}(w). \]
Initialization and discards will incur a bounded total loss, while every operation that leaves the vertex set unchanged will increase \(P\) by a uniform positive amount.
\Cref{alg:weighted-independent-set} gives the complete procedure, including the one-time initialization.
Each pass through its loop performs one operation, then repeats the checks with the current graph and auxiliary weights.
We justify these operations in \Cref{subsec:algorithm-weights,subsec:algorithm-select-or-reweight}, following \ref{step:overview-c-range}--\ref{step:overview-c-select-or-change}, and bound the total error and number of iterations in \Cref{subsec:algorithm-analysis}.

In the pseudocode, steps~\ref{step:weighted-discard} and~\ref{step:weighted-select} implement \ref{step:overview-c-range} by discarding small-weight vertices and selecting isolated or large-weight vertices.
Step~\ref{step:weighted-coordinate} implements \ref{step:overview-c-derivative} by updating one auxiliary weight, with a discard or selection if the update crosses a threshold.
Step~\ref{step:weighted-triangle-test} implements \ref{step:overview-c-select-or-change} by applying the dichotomy algorithm and then either selecting a vertex or reweighting two sets.

\bigskip
\Needspace{6\baselineskip}
\begin{breakablealgorithm}
\caption{Weighted independent-set algorithm}\label{alg:weighted-independent-set}
\noindent\textbf{Input}: An integer \(r\ge4\), an \(n\)-vertex \(K_r\)-free graph \(G\), and rational weights \(s>0\) and \(u\ge0\).

\noindent\textbf{Output}: An independent set of objective weight at least \((c_r/2)\Phi_{G,u}(s)\).

\begin{enumerate}[itemsep=5pt,ref=(\arabic*)]
    \item \textbf{Initialization:}
    \begin{enumerate}[itemsep=5pt,ref=(\arabic{enumi}\alph*)]
        \item If \(u\equiv0\), return \(\emptyset\).
        Record the input length \(L\).
        If \(n=1\), return \(V(G)\).
        Normalize \(u\) once so that \(\sum_{v\in V(G)}u(v)=1\).
        Fix a vertex order for all choices and ties, and retain \(u\) and \(n\) throughout the iterations.

        \item Set \(\lambda\defeq(100n)^{-4}\) and \(c\defeq c_r\), with \(c_r\) from \eqref{eq:constant}.
        Fix \(P_n,J,\eta,\ell,\gamma\) by \eqref{alg:witness-exponent}, \eqref{alg:tolerance-parameters}, and \eqref{alg:coordinate-parameters}, and fix the rounding grid in \eqref{eq:numerics-rounding-grid}.

        \item Set \(I\defeq\emptyset\) and \(H\defeq G[\{v:u(v)>0\}]\).
        Add all isolated vertices of \(H\) to \(I\) and remove them from \(H\).
        If \(V(H)\ne\emptyset\), apply \Cref{lem:lift} to \(H\) with the input weights \(s\) restricted to \(V(H)\), and initialize \(w\) at the certified upper grid approximations constructed in \Cref{subsec:numerics-initialization}.
    \end{enumerate}

    \item \textbf{Iteration:}
    Execute the following steps in order.
    Selecting a vertex \(v\) means adding it to \(I\) and removing \(N_H[v]\) from \(H\).

    \begin{enumerate}[itemsep=5pt,ref=(\arabic{enumi}\alph*)]
        \item\label{step:weighted-check}
        If \(V(H)=\emptyset\), proceed to \ref{step:weighted-output}.

        \item\label{step:weighted-discard} %
        If some \(v\in V(H)\) has \(u(v)<\lambda\) or \(w(v)\le\lambda\), discard one such vertex and return to \ref{step:weighted-check}.

        \item\label{step:weighted-select} %
        If some \(v\in V(H)\) is isolated or has \(w(v)\ge\lambda^{-1}\), select one such vertex and return to \ref{step:weighted-check}.

        \item\label{step:weighted-coordinate} \textbf{Improve one weight:}
        Compute rational approximations \(\widehat\zeta_v\) to \(\zeta_v\) for all \(v\in V(H)\), satisfying \(\abs{\widehat\zeta_v-\zeta_v}\le\eta/8\).
        If \(\abs{\widehat\zeta_v}\le3\eta/4\) at every vertex, proceed to \ref{step:weighted-triangle-test}.
        Otherwise, choose a vertex \(v\) with \(\abs{\widehat\zeta_v}>3\eta/4\) and define \( \widetilde w(v)\defeq w(v)+\ell\operatorname{sgn}(\widehat\zeta_v) \).
        Apply exactly one of the following and return to \ref{step:weighted-check}:
        \begin{enumerate}[itemsep=5pt,ref=(\arabic{enumi}\alph{enumii}-\roman*)]
            \item\label{step:weighted-full-coordinate}
            If \(\lambda\le\widetilde w(v)\le\lambda^{-1}\), replace \(w(v)\) by \(\widetilde w(v)\) and round the positive weights upward to the fixed grid.
            \item\label{step:weighted-coordinate-discard} If \(\widetilde w(v)<\lambda\), replace \(w(v)\) by \(\lambda\) and discard \(v\), without rounding.
            \item\label{step:weighted-coordinate-select} If \(\widetilde w(v)>\lambda^{-1}\), replace \(w(v)\) by \(\lambda^{-1}\) and select \(v\), without rounding.
        \end{enumerate}
        \item\label{step:weighted-triangle-test} \textbf{Select a vertex or reweight two sets:}
        Compute \(\bar\xi_v\) and \(\theta_v\) as in \eqref{eq:numerics-xi-approximation} and \eqref{alg:theta-parameters}, and run the dichotomy algorithm in \Cref{lem:triangle-test-algorithm} on \(H\) with \(w,\theta\).
        Apply exactly one of the following and return to \ref{step:weighted-check}:
        \begin{enumerate}[itemsep=5pt,ref=(\arabic{enumi}\alph{enumii}-\roman*)]
            \item If it verifies \(T_H(w)\le C_\triangle D_\theta\), approximate every \(g_x\) to error at most \(\lambda^3/(100n)\) and select a vertex with largest approximate gain.
            \item\label{step:weighted-two-set-reweighting}
            Otherwise, let \(A,B,p,\rho\) be the returned witness.
            Construct \(w'\) from \Cref{lem:improvement} and replace \(w\) by \(w'\).
            Round the positive weights upward to the fixed grid, leaving zero weights zero, and remove all zero-weight vertices from \(H\).
        \end{enumerate}
    \end{enumerate}

    \item\label{step:weighted-output} \textbf{Output:}
    Return the heavier of \(I\) and a singleton of largest objective weight in \(G\).
\end{enumerate}
\end{breakablealgorithm}
\bigskip

In step~\ref{step:weighted-coordinate-discard}, setting \(w(v)\) to \(\lambda\) before discarding \(v\) serves only to simplify the analysis by treating steps~\ref{step:weighted-coordinate-discard} and~\ref{step:weighted-coordinate-select} together.
This assignment may be omitted in an implementation.

\subsection{Controlling and updating the weights}\label{subsec:algorithm-weights}
We implement \ref{step:overview-c-range} and \ref{step:overview-c-derivative} in order: first restrict the surviving weights to a bounded range, then choose a weight update that increases the potential.

\subsubsection{Control extreme weights}
We justify the discards and selections in steps~\ref{step:weighted-discard} and \ref{step:weighted-select}.
First, we discard vertices with \(u(v)<\lambda\) or \(w(v)\le\lambda\), recalling \(\lambda=(100n)^{-4}\) from \Cref{alg:weighted-independent-set}, and bound the resulting loss in \(P\).
Discarding a vertex with \(u(v)<\lambda\) loses at most \(\lambda\) from \(\mathcal U_{H,u}\), since \(U_{u(v)}(w(v))\le u(v)\).
Discarding one with \(w(v)\le\lambda\) loses at most \[ U_1(\lambda)\le f(1/\lambda)+\lambda \le\lambda\bigl(1+2\log(1/\lambda)\bigr). \]
Here we used \(u(v)\le1\), monotonicity of \(U_a(x)\) in both variables, the bound \(U_1(\lambda)\le f(\lambda^{-1})+\lambda\), obtained by taking \(t=\lambda^{-1}\) in the infimum, and \(f(d)\le2\log d/d\) for \(d\ge2\) from \Cref{app:scalar}.
Removing incident edges can only increase the potential.
Since each vertex is discarded at most once, the total loss in \(\mathcal{U}_{H,u}\), and hence in \(P\) since \(c\le1\), is at most \begin{equation}\label{alg:discard-loss} n \cdot \lambda\bigl(1+2\log(1/\lambda)\bigr)<\frac{1}{100n}. \end{equation}
For \(\lambda=(100n)^{-4}\), the last inequality holds for every \(n\ge2\).

After these discards, a nonisolated vertex \(x\) with \(w(x)\ge\lambda^{-1}\) has a surviving neighbor \(y\) with \(w(y)>\lambda\), so \(2w(x)\Lambda_x\ge2\).
Using \(U_{u(x)}(w(x))\le u(x)\) and retaining only the edge gain from edges incident to \(x\), we obtain \[ g_x\ge(1-c)u(x)-c\sum_{y\in N_H(x)} u(y) +2cw(x)\Lambda_x\ge c>0, \] where we also use \(\sum_{y\in N_H(x)}u(y)\le1\).
An isolated vertex \(x\) also has nonnegative gain, since \(g_x=u(x)-c\,U_{u(x)}(w(x))\ge0\).
Thus, after the discards above, we select an isolated vertex or a vertex \(x\) with \(w(x)\ge\lambda^{-1}\) whenever one exists.
After each discard or selection, we return to step~\ref{step:weighted-check}; we proceed to step~\ref{step:weighted-coordinate} only when neither a discard nor a selection applies.
Whenever none applies, the nonempty current graph has no isolated vertices and satisfies \begin{equation}\label{alg:box} u(v)\ge\lambda,\qquad \lambda<w(v)<\lambda^{-1}. \end{equation}
In particular, \(S\ge\lambda^2\), and \eqref{eq:curvature} gives \(\abs{U''_{u(v)}(x)}\le2/\lambda^2\) for all \(x\in[\lambda,\lambda^{-1}]\).

\subsubsection{Improve one weight}
Let \(\eta>0\) be a tolerance whose value we fix in \Cref{subsec:algorithm-select-or-reweight}.
Define the step length and the required increase by \begin{equation}\label{alg:coordinate-parameters} \ell\defeq\frac{\eta\lambda^2}{16},\qquad \gamma\defeq\frac{\eta^2\lambda^2}{100}. \end{equation}
Compute rational approximations \(\widehat{\zeta}_v\) to the current partial derivatives with \begin{equation}\label{alg:derivative-accuracy} \abs{\widehat\zeta_v-\zeta_v}\le\eta/8. \end{equation}
If \(\abs{\widehat\zeta_v}\le3\eta/4\) at every vertex, then \(\norm{\zeta}_\infty\le7\eta/8<\eta\), and we proceed to the dichotomy algorithm in \Cref{subsec:algorithm-select-or-reweight}.
Otherwise choose \(v\) with \(\abs{\widehat\zeta_v}>3\eta/4\).
Then \(\zeta_v\) has the same sign as \(\widehat\zeta_v\) and absolute value greater than \(5\eta/8\).
Define \(\widetilde w(v)\) as \(w(v)+\ell\) if \(\widehat\zeta_v>0\) and \(w(v)-\ell\) otherwise.
Define \(w'\) on \(V(H)\) by \(w'(y)\defeq w(y)\) for \(y\ne v\) and the following choices at \(v\), labeled by the corresponding steps of \Cref{alg:weighted-independent-set}:
\begin{enumerate}[leftmargin=4em]
\item[\ref{step:weighted-full-coordinate}] if \(\lambda\le\widetilde w(v)\le\lambda^{-1}\), set \(w'(v)\defeq\widetilde w(v)\);
\item[\ref{step:weighted-coordinate-discard}] if \(\widetilde w(v)<\lambda\), set \(w'(v)\defeq\lambda\) and discard \(v\);
\item[\ref{step:weighted-coordinate-select}] if \(\widetilde w(v)>\lambda^{-1}\), set \(w'(v)\defeq\lambda^{-1}\), select \(v\), and remove \(N[v]\).
\end{enumerate}
In step~\ref{step:weighted-full-coordinate}, put \(\delta\defeq w'(v)-w(v)=\pm\ell\).
By Taylor's theorem, for some \(x_*\) between \(w(v)\) and \(w'(v)\), \[ \mathcal U_{H,u}(w')-\mathcal U_{H,u}(w) =\zeta_v\delta+\frac{1}{2}U''_{u(v)}(x_*)\delta^2. \]
Since \(\delta\) has the same sign as \(\zeta_v\) and \(\abs{U''_{u(v)}(x_*)}\le2/\lambda^2\), the potential increase is at least \[ \abs{\zeta_v}\ell-\ell^2/\lambda^2 \ge (5\eta/8)\ell-\ell^2/\lambda^2 = 9\eta^2\lambda^2/256 > \gamma. \]
In steps~\ref{step:weighted-coordinate-discard} and~\ref{step:weighted-coordinate-select}, consider \(w'\) before vertex removal.
Put \(\delta\defeq\abs{w'(v)-w(v)}\), so \(0<\delta<\ell\).
The same Taylor estimate gives \[ \mathcal U_{H,u}(w')-\mathcal U_{H,u}(w) \ge\frac{5\eta}{8}\delta-\frac{\delta^2}{\lambda^2} =\delta\biggl(\frac{5\eta}{8}-\frac{\delta}{\lambda^2}\biggr)>0, \] since \(\delta/\lambda^2<\ell/\lambda^2=\eta/16\).
Thus every coordinate update above increases the potential before any vertex removal.
The subsequent discard loss is included in \eqref{alg:discard-loss}, since \(w'(v)=\lambda\).

\subsection{Selecting a vertex or reweighting two sets}\label{subsec:algorithm-select-or-reweight}
Next, we implement \ref{step:overview-c-select-or-change} when the current graph is nonempty, \eqref{alg:box} holds, and \(\norm{\zeta}_\infty\le\eta\).
We first choose the parameters for the dichotomy algorithm and fix \(\eta\) so that both outcomes give progress.

\subsubsection{Choose the parameters}
Write \(h\defeq-f'\) and \(d_*\defeq d_*(u(v),w(v))\) when working at one vertex.
By \eqref{eq:scalar}, we have \(h(d_*)=w(v)^2/u(v)\ge\lambda^2\).
The bound \(d_*h(d_*)\le1\) from \eqref{eq:scalar-h-bound} then gives \(d_*\le\lambda^{-2}\).
Using \(f(d)-dh(d)\ge1/(3\max\{1,d\})\) from \eqref{eq:scalar-g-bounds}, we obtain \begin{equation}\label{alg:xi-lower} \xi_v=u(v)\bigl(f(d_*)-d_*h(d_*)\bigr)\ge\lambda^3/3. \end{equation}
The construction in \eqref{eq:numerics-xi-approximation}, together with the bound in \eqref{eq:numerics-xi-upper-bound}, gives rational parameters satisfying
\begin{equation}\label{alg:theta-parameters}
    \xi_v\le\bar\xi_v\le2\xi_v,\qquad
    \theta_v\defeq\frac{2r\bar\xi_v}{w(v)}.
\end{equation}
These parameters preserve the hypothesis \(\theta_vw(v)\ge2r\xi_v\) of \Cref{lem:improvement}, and all tests in the dichotomy algorithm use exact rational arithmetic.

We next bound the size of a possible reweighting, so that its derivative error can be absorbed.
Recall from \eqref{eq:scalar} and the definition of \(\zeta_v\) that \(\Lambda_v=w(v)d_*-\zeta_v/2\), which is at most \(w(v)d_*+\eta/2\) by \(\abs{\zeta_v}\le\eta\).
Together with \eqref{alg:xi-lower} and \eqref{alg:theta-parameters}, this gives \[ \frac{\Lambda_v}{\theta_v} =\frac{w(v)\Lambda_v}{2r\bar\xi_v} \le\frac{w(v)^2d_*}{2r\xi_v}+\frac{\eta w(v)}{4r\xi_v} \le\frac{d_*h(d_*)}{2r(f(d_*)-d_*h(d_*))} +\frac{3\eta}{4r\lambda^4}. \]
Here we also used \(w(v)^2=u(v)h(d_*)\) and \(w(v)\le\lambda^{-1}\).

We bound the fraction involving \(d_*\) as follows.
For \(d_*\le1\), use \(d_*h(d_*)\le1\) from \eqref{eq:scalar-h-bound}.
For \(d_*>1\), use \(d_*^2h(d_*)\le\log(1+d_*)+1/2\) from \eqref{eq:scalar-h-log-bound}.
Combining each bound with the lower bound on \(f-dh\) in \eqref{eq:scalar-g-bounds} gives \[ \frac{d_*h(d_*)}{2r(f(d_*)-d_*h(d_*))} \le\begin{cases} \frac{3}{2r},&0<d_*\le1,\\[4pt] \frac{3}{2r}\bigl(\log(1+d_*)+\tfrac{1}{2}\bigr),&d_*>1. \end{cases} \]
To control the error term \(3\eta/(4r\lambda^4)\), we now require \(\eta\le\lambda^4/(100n^2)\).
Using this bound together with \(d_*\le\lambda^{-2}\), \(r\ge4\), and the definition of \(\lambda\), we obtain \begin{equation}\label{alg:load} \frac{\Lambda_v}{\theta_v} \le\frac{3}{4}\log(1/\lambda)+1 \le10\log(e/\lambda)\le400\log n. \end{equation}
If the algorithm in \Cref{lem:triangle-test-algorithm} returns alternative~\ref{case:triangle-test-witness} of \Cref{lem:triangle-test}, it produces sets \(A,B\) and parameters \(p,\rho\) satisfying \eqref{alg:size} and \eqref{alg:strong-witness}.
Put
\[
    b_n\defeq\lceil\log_2(2n)\rceil,\qquad
    q_n\defeq\lceil\log_2(400b_n)\rceil.
\]
Recall from \Cref{lem:triangle-test-algorithm} that, for the current graph \(H\), weights \(w\), and parameters \(\theta\), the dichotomy algorithm uses
\[
    t_0\defeq\max\bigl\{1,\lceil\log_2
    \max_{v\in V(H)}(\Lambda_v/\theta_v)\rceil\bigr\}.
\]
By \eqref{alg:load} and the definitions of \(b_n,q_n\), we have \(t_0\le q_n\).
Thus, if the dichotomy algorithm returns a witness, its parameter \(p\) satisfies
\begin{equation}\label{alg:witness-exponent}
    p\le30t_0^3\le30q_n^3\eqdef P_n.
\end{equation}
Keeping this bound, of order \((1+\log\log(2n))^3\), will give a polynomial degree independent of \(r\).

We can now fix \begin{equation}\label{alg:tolerance-parameters} J\defeq n\lambda^{-1}\bigl(1+2^{(12r-14)P_n}\bigr),\qquad \eta\defeq\min\biggl\{\frac{\lambda^4}{100n^2},\frac{\lambda^3}{100J}\biggr\}. \end{equation}
This also fixes \(\ell\) and \(\gamma\) in \eqref{alg:coordinate-parameters}.
These choices preserve both \(\eta\le\lambda^4/(100n^2)\) and \(\eta J\le\lambda^3/100\), which are the bounds needed for the correctness estimates.
The uniform running time estimates for these parameters are given in \Cref{subsec:numerics-complexity}.
Apply the dichotomy algorithm in \Cref{lem:triangle-test-algorithm} to \(H\), with the current weights \(w\) and the parameters \(\theta_v\) defined in \eqref{alg:theta-parameters}, and consider its two outcomes.

\subsubsection{Select a vertex}
Suppose the algorithm verifies \(T_H(w)\le C_\triangle D_\theta\).
By \Cref{lem:deletion}, the weighted sum of gains satisfies \begin{equation}\label{alg:gain-identity} \sum_xw(x)g_x =(1-c)S-6cT_H(w) -c\sum_v\zeta_v\bigl(w(v)\Lambda_v-\xi_v/2\bigr). \end{equation}
We estimate the triangle term and the derivative term (the last sum on the right-hand side) separately.
Using \(\Lambda_v=w(v)d_{*,v}-\zeta_v/2\), \(d_{*,v}\xi_v\le u(v)\), \(\norm{\zeta}_\infty\le\eta\), and \(\sum_v\xi_v\le1\), we have \begin{equation}\label{alg:budget} D_\theta=2r\sum_v\bar\xi_v\Lambda_v \le4r\sum_v\xi_v\Lambda_v \le4r(S+\eta/2). \end{equation}
Consequently, \[ 6cT_H(w)\le6cC_\triangle D_\theta\le24cC_\triangle r(S+\eta/2). \]
The derivative term satisfies \[ \abs[\bigg]{\sum_v\zeta_v(w(v)\Lambda_v-\xi_v/2)} \le\eta\biggl(\sum_vw(v)\Lambda_v+\frac{1}{2}\sum_v\xi_v\biggr) \le\eta(n^2\lambda^{-2}+1/2). \]
Our choice of \(\eta\), together with \(S\ge\lambda^2\), also gives \[ \eta(n^2\lambda^{-2}+1)\le\frac{\lambda^2}{100}+\frac{\lambda^4}{100n^2} \le\frac{\lambda^2}{50}\le\frac{S}{4}. \]
Substituting these estimates into \eqref{alg:gain-identity}, we obtain \begin{align*} \sum_xw(x)g_x &\ge(1-c-24cC_\triangle r)S -c\eta(n^2\lambda^{-2}+1/2+12C_\triangle r)\\
    &\ge\frac{1}{2}S-\eta(n^2\lambda^{-2}+1) \ge\frac{S}{4} \ge \frac{\lambda^2}{4}. \end{align*}
Here we used \(c(1+24C_\triangle r)\le1/2\), which follows from the choice of \(c\) in \eqref{eq:constant}.
Since \(\sum_xw(x)\le n\lambda^{-1}\), there exists a vertex \(x\) such that \begin{equation}\label{alg:pivot-margin} g_x\ge\frac{\lambda^3}{4n}. \end{equation}
Using \Cref{subsec:numerics-evaluation}, evaluate every gain to additive error at most \(\lambda^3/(100n)\) and select a vertex with largest approximate gain.
Its true gain is within twice this error of the largest true gain, so it is positive.

\subsubsection{Reweight two sets}
Otherwise the algorithm returns \(A,B,p,\rho\) satisfying \eqref{alg:size}, \eqref{alg:strong-witness}, and \eqref{alg:witness-exponent}.
Since \(p\le P_n\) and \(w(V(H))\le n\lambda^{-1}\), the definition in \eqref{eq:Jpw-definition} gives \(J_p(w)\le J\).
Thus \Cref{lem:improvement} gives a vertex weight \(w'\) such that \(\norm{w'-w}_1\le J\).
The witness inequality \eqref{alg:strong-witness} forces \(A\ne\emptyset\), so \eqref{alg:xi-lower} and \eqref{alg:theta-parameters} give \(b(A)\ge2r\lambda^3/3\).
The potential-increase estimate in \Cref{lem:improvement}, together with \(J_p(w)\le J\), gives \[ \mathcal U_{H,u}(w')-\mathcal U_{H,u}(w) >2^{4p}p\,b(A)-\eta J>3\lambda^3, \] where the last inequality uses \(p\ge1\), \(r\ge4\), and \(\eta J\le\lambda^3/100\) by the choice of \(\eta\) in \eqref{alg:tolerance-parameters}.

The weight \(w'\) is constructed on the full vertex set \(V(H)\) by multiplying the weights on a set \(C\) by \(1+\alpha\), setting those on a disjoint nonempty set \(B'\) to zero, and leaving all other weights unchanged.
In step~\ref{step:weighted-triangle-test}, the algorithm subsequently removes the zero-weight vertices \(B'\) from \(H\).
This removal causes no further change in the potential, because their scalar terms and all incident edge products are zero.
Thus either outcome of the dichotomy algorithm in \Cref{lem:triangle-test-algorithm} removes at least one vertex and increases \(P\) before rounding.

\subsection{The algorithm and its analysis}\label{subsec:algorithm-analysis}
The three operations in \Cref{alg:weighted-independent-set} are now justified.
We verify its rounding and initialization guarantees, then combine the progress bounds to establish the running time bound and the independent-set guarantee.

\subsubsection{Rounding and initialization}
\Cref{subsec:numerics-evaluation} gives certified rational procedures for approximating the scalar terms, derivatives, and gains; their uniform costs are bounded in \Cref{subsec:numerics-complexity}.
\Cref{subsec:numerics-rounding} chooses the fixed parameters \(W_{\max}\) and \(\sigma\) in \eqref{eq:numerics-rounding-grid}.
Here \(W_{\max}\) bounds the stored vertex weights, and \(\sigma\) is the grid spacing, with \(1/\sigma\) a positive integer.
In steps~\ref{step:weighted-full-coordinate} and~\ref{step:weighted-two-set-reweighting}, round each positive vertex weight \(w(v)\), for \(v\in V(H)\), upward to this grid, leaving zero vertex weights zero.
By \eqref{eq:numerics-rounding-loss}, each rounding loses less than \(\min\{\gamma/2,\lambda^3/2\}\) from \(\mathcal U_{H,u}\).
Thus, after rounding, step~\ref{step:weighted-full-coordinate} still increases \(\mathcal U_{H,u}\) by at least \(\gamma/2\), and step~\ref{step:weighted-two-set-reweighting} still increases it by at least \(\lambda^3/2\).
Steps~\ref{step:weighted-coordinate-discard} and~\ref{step:weighted-coordinate-select} set \(w(v)\) to the corresponding endpoint and immediately discard or select \(v\), without any rounding.
Thus all surviving weights remain on the same grid throughout the algorithm, with encoding lengths controlled by \(W_{\max}\) and \(1/\sigma\).

Start with \(I\defeq\emptyset\) after removing zero-weight vertices, and select all initial isolated vertices.
If the remaining graph \(H\) is nonempty, apply \Cref{lem:lift} to \(H\), taking \(s\) to be the prescribed input vertex weights restricted to \(V(H)\).
Thus the exact initial weights are \(q(v)\defeq\sqrt{-u(v)f'(d_{H,s}(v))}\) for \(v\in V(H)\).
Initialize \(w\) at the certified upper grid approximations to \(q\) constructed in \Cref{subsec:numerics-initialization}.
By \eqref{eq:numerics-initialization-loss}, this initialization loses less than \(1/(100n)\) from the exact potential.
Since deleting zero-weight vertices can only increase the target sum, and each selected isolate contributes \(u(v)\) in place of its target contribution \(c\, u(v)\), the initial value satisfies \[ P\ge c\Phi_{G,u}(s)-\frac{1}{100n}. \]

\subsubsection{Progress and termination}
Every selection has nonnegative gain, step~\ref{step:weighted-two-set-reweighting} increases \(P\) even after rounding, and step~\ref{step:weighted-full-coordinate} increases \(P\) by at least \(c\gamma/2\).
Steps~\ref{step:weighted-coordinate-discard} and~\ref{step:weighted-coordinate-select} incur no loss from the weight update, so their only possible loss is the discard loss already bounded in \eqref{alg:discard-loss} by \(1/(100n)\).
Consequently, after \(K\) executions of step~\ref{step:weighted-full-coordinate}, \[ P\ge \biggl( c\Phi_{G,u}(s)-\frac{1}{100n} \biggr) +\frac{Kc\gamma}{2}-\frac{1}{100n} \ge\frac{Kc\gamma}{2}-\frac{1}{50n}. \]
Thus \(P\le1\) in \eqref{alg:tracked-potential} gives \(K<4/(c\gamma)\).
Every other iteration removes a vertex, so there are at most \(n\) such iterations, including at most \(n\) calls to the dichotomy algorithm in \Cref{lem:triangle-test-algorithm}.
Consequently, the total number of iterations is less than \(4/(c\gamma)+n\).
\Cref{subsec:numerics-complexity} combines this bound with the rational implementation to prove running time \(2^{\rtconstant r^2}L^{\rtdegree}\).

\subsubsection{The independent-set guarantee}
At termination, \(H\) is empty so that \(P=u(I)\), giving \[ u(I)\ge c\Phi_{G,u}(s)-\frac{1}{50n}, \] where \(1/(50n)\) is the initial loss bound \(1/(100n)\) plus the discard loss bound \(1/(100n)\).
If \(c\Phi_{G,u}(s)\ge2/(50n)\), the preceding lower bound is at least \((c/2)\Phi_{G,u}(s)\).
Otherwise a vertex of largest objective weight has weight at least \(1/n>(c/2)\Phi_{G,u}(s)\).
Returning the heavier set proves the independent-set guarantee in \Cref{thm:algorithm}.

\subsection{Constructing the local-demand distribution}\label{subsec:construct-distribution}

We use \Cref{thm:algorithm} with varying objective weights to prove \Cref{cor:algorithm-local}.
Meeting the inclusion demands is a fractional covering problem over the independent sets; multiplicative-weights methods for such problems are surveyed by Arora, Hazan, and Kale~\cite[Section~3.3.3]{arora2012multiplicative}.
We give a rational update rule and prove its guarantee below.
Let \(n\defeq\abs{V(G)}\).
If \(n=0\), return the list consisting of the empty set, so assume \(n\ge1\).
All calls to \Cref{thm:algorithm} use the original graph \(G\) and degree weights \(s\equiv1\).

We first choose positive rational targets \(a_v\) satisfying
\begin{equation}\label{eq:rational-demands}
    \frac{c_r}{4}f(d_G(v))\le a_v\le\frac{c_r}{2}f(d_G(v)).
\end{equation}
Use \(f_N\) from \eqref{eq:numerics-f-approximation} with \(N\defeq4n\), and put
\[
    a_v\defeq\frac{c_r}{2}f_{4n}(d_G(v)).
\]
For \(0\le d\le n-1\), the integral representation of \(f\) gives \(f(d)\ge\int_0^1\frac{1-x}{n}\,dx=1/(2n)\).
Together with \eqref{eq:numerics-f-approximation-error}, this gives
\[
    \frac{1}{2}f(d_G(v))\le f_{4n}(d_G(v))\le f(d_G(v))\le1,
    \qquad \frac{c_r}{8n}\le a_v\le1.
\]
Thus \eqref{eq:rational-demands} holds, and the targets are computable exactly with the encoding lengths bounded in \Cref{subsec:numerics-complexity}.

Recall \(b_n=\lceil\log_2(2n)\rceil\), set
\begin{equation}\label{eq:rational-mw-parameters}
    m\defeq32\lceil8n/c_r\rceil b_n,
\end{equation}
and initialize \(z_v^{(0)}\defeq1\) for every \(v\in V(G)\).
For each \(t\in\{0,\ldots,m-1\}\), call \Cref{thm:algorithm} with objective weights \(u(v)\defeq z_v^{(t)}\), obtain an independent set \(I_{t+1}\), and update
\begin{equation}\label{eq:rational-mw-update}
    z_v^{(t+1)}\defeq z_v^{(t)}
    \biggl(1+\frac{a_v-\ind_{\{v\in I_{t+1}\}}}{4}\biggr)
    \qquad(v\in V(G)).
\end{equation}
Each update factor is either \(1+a_v/4\) or \(3/4+a_v/4\), hence is a positive rational at least \(3/4\).
All objective weights therefore remain positive rationals.
By \Cref{thm:algorithm} and \eqref{eq:rational-demands}, each returned set satisfies
\[
    \sum_{v\in I_{t+1}}z_v^{(t)}
    \ge\frac{c_r}{2}\sum_vz_v^{(t)}f(d_G(v))
    \ge\sum_va_vz_v^{(t)}.
\]
Consequently, the total weight \(Z_t\defeq\sum_vz_v^{(t)}\) satisfies
\[
    Z_{t+1}
    =Z_t+\frac{1}{4}\biggl(\sum_va_vz_v^{(t)}
    -\sum_{v\in I_{t+1}}z_v^{(t)}\biggr)
    \le Z_t\le Z_0=n.
\]

To bound the inclusion frequencies, use \(\log(1+x)\ge x-x^2\) for \(\abs{x}\le1/4\) and \((a-h)^2\le a+h\) for \(0\le a\le1\) and \(h\in\{0,1\}\), which give
\[
    \log\biggl(1+\frac{a-h}{4}\biggr)
    \ge\frac{a-h}{4}-\frac{(a-h)^2}{16}
    \ge\frac{3a-5h}{16}.
\]
For each vertex \(v\), let \(N_v\defeq\abs{\{j\in\{1,\ldots,m\}:v\in I_j\}}\).
Summing the logarithmic increments in \eqref{eq:rational-mw-update} yields
\[
    \log n\ge\log z_v^{(m)}\ge\frac{3ma_v-5N_v}{16}.
\]
The choice of \(m\) and the lower bound \(a_v\ge c_r/(8n)\) give \(ma_v\ge32b_n\ge32\log n\), so
\[
    \frac{N_v}{m}
    \ge\frac{3a_v}{5}-\frac{16\log n}{5m}
    \ge\frac{a_v}{2}\ge\frac{c_r}{8}f(d_G(v)).
\]
Thus the uniform distribution on the returned list \(I_1,\ldots,I_m\) has the required inclusion probabilities.
The degree comparison in the deduction of \Cref{thm:main} gives \(f(d)\ge\eta_0\log(d+2)/(d+1)\) for an absolute constant \(\eta_0>0\).
Together with \(c_r\ge\kappa_1/r\), this proves the final assertion with \(\kappa_{\rm alg}\defeq\kappa_1\eta_0/8\).

The construction makes \(m=O(rn\log(2n))\) calls to \Cref{thm:algorithm}.
\Cref{subsec:numerics-complexity} bounds their input lengths by a polynomial of absolute degree in \(n,r\) and proves total running time \(2^{\rtconstant r^2}n^{\rtdegree}\).


\section{Neighborhood walks and entropy bounds}\label{sec:triangle-walks}
This section and the next prove the structural dichotomy in \Cref{lem:triangle-test} by iterating the reduction step in \Cref{lem:triangle-reduction}.

Recall that a \emph{vertex splitting} is an operation which partitions the incident edges at a vertex into groups and replaces the vertex by one copy for each group.
Each copy inherits its parent's weight and parameter.
Choose the absolute constants \[ %
 C_{\mathrm{ent}}\defeq20+12\cdot24^2,\qquad C_{\mathrm{loss}}\defeq24(C_{\mathrm{ent}}+1),\qquad t_*\defeq(100C_{\mathrm{loss}})^4,\qquad C_\triangle\defeq2^{t_*}. %
\]
For an integer \(t\ge32\), define the neighborhood-ratio bound \(\ratiobound{t}\) and the rational loss allowance \(\delta(t)\) by \begin{equation}\label{eq:triangle-rational-loss} \ratiobound{t}\defeq2^t, \qquad \delta(t)\defeq\frac{2C_{\mathrm{loss}}\lceil\log_2t\rceil}{\lfloor\sqrt t\rfloor}. \end{equation}
The allowance \(\delta(t)\) bounds the analytic loss estimate by a rational number, so the eventual split can be checked using exact arithmetic.
Write \(T_{\rm in},D_{\rm in}\) for the triangle weight and \(D_\theta\) before a reduction, and \(T_{\rm out},D_{\rm out}\) for their values afterwards.

\begin{lemma}[Reduction step]\label{lem:triangle-reduction}
Suppose \(t>t_*\) is an integer, \(\Lambda_v/\theta_v\le2^t\) at every vertex of the current graph, and \(T_{\rm in}>2D_{\rm in}>0\).
Either there are sets and parameters satisfying \eqref{alg:size} and \eqref{alg:strong-witness}, with \(p\le30t^3\) and \(\rho\) a current vertex parameter, or deletion and splitting produce a graph in which \begin{equation}\label{eq:one-reduction} \begin{gathered} t'\defeq\lceil\sqrt t\rceil,\qquad \Lambda'_v/\theta'_v\le2^{t'}\quad\text{for every remaining vertex},\\ D_{\rm out}\le D_{\rm in},\qquad T_{\rm out}\ge(1-\delta(t))T_{\rm in}-D_{\rm in}/t. \end{gathered} \end{equation}
\end{lemma}

This section prepares the graph and either finds witness sets or establishes the walk estimates in \Cref{lem:triangle-walk-estimates}.
\Cref{sec:triangle-proof} uses these estimates to complete the proof of \Cref{lem:triangle-reduction}, then iterates the reduction to prove \Cref{lem:triangle-test}.
The splitting groups neighbors by walk endpoints chosen using a label coupling: close endpoint distributions allow labels to agree at most triangle corners, preserving triangles.
To keep each group small, we retain a neighbor only when the transition probability from its label back to that neighbor is sufficiently large relative to its weight.
Summing these reverse transition probabilities bounds the retained group weight, while the entropy bound controls triangle loss from rejected edges.
Together, these estimates allow us to reduce the ratios \(w(N(v))/\theta_v\) with controlled triangle loss.
Once these ratios are bounded by an absolute constant, an elementary count bounds the remaining triangle weight by a constant multiple of \(D_\theta\), and the controlled losses give the same conclusion for the original graph.

\subsection{Edge deletion}\label{subsec:edge-deletion}
Starting from the input graph, we repeatedly delete edges and split vertices.
Each copy retains the weight and parameter of its original vertex.
Every surviving edge or triangle corresponds to a distinct original edge or triangle of the same weight.
In particular, neither triangle weight nor \(D_\theta=\sum_{xy\in E}w(x)w(y)(\theta_x+\theta_y)\) can increase under deletion and splitting.

Suppose the current graph satisfies \(\Lambda_v/\theta_v\le2^t\) at every vertex, where \(t\ge32\) is an integer.
Recall that \(T_{\rm in},D_{\rm in}\) denote its triangle weight and \(D_\theta\).
The iteration in \Cref{subsec:triangle-iteration} stops whenever \(T_{\rm in}\le2D_{\rm in}\), so we only need to perform the reduction when \(T_{\rm in}>2D_{\rm in}\).
In this case the graph contains a triangle, so \(D_{\rm in}>0\).
Deleting edges with too little common-neighbor weight costs little triangle weight and bounds the transition probabilities relative to vertex weights.

Repeatedly remove an edge \(xv\) if its current common-neighbor weight is less than \(\max(\theta_x,\theta_v)/t\).
Removing it destroys triangles of total weight less than \(w(x)w(v)\max(\theta_x,\theta_v)/t\).
Since each edge is removed at most once and \(\max(\theta_x,\theta_v)\le\theta_x+\theta_v\), the total triangle loss is at most \(D_{\rm in}/t\).
Discard isolated vertices, and call the resulting graph \(F\).
Until we complete one reduction in \Cref{sec:triangle-proof}, all neighborhoods and local quantities refer to \(F\), and we write \(T\defeq T_F(w)\), \(D\defeq D_\theta(F)\).
We have \begin{equation}\label{eq:triangle-preparation} D\le D_{\rm in},\qquad T\ge T_{\rm in}-D_{\rm in}/t>D_{\rm in}\ge D>0. \end{equation}
Here \(T>0\) ensures that \(F\) has an edge and hence \(D>0\).
Every remaining edge satisfies \begin{equation}\label{eq:triangle-common-neighbors} \Lambda_x\le \ratiobound{t}\theta_x,\qquad C_{xv}\defeq w(N(x)\cap N(v))\ge\max(\theta_x,\theta_v)/t. \end{equation}
Consequently, \begin{equation}\label{eq:triangle-adjacent-scales} \Lambda_x\ge\theta_x/t,\qquad (t\ratiobound{t})^{-1}\le\theta_x/\theta_v\le t\ratiobound{t}\quad(xv\in E(F)). \end{equation}
The ratio bound follows from \(\theta_v/t\le C_{xv}\le\Lambda_x\le \ratiobound{t}\theta_x\) and the same inequalities with \(x,v\) interchanged.

\subsection{Neighborhood walks}\label{subsec:triangle-neighborhood-walks}
For each center vertex \(x\), define \(P_x\) to be the random walk on \(N(x)\) that moves from \(v\) to a neighbor \(z\) within \(F[N(x)]\), chosen with probability proportional to \(w(z)\).
Write \(P_x\) also for its transition matrix, and define its lazy version \(Q_x\) and the kernels \(K_x^j\) by \[ P_x(v,z)\defeq\frac{w(z)\ind_{vz\in E(F)}}{C_{xv}},\qquad Q_x\defeq(I+P_x)/2,\qquad K_x^j\defeq P_xQ_x^j. \]
All denominators are positive by \eqref{eq:triangle-common-neighbors}.
Here \(I\) is the identity matrix, \(j\ge0\) is an integer, and distributions are row vectors multiplied by transition matrices on the right.
Thus \(Q_x\) stays at its current state with probability \(1/2\) and otherwise makes a \(P_x\)-step.
This laziness will let us compare endpoint distributions from adjacent starting vertices through a common mixture.
We write \(K(v,\cdot)\) for a row and \(K(v,A)\defeq\sum_{y\in A}K(v,y)\) for a vertex set \(A\).

To connect these walks with triangle weight, put \[ T_x\defeq\sum_{v\in N(x)}w(v)C_{xv},\qquad \pi_x(v)\defeq\frac{w(v)C_{xv}}{T_x},\qquad \nu(x,v)\defeq\frac{w(x)w(v)C_{xv}}{6T}. \]
Since \(\sum_xw(x)T_x=6T\), the distribution \(\nu\) is a probability distribution on directed edges, with conditional distribution \(\pi_x\) for \(v\) given \(x\).
Sampling \((x,v)\) from \(\nu\) and then \(z\) from \(P_x(v,\cdot)\) gives the ordered-triangle distribution \begin{equation}\label{eq:triangle-sampling} \tau(x,v,z)\defeq\frac{w(x)w(v)w(z)}{6T}. \end{equation}
Thus one walk step completes a weighted random triangle.
Comparing endpoint distributions at \(v\) and \(z\) will therefore control label disagreements at the corners of the triangles we wish to preserve.

The identity \[ \pi_x(v)P_x(v,z)=\frac{w(v)w(z)\ind_{vz\in E(F)}}{T_x} =\pi_x(z)P_x(z,v) \] shows that \(P_x\) is reversible with respect to \(\pi_x\).
Summing over \(v\) gives \(\pi_xP_x=\pi_x\), so \(\pi_x\) is stationary: one step leaves this distribution unchanged.
Thus the same properties hold for \(Q_x\) and \(K_x^j\).
We will use stationarity to average entropy increments and reversibility to analyze the retention test.

The initial \(P_x\)-step gives a density bound that persists through the subsequent lazy steps.
Indeed, \eqref{eq:triangle-common-neighbors} gives \(P_x(z,y)\le tw(y)/\theta_x\).
Since \(P_x\) and \(Q_x\) commute, \(K_x^j=Q_x^jP_x\), so each row of \(K_x^j\) is a mixture of \(P_x\)-rows and \begin{equation}\label{eq:walk-density} K_x^j(v,y)\le tw(y)/\theta_x. \end{equation}
We next bound the averaged entropy of these endpoint distributions, both to find a common walk length with close distributions and to control triangle loss from rejected edges.

\subsection{Sets capturing the walk endpoints}\label{subsec:triangle-capture}
To bound the averaged entropy, we first find sets containing most of the endpoint probability mass; the next subsection will bound the weight of their intersections with \(N(x)\), divided by \(\theta_x\), for most pairs \((x,v)\).
For each starting vertex \(v\), we construct one such set that works simultaneously for every center \(x\in N(v)\).
The following enlargement rule makes the probability escaping at each step small, while keeping the sets within the size bounds required for a witness.
For each starting vertex \(v\), define \begin{align*} S_0(v)&\defeq N(v),\\
    S_{i+1}(v)&\defeq S_i(v)\cup \{y:w(N(y)\cap S_i(v))\ge\theta_v/(t^7\ratiobound{t}^2)\}, \qquad 0\le i<t. \end{align*}

Every vertex of \(S_i(v)\) has distance at most \(i+1\) from \(v\), so its \(\theta\)-value is at most \(\theta_v(t\ratiobound{t})^{i+1}\) by \eqref{eq:triangle-adjacent-scales}.
Each new vertex has at least \(\theta_v/(t^7\ratiobound{t}^2)\) neighbor weight in \(S_i(v)\), so double counting gives \[ w(S_{i+1}(v)) \le\bigl(1+t^7\ratiobound{t}^3(t\ratiobound{t})^{i+1}\bigr)w(S_i(v)). \]
Starting with \(w(S_0(v))\le \ratiobound{t}\theta_v\), we obtain \begin{equation}\label{eq:triangle-capture-size} w(S_t(v))\le\theta_v2^{30t^3},\qquad b(S_t(v))\le\theta_v^22^{60t^3}. \end{equation}

We next prove the promised concentration bound: \begin{equation}\label{eq:walk-leakage} K_x^i(v,N(x)\setminus S_i(v))\le it^{-4} \qquad(x\in N(v),\ 0\le i\le t). \end{equation}
At \(i=0\), the first step ends in \(N(v)=S_0(v)\) with probability one.
For the inductive step, use \(K_x^{i+1}=K_x^iQ_x\).
Mass already outside \(S_i(v)\) contributes at most \(it^{-4}\).
If the walk is in \(S_i(v)\) and stays at its current vertex during the next lazy step, it remains in \(S_i(v)\subseteq S_{i+1}(v)\).
For \(y\notin S_{i+1}(v)\), the probability of moving into \(y\) from inside \(S_i(v)\) is at most \[ \sum_{z\in S_i(v)\cap N(x)}K_x^i(v,z)Q_x(z,y) \le\frac{t^2w(y)}{\theta_x^2}w(N(y)\cap S_i(v)) <\frac{t^2w(y)\theta_v}{\theta_x^2t^7\ratiobound{t}^2}. \]
The first inequality uses \eqref{eq:walk-density} and \(Q_x(z,y)=P_x(z,y)/2\) for \(z\ne y\), together with \(C_{xz}\ge\theta_x/t\) from \eqref{eq:triangle-common-neighbors}; the second inequality follows from the definition of \(S_{i+1}(v)\) since \(y\notin S_{i+1}(v)\).
Summing over \(y\in N(x)\setminus S_{i+1}(v)\), and using \(\Lambda_x\le \ratiobound{t}\theta_x\) and \(\theta_v/\theta_x\le t\ratiobound{t}\), bounds this new escaping mass by \(t^{-4}\).
Adding the two contributions proves \eqref{eq:walk-leakage}.

\subsection{Two tests and exceptional pairs}\label{subsec:triangle-tests}
The concentration bound \eqref{eq:walk-leakage} alone does not give the entropy bound we need, since \(w(N(x)\cap S_t(v))/\theta_x\) might still be large.
We introduce two tests.
If either fails, it supplies a witness satisfying alternative~\ref{case:triangle-test-witness} of \Cref{lem:triangle-test}.
If both hold at every vertex, the pairs \((x,v)\) for which this ratio is large have small probability under \(\nu\).
For each \(v\), define the set of such centers by \[ B_v\defeq\{x\in N(v):w(N(x)\cap S_t(v))>t^{10}\theta_x\}. \]
We test the following inequalities for every vertex \(v\): \begin{align} e_w(N(v),S_t(v)) &\le24\bigl(30t^3b(N(v))+\theta_v^22^{-30t^3}\bigr), \label{alg:test-one}\\
    e_w(B_v,N(v)) &\le24\bigl(2tb(B_v)+\theta_v^22^{-2t}\bigr). \label{alg:test-two} \end{align}
When both tests hold, the first bounds \(b(B_v)\), and the second converts this into a bound on the probability of encountering a pair with \(x\in B_v\) under \(\nu\).

First we check that failure of either test supplies an admissible witness.
For \eqref{alg:test-one}, take \((A,B,p,\rho)\defeq(N(v),S_t(v),30t^3,\theta_v)\).
The size conditions \eqref{alg:size} follow from \eqref{eq:triangle-capture-size} and \(w(N(v))\le \ratiobound{t}\theta_v\).
For \eqref{alg:test-two}, take \((A,B,p,\rho)\defeq(B_v,N(v),2t,\theta_v)\) and use \[ w(B_v),w(N(v))\le \ratiobound{t}\theta_v,\qquad b(N(v))\le t\ratiobound{t}^2\theta_v^2\le2^{4t}\theta_v^2, \] where the bound on \(b(N(v))\) follows from \(\theta_x\le t\ratiobound{t}\theta_v\) for \(x\in N(v)\).
Failure of either test is exactly the strict inequality \eqref{alg:strong-witness} for these choices.

These witnesses can also be returned when the current graph already contains split vertices.
Let \(G\) be the original input graph, let \(G'\) be the current graph, and let \(\sigma\) map every current vertex to its original vertex.
For \(A',B'\subseteq V(G')\), put \(A\defeq\sigma(A')\), \(B\defeq\sigma(B')\).
Then we have \begin{equation}\label{eq:witness-projection} \begin{gathered} e_{G,w}(A,B)\ge e_{G',w}(A',B'),\qquad b_G(A)\le b_{G'}(A'),\qquad b_G(B)\le b_{G'}(B'),\\ w_G(A)\le w_{G'}(A'),\qquad w_G(B)\le w_{G'}(B'). \end{gathered} \end{equation}
Thus a witness in \(G'\) gives one in \(G\) with the same \(p,\rho\), where \(\rho\) is an original \(\theta\)-value by inheritance.

Assume from now on that both tests hold for every \(v\), and call \((x,v)\) exceptional if \(x\in B_v\).
Since \(B_v\subseteq N(v)\), we have \[ e_w(B_v,N(v)) =\sum_{x\in B_v}w(x)w(N(x)\cap N(v)) =\sum_{x\in B_v}w(x)C_{xv}. \]
To compute the probability of an exceptional pair, sum \(\nu(x,v)\) over all directed edges with \(x\in B_v\).
Using the definition of \(\nu\) and the preceding identity gives \begin{equation}\label{eq:exceptional-probability-formula} \begin{aligned} \Prob_{(x,v)\sim\nu}(x\in B_v) =\sum_v\sum_{x\in B_v}\nu(x,v) =\frac{1}{6T}\sum_vw(v)\sum_{x\in B_v}w(x)C_{xv} =\frac{1}{6T}\sum_vw(v)e_w(B_v,N(v)). \end{aligned} \end{equation}
Thus it suffices to bound \(e_w(B_v,N(v))\) for each \(v\).
Multiplying \(w(N(x)\cap S_t(v)) > t^{10}\theta_x\) for each \(x\in B_v\) by \(w(x)\) and summing, then applying \eqref{alg:test-one}, gives \[ t^{10}b(B_v)\le e_w(N(v),S_t(v)) \le24\bigl(30t^3b(N(v))+\theta_v^22^{-30t^3}\bigr). \]
Substitution into \eqref{alg:test-two} yields \begin{equation}\label{eq:exceptional-edge-weight} e_w(B_v,N(v))\le60\cdot24^2t^{-6}b(N(v)) +(2\cdot24^2t^{-9}2^{-30t^3}+24\cdot2^{-2t})\theta_v^2. \end{equation}
The definitions of \(b\) and \(D\) give \(\sum_vw(v)b(N(v))=D\).
Also, \eqref{eq:triangle-adjacent-scales} gives \(\theta_v\le t\Lambda_v\), so \[ \sum_vw(v)\theta_v^2 \le t\sum_vw(v)\theta_v\Lambda_v=tD. \]
Substituting \eqref{eq:exceptional-edge-weight} into \eqref{eq:exceptional-probability-formula} and using these two sums together with \(T\ge D\) and \(t\ge32\), we obtain \begin{equation}\label{eq:triangle-exceptional} \begin{aligned} \Prob_\nu(x\in B_v) &\le10\cdot24^2t^{-6} +(24^2/3)t^{-8}2^{-30t^3}+(24/6)t2^{-2t}\\
    &\le12\cdot24^2t^{-6}. \end{aligned} \end{equation}

\subsection{Entropy bounds and a common walk length}\label{subsec:triangle-entropy}
The entropy bound has two uses: its small total increase supplies a common walk length with close endpoint distributions, and its upper bound will later control the probability that an edge is deleted by the retention test \eqref{eq:incidence-test}.
The use of entropy increments to compare walk distributions has precedents in Erschler and Karlsson~\cite[Lemma~5.1]{erschler2010homomorphisms} and Benjamini, Duminil-Copin, Kozma, and Yadin~\cite[Section~2, Theorem~8]{benjamini2015disorder}.

For a probability vector \(p\) on a finite set and positive weights \(s(y)\), define \[ H_s(p)\defeq\sum_y p(y)\log\frac{s(y)}{p(y)}, \] with zero terms interpreted as zero.
The weights \(s(y)\) need not sum to one, so this entropy may be negative.
For probability vectors \(p,q\) on the same finite set, define their relative entropy by \(\KL(p\|q)\defeq\sum_y p(y)\log(p(y)/q(y))\), which is nonnegative.
Terms with \(p(y)=0\) are zero, and \(\KL(p\|q)=+\infty\) if \(p(y)>0=q(y)\) for some \(y\).
For probability vectors \(p_i\), coefficients \(\lambda_i\ge0\) with \(\sum_i\lambda_i=1\), and \(m\defeq\sum_i\lambda_ip_i\), expanding the definitions gives the mixture identity \begin{equation}\label{eq:mixture} H_s(m)-\sum_i\lambda_iH_s(p_i) =\sum_i\lambda_i\KL(p_i\|m). \end{equation}
Terms with \(\lambda_i=0\) are omitted.

For \(s_x(y)\defeq w(y)/\theta_x\) on \(N(x)\), define the averaged entropy of the walk endpoint distributions by \[ H_j\defeq\E_{(x,v)\sim\nu}H_{s_x}\bigl(K_x^j(v,\cdot)\bigr). \]
We use \(\norm{p-q}_{\TV}\defeq\tfrac{1}{2}\sum_y\abs{p(y)-q(y)}\) for total variation distance.

\begin{lemma}\label{lem:triangle-walk-estimates}
If both \eqref{alg:test-one} and \eqref{alg:test-two} hold for every vertex of the graph \(F\) obtained after edge deletion, then \begin{equation}\label{eq:triangle-entropy-bounds} -\log t\le H_0\le H_1\le\cdots\le H_t\le C_{\mathrm{ent}}\log t. \end{equation}
Moreover, some common index \(0\le j<t\) satisfies \begin{equation}\label{eq:nearby-rows} \E_{(x,v,z)\sim\tau} \norm{K_x^j(v,\cdot)-K_x^j(z,\cdot)}_{\TV} \le2\sqrt{(C_{\mathrm{ent}}+1)\log t/t}. \end{equation}
\end{lemma}
\begin{proof}
By \eqref{eq:walk-density}, each ratio \(s_x(y)/K_x^j(v,y)\) with positive denominator is at least \(1/t\), so \(H_j\ge-\log t\).
To bound \(H_t\) from above, first fix a nonexceptional pair \((x,v)\), so \(x\notin B_v\).
By \eqref{eq:walk-leakage} with \(i=t\), the endpoint distribution \(K_x^t(v,\cdot)\) assigns probability at least \(1-t^{-3}\) to \(N(x)\cap S_t(v)\).
By the definition of \(B_v\), we have \[ s_x(N(x)\cap S_t(v))=w(N(x)\cap S_t(v))/\theta_x\le t^{10}. \]
Also, \(s_x(N(x))=\Lambda_x/\theta_x\le \ratiobound{t}\).
For a general probability vector \(p\) and positive weights \(s\) on the same finite set, concavity of the logarithm gives, for any subset \(A\) with \(p(A)>0\), \[ \sum_{y\in A}p(y)\log\frac{s(y)}{p(y)}\le p(A)\log\frac{s(A)}{p(A)}; \] if \(p(A)=0\), both sides are interpreted as zero.
We apply this inequality with \(p(y)=K_x^t(v,y)\) and \(s(y)=s_x(y)\) on \(N(x)\).
Put \(A\defeq N(x)\cap S_t(v)\) and \(\mu\defeq K_x^t(v,A)\), so the probabilities of \(A\) and \(N(x)\setminus A\) are \(\mu\) and \(1-\mu\le t^{-3}\), respectively.
Apply the preceding inequality separately to these two sets, using \(s_x(A)\le t^{10}\) and \(s_x(N(x)\setminus A)\le \ratiobound{t}\), to obtain \[ \begin{aligned} H_{s_x}\bigl(K_x^t(v,\cdot)\bigr) &\le\mu\log\frac{t^{10}}{\mu}+(1-\mu)\log\frac{\ratiobound{t}}{1-\mu}\\
    &= -\mu\log\mu-(1-\mu)\log(1-\mu)+10\mu\log t+(1-\mu)\log \ratiobound{t}\\
    &\le \log 2 + 10\log t + t^{-2}. \end{aligned} \]
For exceptional pairs, the entropy is still at most \(\log s_x(N(x))\le \log \ratiobound{t}\le t\).
By \eqref{eq:triangle-exceptional}, their probability under \(\nu\) is at most \(12\cdot24^2t^{-6}\), so their contribution to \(H_t\) is at most \(12\cdot24^2t^{-5}\).
Averaging over both types of pairs therefore gives \[ H_t\le\log2+10\log t+t^{-2}+12\cdot24^2t^{-5}\le C_{\mathrm{ent}}\log t. \]

We now prove that \(H_j\) is nondecreasing.
Laziness of \(Q_x\) provides a common mixture containing the rows we wish to compare: \[ K_x^{j+1}(v,\cdot) =\frac{1}{2}K_x^j(v,\cdot) +\frac{1}{2}\sum_zP_x(v,z)K_x^j(z,\cdot). \]
For fixed \(x,v\), apply \eqref{eq:mixture} to the distributions \(K_x^j(z,\cdot)\) with mixture weights \(Q_x(v,z)\), obtaining \[ H_{s_x}\bigl(K_x^{j+1}(v,\cdot)\bigr) -\sum_zQ_x(v,z)H_{s_x}\bigl(K_x^j(z,\cdot)\bigr) =\sum_zQ_x(v,z)\KL\bigl(K_x^j(z,\cdot)\,\|\,K_x^{j+1}(v,\cdot)\bigr). \]
Now average over \((x,v)\sim\nu\).
The first term on the left averages to \(H_{j+1}\) by definition.
Conditional on \(x\), we have \(v\sim\pi_x\), and stationarity gives \(\sum_v\pi_x(v)Q_x(v,z)=\pi_x(z)\).
Thus sampling \(z\sim Q_x(v,\cdot)\) leaves the conditional distribution unchanged, so \((x,z)\) also has distribution \(\nu\).
The second term on the left therefore averages to \(H_j\), giving \begin{equation}\label{eq:entropy-increment} H_{j+1}-H_j =\E_{(x,v)\sim\nu}\E_{z\sim Q_x(v,\cdot)} \KL\bigl(K_x^j(z,\cdot)\,\|\,K_x^{j+1}(v,\cdot)\bigr) \ge0. \end{equation}
This proves monotonicity and hence \eqref{eq:triangle-entropy-bounds}.

Since \(H_t-H_0\le(C_{\mathrm{ent}}+1)\log t\), choose \(j<t\) with \(h\defeq H_{j+1}-H_j\le(C_{\mathrm{ent}}+1)\log t/t\).
Expanding \(Q_x=(I+P_x)/2\) in \eqref{eq:entropy-increment} gives \[ h =\frac{1}{2} \biggl( \E_{(x,v)\sim\nu} \KL\bigl(K_x^j(v,\cdot)\,\|\,K_x^{j+1}(v,\cdot)\bigr) + \E_{(x,v)\sim\nu}\E_{z\sim P_x(v,\cdot)} \KL\bigl(K_x^j(z,\cdot)\,\|\,K_x^{j+1}(v,\cdot)\bigr) \biggr). \]
Both expectations are nonnegative, so each is at most \(2h\).
In the second expectation, \((x,v,z)\) has distribution \(\tau\), by \eqref{eq:triangle-sampling}.
Pinsker's inequality \(\norm{p-q}_{\TV}\le\sqrt{\KL(p\|q)/2}\) and \(\E\sqrt Z\le\sqrt{\E Z}\) therefore give \[ \begin{aligned} \E_{(x,v)\sim\nu}\norm{K_x^j(v,\cdot)-K_x^{j+1}(v,\cdot)}_{\TV} &\le\sqrt h,\\
    \E_{(x,v,z)\sim\tau}\norm{K_x^j(z,\cdot)-K_x^{j+1}(v,\cdot)}_{\TV} &\le\sqrt h. \end{aligned} \]
Since the \((x,v)\)-marginal of \(\tau\) is \(\nu\), the first expectation can also be taken over \(\tau\).
The triangle inequality through \(K_x^{j+1}(v,\cdot)\) then bounds \(\E_{(x,v,z)\sim\tau} \norm{K_x^j(v,\cdot)-K_x^j(z,\cdot)}_{\TV}\), the left side of \eqref{eq:nearby-rows}, by \(2\sqrt h\), as required.
A proof of the version of Pinsker's inequality used here appears in \Cref{app:probability}.
\end{proof}

The same index \(j\) is used at every center \(x\), and \(H_j\le H_t\) retains the entropy bound needed to control rejection by the test \eqref{eq:incidence-test}.
For rational input, we can choose \(j\) without evaluating entropies: compute the average on the left of \eqref{eq:nearby-rows} for every \(0\le j<t\), and minimize it.
This choice still satisfies both \eqref{eq:nearby-rows} and \(H_j\le C_{\mathrm{ent}}\log t\).


\section{Vertex splitting and the triangle dichotomy}\label{sec:triangle-proof}
We now use the walk estimates from \Cref{sec:triangle-walks} to complete the proof of \Cref{lem:triangle-reduction}, then iterate the reduction to prove \Cref{lem:triangle-test}.

The next two subsections complete the proof of \Cref{lem:triangle-reduction}.
Consider a reduction with parameter \(t>t_*\), and write \(T_{\rm in},D_{\rm in}\) for the triangle weight and \(D_\theta\) of its input graph.
First perform the edge deletions of \Cref{subsec:edge-deletion}, obtaining the prepared graph \(F\) and \eqref{eq:triangle-preparation}.
If a test fails, \Cref{subsec:triangle-tests} gives the witness alternative, including its transfer to the original input graph.
Otherwise fix an index \(j\) satisfying \Cref{lem:triangle-walk-estimates}, abbreviate \(K_x^j\) to \(K_x\), and put \(t'\defeq\lceil\sqrt t\rceil\).
Throughout these two subsections, recall that \(T=T_F(w)\) and \(D=D_\theta(F)\), and take all neighborhoods in \(F\).
Our goal is to split \(F\) so that every new vertex has neighborhood weight at most \(2^{t'}\) times its parameter, while preserving most of the triangle weight.

\subsection{Label coupling and vertex splitting}\label{subsec:finite-coupling}
At each center \(x\), assign a random label \(Y_{xv}\in N(x)\) to each \(v\in N(x)\), with distribution \(K_x(v,\cdot)\).
A \emph{label coupling} is a joint choice of random labels that preserves the specified distribution of each label.
We choose these labels jointly so that close distributions are likely to give the same label.
We use the shared rejection sampling construction of Kleinberg and Tardos~\cite[Section~3]{kleinberg2002approximation}, with the disagreement estimate of Angel and Spinka~\cite[Section~2.1]{angel2021pairwise}.

\begin{lemma}\label{lem:triangle-coupling}
For a finite family of probability vectors \(p_i\) on a common finite set \(\mathcal{L}\), there is a joint choice of labels \(Y_i\) with these marginal distributions such that \[ \Prob(Y_i\ne Y_k)\le2\norm{p_i-p_k}_{\TV} \qquad\text{for all }i,k. \]
\end{lemma}
We prove this lemma in \Cref{subsec:coupling-proof} and give a finite implementation for deterministic splitting in \Cref{subsec:choosing-labels}.

\step{Constructing the split.}
At each center \(x\), apply \Cref{lem:triangle-coupling} with \(\mathcal{L}\defeq N(x)\) and \(p_v\defeq K_x(v,\cdot)\), independently at different centers, yielding labels \(Y_{xv}\in N(x)\) for \(v\in N(x)\).
For each \(y\in N(x)\), the \emph{label group} with label \(y\) at \(x\) is the set \(\{v\in N(x):Y_{xv}=y\}\).
The label coupling controls agreement at triangle corners, but a label group could still have large weight.
The following retention test ensures that the retained neighbors with any fixed label \(y\) have total weight at most \(2^{t'}\theta_x\).
An incidence \((x,v)\) means the edge \(xv\) viewed at its endpoint \(x\).
If \(Y_{xv}=y\), retain this incidence precisely when \begin{equation}\label{eq:incidence-test} K_x(y,v)\ge(w(v)/\theta_x)2^{-t'}. \end{equation}
We say \((x,v)\) is \emph{rejected} if it fails the test.
Let \(R_{x,y}\) be the subset of the label group with label \(y\) at \(x\) consisting of neighbors whose incidences at \(x\) are retained.
We obtain the promised vertex splitting of \(F\) by replacing each vertex \(x\in V(F)\) with one copy \((x,y)\) for each \(y\) with \(R_{x,y}\ne\emptyset\).
Each copy inherits weight \(w(x)\) and parameter \(\theta_x\).
An edge \(xv\) survives if both incidences are retained; it joins \((x,Y_{xv})\) and \((v,Y_{vx})\).
Discard isolated copies.
Let \(T'\) be the triangle weight of the resulting random graph, and write \(D_{\rm out}\) for its \(D_\theta\) value, as in \Cref{lem:triangle-reduction}.
Once a label assignment is chosen, \(T_{\rm out}\) denotes the corresponding value of \(T'\).

Summing \eqref{eq:incidence-test} over \(v\in R_{x,y}\) gives \begin{equation}\label{eq:split-group-weight} w(R_{x,y}) \le2^{t'}\theta_x K_x(y,R_{x,y}) \le2^{t'}\theta_x. \end{equation}
Deleting an edge at its other endpoint can only reduce a copy's neighborhood weight, so \eqref{eq:split-group-weight} implies that every new vertex satisfies \(\Lambda'_v/\theta'_v\le2^{t'}\).
Every surviving edge comes from a distinct edge of \(F\), and hence \(D_{\rm out}\le D\).

\subsection{Triangle loss and completion of the reduction}\label{subsec:triangle-split-loss}
A triangle survives if its two incidences receive the same label at each vertex and all six incidences pass \eqref{eq:incidence-test}.
Thus we bound the loss in two parts: corner disagreements and incidence rejections; we then combine them by a union bound.

\step{Corner disagreements.}
Fix an ordered triangle \((x,v,z)\) of \(F\).
To preserve this triangle after splitting, \(v\) and \(z\) must receive the same label in the label coupling at \(x\).
We first bound the probability that \(Y_{xv}\ne Y_{xz}\).
The labels \(Y_{xv}\) and \(Y_{xz}\) have distributions \(K_x(v,\cdot)\) and \(K_x(z,\cdot)\), respectively.
By \Cref{lem:triangle-coupling}, their label coupling gives \begin{equation}\label{eq:fixed-triangle-label-disagreement} \Prob_{\mathrm{labels}}(Y_{xv}\ne Y_{xz}) \le2\norm{K_x(v,\cdot)-K_x(z,\cdot)}_{\TV}. \end{equation}
Now sample the ordered triangle with probability \(\tau(x,v,z)=w(x)w(v)w(z)/(6T)\), independently of the labels.
Averaging \eqref{eq:fixed-triangle-label-disagreement} over this choice and applying \eqref{eq:nearby-rows} gives \begin{equation}\label{eq:corner-disagreement} \E_{(x,v,z)\sim\tau}\Prob_{\mathrm{labels}}(Y_{xv}\ne Y_{xz}) \le4\sqrt{(C_{\mathrm{ent}}+1)\log t/t}. \end{equation}
Here the inner probability concerns the random labels for a fixed triangle, and the outer expectation averages over triangles.
Since \(\tau\) is unchanged by permuting the three vertices, the same bound applies to \(Y_{vx}\ne Y_{vz}\) at \(v\) and \(Y_{zx}\ne Y_{zv}\) at \(z\).

\step{Incidence rejections.}
We use the upper bound on \(H_j\) from \Cref{lem:triangle-walk-estimates} to bound the probability that an incidence fails the retention test \eqref{eq:incidence-test}.

Sample \((x,v)\) from \(\nu\), independently of the labels, and put \(y\defeq Y_{xv}\).
Conditional on \((x,v)\), this label has distribution \(K_x(v,\cdot)\).
Conditional on \(x\), the starting vertex \(v\) has distribution \(\pi_x\), and reversibility gives \begin{equation}\label{eq:split-label-reversibility} \pi_x(v)K_x(v,y)=\pi_x(y)K_x(y,v). \end{equation}
Thus exchanging \(v,y\) preserves their joint law conditional on \(x\), and hence also preserves the joint law of \((x,v,y)\) after averaging over \(x\).
Since \(\pi_x\) is positive on \(N(x)\), \eqref{eq:split-label-reversibility} also shows that \(K_x(y,v)>0\) whenever the sampled triple has positive probability.
Define the random variable \[ Z\defeq\log\frac{t\,w(v)}{\theta_xK_x(y,v)}\ge0, \] which is well-defined on the sampling support and nonnegative by \eqref{eq:walk-density}.
With \(K_x=K_x^j\) and \(s_x(y)=w(y)/\theta_x\), expanding the definition of \(H_j\) gives \begin{equation}\label{eq:entropy-sampling-identity} H_j =\E_{(x,v)\sim\nu}\sum_{y\in N(x)}K_x(v,y) \log\frac{w(y)/\theta_x}{K_x(v,y)} =\E_{(x,v)\sim\nu} \E_{y\sim K_x(v,\cdot)} \log\frac{w(y)}{\theta_xK_x(v,y)}. \end{equation}
Exchanging \(v,y\) by \eqref{eq:split-label-reversibility} and using \eqref{eq:entropy-sampling-identity} together with \(H_j\le C_{\mathrm{ent}}\log t\) from \eqref{eq:triangle-entropy-bounds} gives \begin{equation}\label{eq:reverse-kernel-entropy} \E Z=\E\log\frac{t\,w(y)}{\theta_xK_x(v,y)} =H_j+\log t\le(C_{\mathrm{ent}}+1)\log t, \end{equation} where both expectations are taken over \((x,v)\sim\nu\) and \(y\sim K_x(v,\cdot)\).
The incidence is rejected by \eqref{eq:incidence-test} exactly when \(Z>t'\log2+\log t\).
Markov's inequality and \eqref{eq:reverse-kernel-entropy} therefore give \begin{equation}\label{eq:incidence-rejection-probability} \Prob_{\nu,\mathrm{labels}}\bigl( (x,v)\text{ is rejected}\bigr) \le\frac{(C_{\mathrm{ent}}+1)\log t}{t'\log2+\log t} \le\frac{2(C_{\mathrm{ent}}+1)\log t}{\sqrt t}. \end{equation}

\step{Total triangle loss.}
Each of the six directed incidences of the sampled triangle has marginal \(\nu\), so \eqref{eq:incidence-rejection-probability} applies to each one.
The split creates no new triangles, and each surviving triangle corresponds to a unique triangle of \(F\) with unchanged weight.
Consequently, with probability taken over both the ordered triangle and the labels, \begin{equation}\label{eq:split-destruction-probability} \frac{\E(T-T')}{T} =\Prob_{\tau,\mathrm{labels}}\bigl(\text{the sampled triangle is destroyed}\bigr). \end{equation}
A triangle is destroyed if and only if one of its three corners disagrees or one of its six incidences is rejected.
A union bound using \eqref{eq:corner-disagreement}, \eqref{eq:incidence-rejection-probability}, and \eqref{eq:split-destruction-probability} gives \begin{equation}\label{eq:expected-split-loss} \E(T-T') \le\biggl(12\sqrt{(C_{\mathrm{ent}}+1)\log t/t} +12(C_{\mathrm{ent}}+1)\frac{\log t}{\sqrt t}\biggr)T \le C_{\mathrm{loss}}\frac{\log t}{\sqrt t}\,T. \end{equation}
The last inequality in \eqref{eq:expected-split-loss} uses \((C_{\mathrm{ent}}+1)\log t\ge1\) and the choice \(C_{\mathrm{loss}}=24(C_{\mathrm{ent}}+1)\).

\step{Completing the reduction.}
To compare this analytic estimate with the rational allowance \(\delta(t)\) in \eqref{eq:triangle-rational-loss}, note that for \(t\ge32\), \begin{equation}\label{eq:triangle-loss-comparison} 2C_{\mathrm{loss}}\frac{\log t}{\sqrt t}\le\delta(t) \le8C_{\mathrm{loss}}\frac{\log t}{\sqrt t}. \end{equation}
Our choice of \(t_*\) ensures \(\delta(t)<1/2\) for \(t>t_*\).
Using \eqref{eq:expected-split-loss} and \eqref{eq:triangle-loss-comparison}, we obtain an outcome satisfying \[ %
 T_{\rm out}\ge\E T' \ge\biggl(1-C_{\mathrm{loss}}\frac{\log t}{\sqrt t}\biggr)T \ge(1-\delta(t))T. \] %
Combining this with \(T\ge T_{\rm in}-D_{\rm in}/t\) from \eqref{eq:triangle-preparation}, and using \(0\le\delta(t)<1\), gives \[ T_{\rm out}\ge(1-\delta(t))(T_{\rm in}-D_{\rm in}/t) \ge(1-\delta(t))T_{\rm in}-D_{\rm in}/t. \]
Together with \(\Lambda'_v/\theta'_v\le2^{t'}\) at every remaining vertex and \(D_{\rm out}\le D\le D_{\rm in}\), this completes the proof of \Cref{lem:triangle-reduction}.
For rational input, \Cref{subsec:choosing-labels} constructs such an outcome deterministically.

\subsection{Iteration and the final triangle count}\label{subsec:triangle-iteration}
We now iterate the reduction to prove \Cref{lem:triangle-test}.
Write \(T_0,D_0\) for the input quantities.
If the input graph has no edges, then \(T_0=0\) and the triangle bound holds.
Otherwise, put \(t_0\defeq\max\{1,\lceil\log_2\max_v(\Lambda_v/\theta_v)\rceil\}\), as in \Cref{lem:triangle-test-algorithm}, so that \(\Lambda_v/\theta_v\le2^{t_0}\) at every vertex.
Let \(G_i\) be the graph after \(i\) reductions, with triangle weight \(T_i\), quantity \(D_i\), and parameter \(t_i\), maintaining \(\Lambda_v/\theta_v\le2^{t_i}\).
Stop if \(T_i\le2D_i\) or \(t_i\le t_*\); otherwise apply \Cref{lem:triangle-reduction}.
If it returns a witness, transfer it to the input graph using \eqref{eq:witness-projection} and stop.
In the other outcome, put \(t_{i+1}\defeq\lceil\sqrt{t_i}\rceil\) and continue.
An empty graph meets the first stopping condition.

Suppose no witness is found.
The integers \(t_i\) strictly decrease while reductions are made, so the process stops at some index \(k\).

\step{Summing the losses.}
Note that consecutive parameters \(s,t\) with \(t=\lceil\sqrt s\rceil\) satisfy \(s>(t-1)^2\ge t^{3/2}\) for \(t\ge16\).
The function \(\psi(t)\defeq\log t/\sqrt t\) is decreasing for \(t>e^2\), and for \(t\ge t_*\ge256\), \begin{equation}\label{eq:reduction-loss-decay} \psi(s)\le\psi(t^{3/2}) =\frac{3}{2}t^{-1/4}\psi(t)\le\frac{1}{2}\psi(t). \end{equation}
For \(k\ge1\), apply \eqref{eq:reduction-loss-decay} to the reduction parameters \(t_0,\ldots,t_{k-1}\), all exceeding \(t_*\).
Their reciprocals decrease by a factor of at least two at each step when the parameters are read in reverse order, since \(s>t^{3/2}\) gives \(1/s<t^{-1/2}/t\le1/(2t)\) for \(t\ge t_*\).
Summing the two geometric series gives \begin{equation}\label{eq:reduction-geometric-sums} \sum_{i<k}\psi(t_i) \le 2\psi(t_{k-1})\le2\psi(t_*), \qquad \sum_{i<k}1/t_i \le 2/t_{k-1}\le2/t_*. \end{equation}
Combining \eqref{eq:reduction-geometric-sums} with \eqref{eq:triangle-loss-comparison} gives \begin{equation}\label{eq:reduction-total-loss} \sum_{i<k}\delta(t_i)\le16C_{\mathrm{loss}}\psi(t_*)\le\frac{1}{2}, \qquad \sum_{i<k}1/t_i\le2/t_*\le1. \end{equation}
Here \(16C_{\mathrm{loss}}\psi(t_*)=64\log(100C_{\mathrm{loss}})/(10000C_{\mathrm{loss}})<1/2\); these estimates also justify the earlier bound \(\delta(t)<1/2\) for \(t>t_*\).
Since \(D_i\le D_0\), iterating \(T_{i+1}\ge(1-\delta(t_i))T_{i}-D_{i}/t_i\) from \eqref{eq:one-reduction} and using \eqref{eq:reduction-total-loss} gives \begin{equation}\label{eq:iterated-triangle-retention} T_k\ge\biggl(\,\prod_{i<k}(1-\delta(t_i))\biggr)T_0 - D_0\sum_{i<k}\frac{1}{t_i} \prod_{j=i+1}^{k-1}(1-\delta(t_j)) \ge\frac{1}{2}T_0-D_0, \end{equation} where the second inequality uses \(\prod_{i<k}(1-\delta(t_i)) \ge 1-\sum_{i<k}\delta(t_i) \ge 1/2\).

\step{The terminal graph.}
It remains to count triangles in the terminal graph \(G_k\).
If \(T_k\le2D_k\), \eqref{eq:iterated-triangle-retention} implies \(T_0\le6D_0\).
Otherwise \(t_k\le t_*\), so every remaining vertex satisfies \(\Lambda_v\le2^{t_*}\theta_v\).
Taking all local quantities in \(G_k\), we obtain \[ 3T_k=\sum_{xy\in E(G_k)}w(x)w(y)C_{xy} \le\sum_{xy\in E(G_k)}w(x)w(y)\min(\Lambda_x,\Lambda_y) \le\frac{2^{t_*}}{2}D_k \le\frac{2^{t_*}}{2}D_0. \]
Thus \(T_k\le2^{t_*}D_0/6\), whence \eqref{eq:iterated-triangle-retention} yields \(T_0\le(2+2^{t_*}/3)D_0<C_\triangle D_0\).
Thus, if no witness is found, either stopping condition gives \(T_0\le C_\triangle D_0\), since \(C_\triangle=2^{t_*}>6\).
This proves \Cref{lem:triangle-test}.
The remaining task is to choose the labels deterministically and bound the bit complexity, which we do next.


\section{Derandomizing the vertex splitting}\label{sec:derandomization}
We prove that \Cref{lem:triangle-reduction} can be implemented deterministically for rational input weights and parameters, in time polynomial in the binary input length and \(t\).
We choose the split labels by conditional expectation, retaining the triangle-loss estimate \eqref{eq:expected-split-loss}.
The method is standard (see, for example, Raghavan~\cite{raghavan1988probabilistic}), and Kleinberg and Tardos~\cite[Section~5]{kleinberg2002approximation} apply it to the labeling procedure underlying our coupling; here we compute the conditional objective exactly for the dependent label coupling in \Cref{subsec:finite-coupling}.
We first prove the coupling lemma (\Cref{lem:triangle-coupling}) using a shared sequence of random proposals.
We then express the same coupling as a bounded sequence of finite choices and choose each outcome deterministically by conditional expectation.
For the algorithmic construction, assume that \(w,\theta\) are positive rational vertex weights.

\subsection{The probabilistic coupling}\label{subsec:coupling-proof}
\begin{proof}[Proof of \Cref{lem:triangle-coupling}]
Let \(\mathcal J\) index the given family of probability vectors, recall that \(\mathcal L\) is their common finite label set, and put \(m\defeq\abs{\mathcal L}\).
For each \(n\ge1\), choose a random label \(L_n\) uniformly from \(\mathcal{L}\) and a random number \(U_n\) uniformly from \([0,1]\), with all these random variables mutually independent.
For each \(i\), define \[ T_i\defeq\min\{n\ge1:U_n\le p_i(L_n)\},\qquad Y_i\defeq L_{T_i}. \]
Thus \(T_i\) is the first index \(n\) for which \(U_n\le p_i(L_n)\), and \(Y_i\) is the label at that index.
For each \(n\), the inequality \(U_n\le p_i(L_n)\) holds with probability \(\sum_y p_i(y)/m=1/m\), so all \(T_i\) are finite almost surely.
For \(y\in\mathcal{L}\), independence gives the required marginal: \[ \Prob(Y_i=y) =\sum_{n\ge1}(1-1/m)^{n-1}\frac{p_i(y)}{m} =p_i(y). \]

Fix two members \(p_i,p_k\) of the family.
Conditional on \(L_n=y\), both inequalities \(U_n\le p_i(L_n)\) and \(U_n\le p_k(L_n)\) hold precisely when \(U_n\le\min\{p_i(y),p_k(y)\}\), and at least one holds precisely when \(U_n\le\max\{p_i(y),p_k(y)\}\).
Put \(N\defeq\min\{T_i,T_k\}\), the first index at which at least one inequality holds.
By independence across indices, \((L_N,U_N)\) has the law of \((L_1,U_1)\) conditioned on \(U_1\le\max\{p_i(L_1),p_k(L_1)\}\).
The equality \(T_i=T_k\) holds precisely when \(U_N\le\min\{p_i(L_N),p_k(L_N)\}\), so \begin{equation}\label{eq:coupling-common-acceptance} \begin{aligned} \Prob(T_i=T_k) &=\Prob\bigl( U_1\le\min\{p_i(L_1),p_k(L_1)\} \,\bigm|\,U_1\le\max\{p_i(L_1),p_k(L_1)\}\bigr)\\
    &=\frac{\frac{1}{m}\sum_{y\in\mathcal{L}}\min\{p_i(y),p_k(y)\}} {\frac{1}{m}\sum_{y\in\mathcal{L}}\max\{p_i(y),p_k(y)\}} =\frac{1-\norm{p_i-p_k}_{\TV}}{1+\norm{p_i-p_k}_{\TV}}. \end{aligned} \end{equation}
The last equality uses \(\min\{s,t\}=(s+t-\abs{s-t})/2\) and \(\max\{s,t\}=(s+t+\abs{s-t})/2\).
If \(T_i=T_k\), then \(Y_i=Y_k\), so \eqref{eq:coupling-common-acceptance} yields \[ \Prob(Y_i\ne Y_k) \le\Prob(T_i\ne T_k) %
 =\frac{2\norm{p_i-p_k}_{\TV}}{1+\norm{p_i-p_k}_{\TV}}\le2\norm{p_i-p_k}_{\TV}. \]
Since \(i,k\) were arbitrary, this proves the required disagreement bound.
\end{proof}

\subsection{Deterministic splitting by conditional expectation}\label{subsec:choosing-labels}\label{subsec:finite-realization}
We express the coupling as a bounded sequence of finite choices, then fix these choices by conditional expectation.
The resulting split retains at least the expected triangle weight of the random construction.

\step{Finite choices.}
Use the family \((p_i)_{i\in\mathcal J}\) and label set \(\mathcal L\) from the preceding proof.
A row \(i\) is \emph{pending} until it accepts a proposal and its label \(Y_i\) is assigned; thereafter the label remains fixed.
Skip every proposal accepted by no pending row, and group threshold values into intervals that assign the same labels.
Initially the set of pending rows is \(\mathcal P\defeq\mathcal J\).

While \(\mathcal P\ne\emptyset\), put \[ M_y\defeq\max_{i\in\mathcal P}p_i(y),\qquad Z\defeq\sum_{y\in\mathcal L}M_y. \]
\begin{enumerate}[label=(\arabic*)]
\item Choose a label \(y\) with probability \(M_y/Z\), among those with \(M_y>0\).
\item Choose an interval \((a,c]\) between successive distinct values in \(\{0\}\cup\{p_i(y):i\in\mathcal P\}\), with probability \((c-a)/M_y\).
\item Assign label \(y\) to every pending row \(i\) with \(p_i(y)\ge c\), and remove those rows from \(\mathcal P\).
\end{enumerate}
This is the shared proposal rule conditioned on acceptance by at least one pending row, so it preserves the label coupling.
Each outcome \((y,a,c)\) has probability \((c-a)/Z\) and assigns a label to at least one row, so the procedure takes at most \(\abs{\mathcal J}\) decisions.
After each decision, the pending rows continue with a fresh shared proposal sequence and their original thresholds \(p_i\).
Conditional on the current state, future proposals are independent of past proposals and follow the original sampling rule.
Thus each pending row \(i\) receives its final label with distribution \(p_i\), and any two pending rows have the joint label distribution defined by the original shared-proposal construction.

\step{Applying the sampler to a reduction.}
Fix a reduction satisfying the hypotheses of \Cref{lem:triangle-reduction}, with integer parameter \(t>t_*\), and put \(t'\defeq\lceil\sqrt t\rceil\).
Let \(F\) be the prepared graph from \Cref{subsec:edge-deletion}, and suppose that both tests \eqref{alg:test-one} and \eqref{alg:test-two} hold at every vertex.
Choose a common index \(0\le j<t\) minimizing the rational average on the left of \eqref{eq:nearby-rows}, as described in \Cref{subsec:triangle-entropy}.
This choice satisfies \Cref{lem:triangle-walk-estimates}; abbreviate \(K_x^j\) to \(K_x\).

At center \(x\), use the finite sampler with label set \(N_F(x)\) and one row \(K_x(v,\cdot)\) for each \(v\in N_F(x)\), denoting its final label by \(Y_{xv}\).
A local state \(\sigma_x\) records the assigned labels and the set \(\mathcal P\) of pending rows.
Let \(\Prob_{\sigma_x}\) denote probability when the sampler is continued randomly from state \(\sigma_x\), keeping all previously assigned labels fixed.

\step{Local survival probabilities.}\label{subsec:conditional-probability}
To evaluate a possible next decision, we first compute the probability that a corner survives under the resulting random continuation.
For distinct \(v,z\in N_F(x)\), let \(G_x(v,z)\) be the set of labels passing the test \eqref{eq:incidence-test} for both incidences \((x,v)\) and \((x,z)\).
The probability that the corner at \(x\) survives is
\begin{equation}\label{eq:corner-survival-probability}
    q_x^{\sigma_x}(v,z)\defeq\Prob_{\sigma_x}\bigl(Y_{xv}=Y_{xz}\in G_x(v,z)\bigr).
\end{equation}
We compute this probability in three cases:
\begin{enumerate}[label=(\roman*)]
\item If both labels are assigned, the value is \(1\) when they agree and their common label lies in \(G_x(v,z)\), and \(0\) otherwise.
\item If only the label of \(v\) is assigned, say to \(y\), the value is \(K_x(z,y)\) if \(y\in G_x(v,z)\), and \(0\) otherwise; the case with only \(z\) assigned is symmetric.
\item If both rows are pending, sum \eqref{eq:conditional-pair} below over \(y\in G_x(v,z)\).
\end{enumerate}

For the last case, put \(p\defeq K_x(v,\cdot)\), \(q\defeq K_x(z,\cdot)\), and \(Z_{p,q}\defeq\sum_y\max\{p(y),q(y)\}\).
The two-row coupling gives \begin{equation}\label{eq:conditional-pair} \Prob_{\sigma_x}(Y_{xv}=Y_{xz}=y) =\frac{\min\{p(y),q(y)\}\bigl(1+\abs{p(y)-q(y)}\bigr)}{Z_{p,q}}. \end{equation}
Indeed, at the first future proposal accepted by either row, both accept \(y\) with probability \(\min\{p(y),q(y)\}/Z_{p,q}\).
At that proposal, only the \(p\)-row accepts \(y\) precisely when the proposed label is \(y\) and the threshold satisfies \(q(y)<U\le p(y)\), an event of probability \(\max\{p(y)-q(y),0\}/Z_{p,q}\).
The pending \(q\)-row later receives \(y\) with probability \(q(y)\), giving the contribution \(\max\{p(y)-q(y),0\}q(y)/Z_{p,q}\).
By symmetry, the case in which only the \(q\)-row accepts \(y\) at that proposal contributes \(\max\{q(y)-p(y),0\}p(y)/Z_{p,q}\).
Adding the three contributions proves \eqref{eq:conditional-pair}.

\step{Choosing the outcomes deterministically.}
Process the centers in a fixed order, assigning all labels at one center before proceeding to the next.
For a global state \(\sigma\defeq(\sigma_x)_x\), continue the local samplers independently and let \(\E_\sigma\) denote expectation under these continuations.
The expected triangle weight is \begin{equation}\label{eq:conditional-triangles} \mathcal F(\sigma)\defeq\E_\sigma T' =\sum_{\{x,v,z\}\in\mathcal T(F)} w(x)w(v)w(z)q_x^{\sigma_x}(v,z)q_v^{\sigma_v}(x,z)q_z^{\sigma_z}(x,v). \end{equation}
Indeed, a triangle survives precisely when all three corners survive, and the continuations at distinct centers are independent.

At each decision, enumerate the finite sampler's outcomes \((y,a,c)\) and choose one maximizing \(\mathcal F\) at the resulting state, breaking ties in a fixed order.
By the law of total expectation, the resulting values, weighted by the sampler's outcome probabilities \((c-a)/Z\), average to the current value \(\mathcal F(\sigma)\).
Thus \(\mathcal F\) never decreases.
When all labels are assigned, construct the split graph \(F_*\) as in \Cref{subsec:finite-coupling}; its triangle weight is the terminal value of \(\mathcal F\).
Writing \(\emptyset\) for the initial state with no assigned labels, we obtain \begin{equation}\label{eq:deterministic-split} T_{F_*}(w)\ge\mathcal F(\emptyset)=\E T' \ge\biggl(1-C_{\mathrm{loss}}\frac{\log t}{\sqrt t}\biggr)T_F(w) \ge(1-\delta(t))T_F(w), \end{equation}
where the second inequality uses \eqref{eq:expected-split-loss}.
The remaining bounds in \Cref{lem:triangle-reduction} hold by the same argument as in its proof, completing the correctness proof of the deterministic reduction.

\subsection{Bit complexity and completion}\label{subsec:split-complexity}
Let \(m_0\) be the number of edges in the original input to \Cref{lem:triangle-test-algorithm}, and let \(t_0\) be its initial integer parameter.
All polynomial bounds below are in the encoding length of that input, with absolute degrees, since nothing in this section depends on \(r\).

\step{Number of operations.}
By edge injectivity, every intermediate graph has at most \(m_0\) edges and \(2m_0\) nonisolated vertices.
Since \(t_0\) is polynomially bounded and decreases at each reduction, there are polynomially many reductions.
Apart from splitting, each reduction performs the edge deletions, the \(t\) enlargement steps defining the sets \(S_t(v)\), the two tests, the computation of the kernels \(K_x^j\) for \(j<t\), and the choice of \(j\) by comparing \(t\) rational averages, all of which use polynomially many arithmetic operations.
Each splitting decision assigns a label to at least one previously unassigned incidence and has at most \(d^2\) candidate outcomes at a center of degree \(d\); the value of \(\mathcal F\) at each resulting state is computed by summing over triangles and labels using \eqref{eq:conditional-triangles} and the formulas for \eqref{eq:corner-survival-probability}.
Thus the entire procedure uses polynomially many arithmetic operations and comparisons.

\step{Encoding length.}
The transition entries and thresholds have polynomial encoding length, since weights and parameters are inherited from the input and all exponents are polynomially bounded in \(t_0\).
Choose a common denominator \(d_x\) for the entries of \(P_x\).
Then \(Q_x\) has common denominator \(2d_x\), and \(K_x^j=P_xQ_x^j\) has common denominator \(d_x(2d_x)^j\).
The bit length of this denominator grows linearly in \(j<t_0\), so the matrix powers have polynomial encoding length.
The explicit formulas for \eqref{eq:corner-survival-probability} and for \(\mathcal F\) involve polynomially many summands and bounded-depth arithmetic expressions, so their encoding lengths remain polynomial.
For each candidate outcome, \(\mathcal F\) is recomputed from the kernels \(K_x\), the assigned labels, and the sets of pending rows recorded in the resulting state, using \eqref{eq:conditional-triangles} and the formulas for \eqref{eq:corner-survival-probability}.
Its arithmetic depth therefore does not grow with the number of previous decisions.
Together with the operation count, this proves polynomial bit complexity for each reduction.

\medskip

We conclude the section by proving \Cref{lem:triangle-test-algorithm}, which gives the dichotomy algorithm.

\begin{proof}[Proof of \Cref{lem:triangle-test-algorithm}]
Apply the iteration used in the proof of \Cref{lem:triangle-test}, implementing each reduction by the deterministic procedure above.
The same stopping and witness-projection arguments apply.
The witness constructions for \eqref{alg:test-one} and \eqref{alg:test-two} give \(p=30t_i^3\) or \(p=2t_i\), respectively.
The operation and encoding-length bounds above give polynomial running time with an absolute polynomial degree.
\end{proof}


\section{Concluding remarks}\label{sec:conclusion}

In this paper, we confirm the conjecture of Kelly and Postle on fractional coloring of $K_r$-free graphs with local demands, showing that every $K_r$-free graph admits a $\phi$-coloring with $\phi(v)=\Omega\bigl(\log(d_G(v)+2)/(r(d_G(v)+1))\bigr)$.
This follows from a weighted bound on independent sets (\Cref{thm:existence}), which also yields a degree-sequence strengthening of \Cref{conj:aeks}.
Moreover, our argument is constructive: the required distributions can be computed deterministically in time \(2^{\rtconstant r^2}n^{\rtdegree}\), where \(n\) is the number of vertices and \(\rtconstant,\rtdegree\) are absolute constants.
We close with several open problems.

Our results, together with the work of OpenAI~\cite{OAI2026logarithmic} and Kang~\cite{kang2026ode}, show that \Cref{conj:aeks} holds with $c_r=\Omega(1/r)$.
For triangle-free graphs, Shearer's bound~\eqref{eq: shearer function} gives the much sharper leading constant $1-o(1)$.
More generally, Dhawan, Janzer, and Methuku~\cite{dhawan2025independent} showed that for every $3$-colorable graph $F$, every $n$-vertex $F$-free graph of average degree $d$ satisfies \[ \alpha(G)\ge(1-o(1))\frac{n\log d}{d}. \]
In the other direction, there exist $d$-regular graphs of arbitrarily large girth with independence number at most $(2+o(1))n\log d/d$~\cite{BollobasIndependence}, so the leading constant $1$ cannot be improved beyond a factor of $2$ for any forbidden graph containing a cycle.
Thus, forbidding a $3$-colorable graph costs nothing in the leading constant compared with forbidding a triangle, and $K_4$ is the smallest forbidden graph for which it is unknown whether the constant $1$ is attainable.

\begin{question}\label{q:leading-constant}
Let $r\ge4$ be fixed.
Does every $n$-vertex $K_r$-free graph with average degree $d$ satisfy \[ \alpha(G)\ge(1-o(1))\frac{n\log d}{d} \] as $d\to\infty$?
More strongly, is there a function $\varepsilon_r\colon\mathbb Z_{\ge1}\to[0,1]$ with $\varepsilon_r(d)\to0$ as $d\to\infty$ such that every $K_r$-free graph $G$ without isolated vertices admits a $\phi$-coloring for $\phi(v)\defeq\bigl(1-\varepsilon_r(d_G(v))\bigr)\log d_G(v)/d_G(v)$?
\end{question}

\Cref{q:leading-constant} bears directly on the off-diagonal Ramsey numbers $R(r,t)$.
Let $G$ be a $K_r$-free graph on $R(r,t)-1$ vertices with $\alpha(G)<t$.
The neighborhood of each vertex induces a $K_{r-1}$-free graph with no independent set of size $t$, so every degree in $G$ is less than $R(r-1,t)$.
Since $\log d/d$ is decreasing for $d\ge e$, a bound of the form $\alpha(G)\ge c_rn\log d/d$ therefore gives \[ R(r,t)\le(1+o(1))\,\frac{t\,R(r-1,t)}{c_r\log R(r-1,t)}. \]
Iterating this from Shearer's bound $R(3,t)\le(1+o(1))t^2/\log t$~\cite{shearer1983note} with $c_r=\Omega(1/r)$ recovers the order of magnitude $R(r,t)=O_r\bigl(t^{r-1}/(\log t)^{r-2}\bigr)$ of Ajtai, Koml\'os, and Szemer\'edi~\cite{ajtai1980note}.
However, the resulting leading constant is weaker than in the bound $R(r,t)\le(1+o(1))t^{r-1}/(\log t)^{r-2}$ of Li, Rousseau, and Zang~\cite{li2001asymptotic}.
In contrast, a positive answer to \Cref{q:leading-constant} would give \[ R(r,t)\le(1+o(1))\frac{t^{r-1}}{(r-2)!\,(\log t)^{r-2}}, \] improving the bound of Li, Rousseau, and Zang by a factor of $(r-2)!$.

For comparison, the case $r=3$ is now understood up to a factor of $2+o(1)$.
Shearer's upper bound is complemented by the lower bound $R(3,t)\ge(1/2+o(1))t^2/\log t$ of Hefty, Horn, King, and Pfender~\cite{hefty2025two}, which improves the earlier lower bounds in~\cite{bohman2013dynamic, fizpontiveros2020triangle, campos2025new}.
For \(r=4\), Mattheus and Verstraete~\cite{mattheus2024asymptotics} proved \(R(4,t)=\Omega(t^3/\log^4t)\), and more recently, Chen~\cite{chen2026logarithmic} and Lin and Zhang~\cite{linzhang2026logarithmic} determined the logarithmic exponent, showing that \[ R(4,t)=\frac{t^3}{(\log t)^{2+o(1)}}. \]
Both preprints rely on OpenAI's lower bound on \(R(5,t)\)~\cite{OAI:The-Sharp-Logarithmic-Exponent-of-r-5-t-September-24-2026}; Chen also uses OpenAI's lower bounds for fixed clique orders at least six~\cite{OAI:Sharp-Logarithmic-Exponents-for-Fixed-Off-Diagonal-Ramsey-Numbers-September-24-2026}, together with the independence bound in~\cite{OAI2026logarithmic}.
It remains open whether \(R(4,t)=\Theta(t^3/(\log t)^2)\); a positive answer to \Cref{q:leading-constant} would improve the leading constant in the upper bound \(R(4,t)\le(1+o(1))t^3/(\log t)^2\) by a factor of $2$.

Another question concerns the dependence on \(r\).
\Cref{thm:existence} gives \(c_r=\Omega(1/r)\), matching Kang~\cite{kang2026ode}, and it would be interesting to determine the optimal dependence of \(c_r\) on \(r\) when \(r\) is allowed to grow with \(d\).
On the algorithmic side, the running times in \Cref{thm:algorithm,cor:algorithm-local} are polynomial with degree independent of \(r\), up to a factor of \(2^{O(r^2)}\) (see \Cref{subsec:numerics-complexity}); it would also be interesting to improve this dependence on \(r\).

\subsection*{Acknowledgments}
This work was initiated during the workshop ``Flags in the Mountains,'' held at the Colorado State University Mountain Campus in June 2026.
We thank J\'{o}zsef Balogh, Bernard Lidick\'{y}, and Florian Pfender for organizing the workshop and for their invitation, and we thank all of the participants for stimulating discussions.
We are especially grateful to Clayton Mizgerd and Sam Spiro for invaluable conversations.

\subsection*{Funding information}
Ewan Davies is supported in part by NSF grant CCF-2309707.
Abhishek Dhawan is supported by NSF RTG grant DMS-1937241.
Abhishek Methuku is supported by Simons award SFI-MPS-TSM-00025583 and a UIUC Campus Research Board Award RB25050.
Jaehyeon Seo is supported by Samsung STF Grant SSTF-BA2201-02 and the National Research Foundation of Korea (NRF) grant MSIT NRF-2022R1C1C1010300.

\subsection*{Statement of AI use}
In the earlier version of this work, we proved the main result with an additional $\log\log$ factor in the denominator of the local demands.
GPT-5.6 Sol was used in a limited capacity during that work, and we completed the draft at the end of September 2026.

During the final polishing stage, OpenAI released the preprints~\cite{OAI2026logarithmic,OAI2026correspondence}.
Building on our earlier version of the work and incorporating ideas from those preprints, we found that we could remove the remaining $\log\log$ factor in our main result.
This led to our first proof of the Kelly--Postle conjecture, an existence proof with \(c_r=c/r^2\), where \(c>0\) is an absolute constant.
Using ideas from~\cite{kang2026ode}, we then obtained an algorithmic proof with the improved constant \(c_r=c/r\).
In preparing the present version, we received substantial assistance from GPT-6 Astra and GPT-6.1 Sol, as well as help from GPT and Claude in polishing the exposition.
The authors have carefully reviewed and revised the manuscript and take full responsibility for its content.


\appendix
\crefalias{section}{appendix}
\crefalias{subsection}{subappendix}
\crefalias{subsubsection}{subsubappendix}

\section{Properties of \texorpdfstring{$f$ and $U_a$}{f and U(a)}}\label{app:scalar}
Write \(U\defeq U_1\), \(\xi\defeq \xi_1\), and \(h\defeq -f'\).
The scaling identities \[ U_a(x)=aU(x/\sqrt a),\qquad \xi_a(x)=a\xi(x/\sqrt a)\qquad(a,x>0) \] follow from the definitions and reduce the assertions about \(U_a\) and \(\xi_a\) to the case \(a=1\).
The integral formulas \[ f(d)=\int_0^\infty\frac{t\,dt}{(1+t)^2(d+t)},\qquad h(d)=\int_0^\infty\frac{t\,dt}{(1+t)^2(d+t)^2} \] hold for \(d>0\).
For the first, the substitution \(t=(1-z)/z\) gives \(\int_0^1(1-z)/(1+(d-1)z)\,dz\), whose evaluation gives the formula for \(f\) in the introduction, including \(f(1)=1/2\).
Differentiation under the integral sign gives the formula for \(h\) and shows that \(f\) is smooth, positive, strictly decreasing, and strictly convex on \((0,\infty)\).
All such differentiations are justified on compact intervals of positive \(d\) by integrable bounds on the differentiated integrands.
The explicit formula gives \(f(d)\to f(0)=1\) as \(d\downarrow0\) and \(f(d)\sim\log d/d\) as \(d\to\infty\).
Continuity at zero extends the monotonicity and convexity of \(f\) to \([0,\infty)\).
Comparing the integral kernels gives \begin{equation}\label{eq:scalar-h-bound} dh(d)\le f(d)\le1\qquad(d>0). \end{equation}
By monotone convergence, \(h(d)\to\infty\) as \(d\downarrow0\), since its integrand increases to \(1/(t(1+t)^2)\), which is not integrable at zero.

Put \[ g(d)\defeq f(d)-dh(d) =\int_0^\infty\frac{t^2\,dt}{(1+t)^2(d+t)^2}>0. \]
Thus \(g\) is strictly decreasing, \(g(d)\le f(d)\le1\), and \(g(1)=1/3\).
Substituting \(t=ds\) in the integrals for \(g\) and \(h\) gives \(dg(d)=g(1/d)\) and \(d^2h(d)=h(1/d)\), so \(dg(d)\) is increasing and at most one.
Comparing \(g(d)\) for \(d\le1\) and \(dg(d)\) for \(d\ge1\) with \(g(1)=1/3\) gives \begin{equation}\label{eq:scalar-g-bounds} dg(d)\le1,\qquad g(d)\ge\frac{1}{3\max\{1,d\}}\qquad(d>0). \end{equation}
Differentiating the explicit formula for \(f\) gives \[ d(d-1)f'(d)+(d+1)f(d)=(1+d)g(d)+2dh(d)=1 \] for \(d\ne1\), and continuity gives the same identity at \(d=1\).
Solving this identity for \(h(d)\) when \(d>1\) and using the asymptotic for \(f\) gives \[ h(d)=\frac{(d+1)f(d)-1}{d(d-1)} \sim\frac{\log d}{d^2}\qquad(d\to\infty). \]

Since \(h'=-f''<0\) and \(h\) decreases from infinity to zero, \(h(d(x))=x^2\) defines a unique smooth \(d(x)>0\) for \(x>0\), by the inverse function theorem.
It uniquely minimizes \(f(t)+tx^2\) over \(t>0\), since the derivative \(x^2-h(t)\) changes from negative to positive at \(t=d(x)\).
Writing \(d=d(x)\), differentiation gives \begin{gather*} d'(x)=-\frac{2x}{f''(d)},\qquad U(x)=f(d)+dh(d),\qquad U'(x)=2dx,\\ \xi(x)=g(d),\qquad U''(x)=2d-\frac{4h(d)}{f''(d)}. \end{gather*}
In particular, \(U'(x)>0\), and taking \(t\downarrow0\) in the defining infimum gives \(U(x)\le1\), so \(0<\xi(x)\le U(x)\le1\).
The extension \(U(0)=0\) is continuous because \[ 0\le U(x)\le f(1/x)+x\longrightarrow0\qquad(x\downarrow0). \]

Differentiating the integral for \(h\) and using \(1/(d+t)<1/d\) gives \(f''(d)<2h(d)/d\).
Cauchy--Schwarz gives \[ h(d)^2\le \biggl(\int_0^\infty\frac{t\,dt}{(1+t)^2(d+t)}\biggr) \biggl(\int_0^\infty\frac{t\,dt}{(1+t)^2(d+t)^3}\biggr) =\frac{f(d)f''(d)}{2}. \]
Consequently, \[ 0<-U''(x)=\frac{4h(d)}{f''(d)}-2d \le\frac{2g(d)}{h(d)}=\frac{2\xi(x)}{x^2}. \]
Thus \(U\) is strictly concave on \((0,\infty)\), and continuity extends concavity to \([0,\infty)\).
Also, \(0<x\xi'(x)=-x^2U''(x)\le2\xi(x)\) for \(x>0\).
Scaling, with \(d_*(a,x)=d(x/\sqrt a)\), proves \eqref{eq:scalar} and \eqref{eq:scalar-ode} and gives \begin{equation}\label{eq:curvature} 0\le-U_a''(x)\le2a/x^2\qquad(a,x>0). \end{equation}
For \(0<a\le b\) and \(x>0\), taking infima over \(t>0\) in \(af(t)+tx^2\le bf(t)+tx^2\) gives \(U_a(x)\le U_b(x)\); this monotonicity in \(a\) also holds at \(x=0\), where both sides vanish.

For the tangent-loss estimate, fix \(0<x\le z\).
Since \(d(y)\) decreases with \(y\) and \(dg(d)\) increases with \(d\), the product \(d(y)\xi(y)=d(y)g(d(y))\) decreases with \(y\).
The curvature estimate therefore gives, for \(y\ge x\), \[ -y\frac{d}{dy}\bigl(U'(y)^2\bigr) =2yU'(y)(-U''(y)) \le8d(y)\xi(y)\le8d(x)\xi(x). \]
Integrating with respect to \(dy/y\) from \(x\) to \(z\) and dividing by \(U'(x)+U'(z)\ge2d(x)x>0\) gives \[ x\bigl(U'(x)-U'(z)\bigr)\le4\xi(x)\log(z/x). \]
Monotonicity of \(U'\) and the preceding bound give \[ U(x)+U'(x)(z-x)-U(z) =\int_x^z\bigl(U'(x)-U'(y)\bigr)\,dy \le4\xi(x)\frac{z-x}{x}\log(z/x). \]
If instead the target is zero, then for every \(x>0\) the loss is exactly \(U(x)-xU'(x)-U(0)=\xi(x)\).
Scaling proves \eqref{eq:scalar-loss} and the identity at target zero for every \(a>0\).

We finish with bounds used in the algorithm.
The integral representation on \([0,1]\) gives \[ f(d)\le\int_0^1\frac{dz}{1+(d-1)z} =\frac{\log d}{d-1}\le\frac{2\log d}{d}\qquad(d\ge2). \]
For \(d\ge1\), splitting the integral for \(d^2h(d)\) at \(t=d\) gives \begin{equation}\label{eq:scalar-h-log-bound} \begin{aligned} d^2h(d)&\le\int_0^d\frac{t\,dt}{(1+t)^2} +\int_d^\infty\frac{d^2\,dt}{t^3}\\
    &=\log(1+d)-\frac{d}{1+d}+\frac{1}{2} \le\log(1+d). \end{aligned} \end{equation}
For \(0<d\le1\), reciprocity and \eqref{eq:scalar-h-log-bound} at \(1/d\) give \(h(d)\le\log(1+1/d)\le\log(1/d)+\log2\).
For \(d\ge1\), monotonicity gives \(h(d)\le h(1)=1/6\), so, since \(\log2<1\), in all cases \begin{equation}\label{eq:scalar-h-log-plus-bound} h(d)\le1+\log^+(1/d)\quad(d>0),\qquad \log^+x\defeq \max\{0,\log x\}\quad(x>0). \end{equation}

\section{Bounding total variation distance by relative entropy}\label{app:probability}
We prove the form of Pinsker's inequality used in \Cref{sec:triangle-walks}.
For probability vectors \(p,q\) on a finite set, \[ \norm{p-q}_{\TV}=\frac{1}{2}\sum_y\abs{p(y)-q(y)} \le\sqrt{\KL(p\|q)/2}. \]
The claim is immediate if \(p=q\) or \(\KL(p\|q)=+\infty\).
Otherwise, finiteness implies \(p(y)=0\) whenever \(q(y)=0\), so these coordinates may be omitted.
Put \(A\defeq \{y:p(y)>q(y)\}\), \(a\defeq p(A)\), and \(b\defeq q(A)\).
Then \(a-b=\norm{p-q}_{\TV}>0\), and the support condition gives \(0<b<a\le1\).
Define \[ D_b(s)\defeq s\log\frac{s}{b}+(1-s)\log\frac{1-s}{1-b} \qquad(0\le s\le1), \] with zero terms interpreted by continuity.
Jensen's inequality for \(z\log z\), extended continuously by zero at \(z=0\), with weights \(q(y)/b\) on \(A\) and \(q(y)/(1-b)\) on its complement gives \(\KL(p\|q)\ge D_b(a)\).
Since \(D_b(b)=D_b'(b)=0\) and \(D_b''(s)=1/(s(1-s))\ge4\) for \(0<s<1\), integrating twice gives \[ D_b(a)=\int_b^a(a-s)D_b''(s)\,ds \ge2(a-b)^2\qquad(b<a<1). \]
Continuity extends this bound to \(a=1\), proving the claimed inequality in every case.

\section{Rational implementation and uniform running time}\label{app:numerics}
We construct certified rational evaluations of the potential and its derivatives without representing \(d_*\), justify the fixed rounding grid, and prove uniform running time bounds.
For the independent-set algorithm of \Cref{thm:algorithm}, \(n\ge2\), \(r\ge4\) is an integer parameter, and \(L\) is the binary input length of the graph and prescribed weights, so \(n\le L\).
The scalar routines use \(0<a\le1\) and \(\lambda\defeq(100n)^{-4}\); the remaining parameters are from \eqref{alg:coordinate-parameters}, \eqref{alg:witness-exponent}, and \eqref{alg:tolerance-parameters}.
All implied constants and polynomial degrees in this appendix are absolute.
The evaluation routines take rational inputs and tolerances and use exact rational arithmetic and comparisons.
A requested accuracy for \(U_a\) or \(U_a'\) means an absolute-error bound for the returned rational midpoint of a certified interval.
In particular, an interval of width at most \(2\varepsilon\) suffices for accuracy \(\varepsilon>0\).

\subsection{Parameter bounds and rounding}\label{subsec:numerics-rounding}
We choose a fixed rational grid for the auxiliary weights \(w\) and bound the decrease in \(\mathcal U_{H,u}\) caused by upward rounding.
The weight bound and grid spacing together control the encoding length of every stored coordinate.

The witness exponent is bounded by \(p\le P_n\) as established in \eqref{alg:witness-exponent}.
At initialization, each prescribed degree weight lies between \(2^{-L}\) and \(2^L\), so a nonisolated vertex satisfies
\[
    2^{-2L}\le d\defeq d_{H,s}(v)\le n\,2^{2L}.
\]
This rational \(d\) has encoding length polynomial in \(L\), since it is a sum of prescribed degree weights divided by one such weight.
By \eqref{eq:scalar-h-log-plus-bound}, \(u(v)\le1\), and \(d\ge2^{-2L}\), we have \(u(v)h(d)\le1+2L\log2\), so the exact initial weight satisfies \(q(v)\defeq\sqrt{u(v)h(d)}\le2(L+1)\).
Choose \begin{equation}\label{eq:numerics-rounding-grid} W_{\max}\defeq J+2(L+1)+\lambda^{-1}+10,\qquad \sigma\defeq \frac{1}{\lceil1000n^2W_{\max}/\gamma\rceil}. \end{equation}
The magnitudes and encoding lengths of these quantities are bounded uniformly in \Cref{subsec:numerics-complexity}.
By \eqref{alg:coordinate-parameters} and \eqref{alg:tolerance-parameters}, \(\gamma\le\lambda^{10}/(10^6n^4)<\min\{\lambda^3,1/n\}\), so \(0<\sigma\le\gamma/(1000n^2W_{\max})<1\).

Before each rounding in the loop, all coordinates are at most \(W_{\max}-1\): one-coordinate updates stay in \([\lambda,\lambda^{-1}]\), and two-set updates satisfy \(w'(v)\le(1+2^{(12r-14)P_n})\lambda^{-1}\le J\).
The initial grid weights constructed in \Cref{subsec:numerics-initialization} are at most \(2(L+1)\).
For a current graph \(H\), nonnegative weights \(w(v)\le W_{\max}\), and increments \(0\le e(v)\le\tau\) with \(\tau\ge0\), monotonicity of the scalar terms and expansion of the edge products give \begin{equation}\label{eq:numerics-upward-perturbation} \begin{aligned} \mathcal U_{H,u}(w+e)-\mathcal U_{H,u}(w) &\ge-2\sum_{xy\in E(H)} \bigl(w(x)e(y)+w(y)e(x)+e(x)e(y)\bigr)\\
    &\ge-n^2(2W_{\max}\tau+\tau^2), \end{aligned} \end{equation} since \(H\) has at most \(n^2/2\) edges.
Round positive coordinates upward to multiples of \(\sigma\), leaving zero coordinates zero, and use this same grid throughout the run.
Applying \eqref{eq:numerics-upward-perturbation} with \(\tau=\sigma\), using \(\sigma\le1\le W_{\max}\), bounds the decrease in \(\mathcal U_{H,u}\) caused by this rounding: \begin{equation}\label{eq:numerics-rounding-loss} \text{rounding loss}\le4n^2W_{\max}\sigma \le\frac{\gamma}{250} <\min\{\gamma/2,\lambda^3/2\}. \end{equation}
If a proposed coordinate lies outside \([\lambda,\lambda^{-1}]\), set it to the corresponding endpoint and immediately discard or select the vertex as in \Cref{subsec:algorithm-weights}, without rounding.
The changed coordinate is removed, and all surviving coordinates were already on the grid.
Thus stored coordinates have the form \(k\sigma\) with \(0\le k\le\lceil W_{\max}/\sigma\rceil\), so each coordinate uses \(O(1+\log W_{\max}+\log(1/\sigma))\) bits.

\subsection{Rational evaluation and error bounds}\label{subsec:numerics-evaluation}
For rational \(0<a\le1\), \(x\ge\lambda/2\), and a prescribed rational \(0<\varepsilon\le1\), we compute \(U_a(x)\) and \(U_a'(x)\) to certified absolute error at most \(\varepsilon\), without representing the minimizer \(d_*(a,x)\).
We then specify the scalar accuracies needed for \(\zeta_v,\xi_v\), potential values, and selection gains in the algorithm.

\step{Function values.}
Recall that, for \(x>0\),
\[
    U_a(x)=\min_{d>0}\bigl(af(d)+dx^2\bigr),
\]
with minimizer \(d_*(a,x)\).
We approximate \(f\) by a rational sum and then minimize over a finite rational grid.
For an integer \(N\ge2\), set \begin{equation}\label{eq:numerics-f-approximation} f_N(d)\defeq \frac{1}{N}\sum_{j=1}^{N-1}\frac{N-j}{N-j+dj}. \end{equation}
In the representation \[ f(d)=\int_0^1\frac{1-z}{1+(d-1)z}\,dz, \] the integrand decreases from one to zero for \(d>0\), and \(f_N\) is its right Riemann sum.
The difference between the left and right sums is \(1/N\), so, also at \(d=0\) by continuity, \begin{equation}\label{eq:numerics-f-approximation-error} 0\le f(d)-f_N(d)\le1/N\qquad(d\ge0). \end{equation}
Differentiating the finite sum gives \(\abs{f_N'(d)}\le N\) for all \(d\ge0\).

Suppose \(x\ge\lambda/2\).
The identity \(ah(d_*)=x^2\) and the bound \(d_*h(d_*)\le1\) give \(d_*\le a/x^2\le4/\lambda^2\), so we may minimize over \([0,4/\lambda^2]\), interpreting \(f(0)=1\).
For a prescribed rational \(\varepsilon>0\), choose
\[
    N\defeq\max\{2,\lceil1/\varepsilon\rceil\},\qquad
    \delta\defeq\frac{\varepsilon}{aN+x^2+1}.
\]
Let \(m\) be the minimum of \(af_N(d)+dx^2\) on the equally spaced grid including both endpoints of \([0,4/\lambda^2]\) with \(\lceil4/(\lambda^2\delta)\rceil\) intervals, whose mesh is at most \(\delta\).
Replacing \(f\) by \(f_N\) decreases the objective by at most \(a/N\le1/N\).
The derivative of \(af_N(d)+dx^2\) has absolute value at most \(aN+x^2\), so moving from \(d_*\) to a grid point within \(\delta\) changes this objective by at most \((aN+x^2)\delta\).
Comparison at that grid point gives \(m\le U_a(x)+(aN+x^2)\delta\), while comparison at a grid minimizer gives \(U_a(x)\le m+1/N\).
Thus the approximation error \(1/N\) and the discretization error \((aN+x^2)\delta\) give the enclosure
\[
    m-(aN+x^2)\delta\le U_a(x)\le m+1/N.
\]
Its width is at most \((aN+x^2)\delta+1/N\le2\varepsilon\), so returning its rational midpoint gives absolute error at most \(\varepsilon\).

\step{Derivatives.}
For \(x\ge\lambda/2\), the curvature bound \eqref{eq:curvature} gives \(\abs{U_a''}\le8/\lambda^2\).
Thus the forward difference \((U_a(x+b)-U_a(x))/b\) has error at most \(4b/\lambda^2\).
Approximating each function value to error at most \(\varepsilon_{\mathrm{val}}\) adds at most \(2\varepsilon_{\mathrm{val}}/b\).
For \(0<\varepsilon\le1\), take
\[
    b\defeq\frac{\varepsilon\lambda^2}{16},\qquad
    \varepsilon_{\mathrm{val}}\defeq\frac{\varepsilon b}{8},\qquad
    \frac{4b}{\lambda^2}+\frac{2\varepsilon_{\mathrm{val}}}{b}=\frac{\varepsilon}{2}.
\]
The computed forward difference plus or minus \(\varepsilon/2\) encloses \(U_a'(x)\) in an interval of width \(\varepsilon\).

\step{Error propagation.}
Under \eqref{alg:box}, evaluating the scalar values, derivatives, and potentials to the accuracies specified in the middle column of the following table gives the error bounds in its right column.
Here \(U,U'\) denote the scalar value and derivative at each relevant vertex; edge terms, \(\Lambda_v\), \(u\), and \(c\) are rational and exact, with \(0<c\le1\); \(\eta\) is defined in \eqref{alg:tolerance-parameters}.
\[ \begin{array}{c|c|c} \text{Quantity}&\text{Sufficient evaluation error bounds}&\text{Resulting error bound}\\ \hline \zeta_v&U':\eta/8&\eta/8\\ \xi_v&U:\lambda^3/12,\quad U':\lambda^4/12&\lambda^3/6\\ \mathcal U_{K,u}&\text{each scalar term}:\varepsilon/(4n)&\varepsilon/4\\ g_x&\text{each potential}:\varepsilon/4&c\varepsilon/2\le\varepsilon/2 \end{array} \]
The error in \(\xi_v=U-w(v)U'\) is at most \(\lambda^3/12+w(v)\lambda^4/12\le\lambda^3/6\).
For any induced subgraph \(K\) of \(H\), potential errors sum over at most \(n\) scalar terms, and the defining expression for \(g_x\) uses two potentials with coefficients \(c\) and \(-c\).
To obtain rational parameters \(\theta_v=2r\bar\xi_v/w(v)\), set \begin{equation}\label{eq:numerics-xi-approximation} \abs{\widehat\xi_v-\xi_v}\le\lambda^3/6,\qquad \bar\xi_v\defeq\widehat\xi_v+\lambda^3/6. \end{equation}
The lower bound \eqref{alg:xi-lower} gives \begin{equation}\label{eq:numerics-xi-upper-bound} \xi_v\le\bar\xi_v\le\xi_v+\lambda^3/3\le2\xi_v. \end{equation}
For selection gains, use \(\varepsilon\defeq\lambda^3/(100n)\); the approximation to \(\mathcal U_{H,u}(w)\) can be reused for every vertex.

The discard and selection tests in steps~\ref{step:weighted-discard} and~\ref{step:weighted-select} establish \eqref{alg:box} before any call to \(U_a\) or \(U_a'\); initialization instead evaluates \(h\).
Thus every positive evaluation argument is at least \(\lambda\), and zero coordinates use \(U_a(0)=0\) directly.
All these arguments, including the additional evaluation points \(x+b\), are at most \(W_{\max}+1\), since stored coordinates are at most \(W_{\max}+\sigma\) with \(\sigma<1\) and \(b\le\lambda^2/16\).
The number of grid points and the bit cost of the rational sums depend polynomially on the input encoding lengths, \(\lambda^{-1}\), \(W_{\max}\), and the reciprocals of the requested accuracies, with absolute degrees.
\Cref{subsec:numerics-complexity} bounds all these quantities for the actual calls, including the function-value accuracies used in forward differences.

\subsection{Initializing from the prescribed weights}\label{subsec:numerics-initialization}
For each nonisolated vertex \(v\), we construct a positive initial weight \(w(v)\) on the fixed \(\sigma\)-grid satisfying \(q(v)\le w(v)<q(v)+2\sigma\), where \(q(v)=\sqrt{u(v)h(d_{H,s}(v))}\) is the exact initial weight from \Cref{lem:lift}.
The construction takes time polynomial in \(L+\log(1/\sigma)\), and its potential loss is bounded by \eqref{eq:numerics-initialization-loss}.

\step{Certified evaluation of \(h\).}
For rational \(d>0\), we first enclose \(h(d)\) in a rational interval of prescribed width at most \(\varepsilon>0\).
Differentiating the formula for \(f\) gives
\[
    h(d)=\frac{(d+1)\log d-2(d-1)}{(d-1)^3}\quad(d\ne1),
    \qquad h(1)=\frac{1}{6}.
\]
For \(d=1\), return the exact value \(1/6\).
For \(0<\abs{d-1}\le1/2\), put \(z\defeq d-1\).
The geometric expansion of \((1-t)/(1+zt)\) and its derivative in \(z\) converge uniformly for \(0\le t\le1\) and \(\abs{z}\le1/2\), so termwise integration and differentiation give
\[
    h(1+z)=\sum_{j=0}^{\infty}
    \frac{(-1)^j(j+1)}{(j+2)(j+3)}z^j.
\]
For any integer \(m\ge1\), the tail after \(0\le j<m\) has absolute value at most \(2^{1-m}\), since each coefficient has absolute value at most one.
The truncated sum plus or minus this bound has width \(2^{2-m}\); choose \(m\) to make this at most \(\varepsilon\).

For \(\abs{d-1}>1/2\), write \(d=2^kz\), where \(k\) is an integer and \(1\le z<2\).
Use \(\log d=k\log2+\log z\) and
\[
    \log z=2\sum_{j=0}^{\infty}
    \frac{1}{2j+1}\biggl(\frac{z-1}{z+1}\biggr)^{2j+1}.
\]
For \(1\le z\le2\) and any integer \(m\ge1\), the terms are nonnegative and the tail after \(m\) terms is at most \(2\cdot3^{-(2m+1)}/(1-1/9)\).
The partial sum and the partial sum plus this bound provide an interval for \(\log z\) and, on setting \(z=2\), for \(\log2\).
Since
\[
    \abs[\bigg]{\frac{d+1}{(d-1)^3}}
    \le\frac{1}{\abs{d-1}^2}+\frac{2}{\abs{d-1}^3}
    \le20,
\]
intervals of width at most \(\varepsilon/[20(\abs{k}+1)]\) for the two logarithms give width at most \(\varepsilon/20\) for \(\log d\), and hence at most \(\varepsilon\) for \(h(d)\).
Use rational interval arithmetic, reversing endpoints when multiplying by a negative number, as can occur for \(k<0\) or for the coefficient \((d+1)/(d-1)^3\) when \(d<1\).
The remaining term \(-2/(d-1)^2\) in the formula for \(h\) is rational and exact.
If \(\beta\) is the binary encoding length of \(d\), then \(\abs{k}=O(\beta)\) in the logarithmic case.
In both cases, the geometric tail bounds give a number of terms and rational encoding lengths polynomial in the input lengths and \(\log^+(1/\varepsilon)\), proving polynomial bit complexity for the interval evaluation of \(h\).

\step{Choosing the initial grid weights.}
At each nonisolated vertex, evaluate \(h(d)\) with interval width at most \(\sigma^2\), multiply the interval by \(u(v)\), and intersect it with \([0,4(L+1)^2]\).
By \(0<u(v)\le1\) and the bound on \(q(v)\) in \Cref{subsec:numerics-rounding}, the resulting nonempty interval \([\underline y,\overline y]\) contains \(q(v)^2=u(v)h(d)\) and has width at most \(\sigma^2\).
Use integer binary search to find the least \(j\) with \(0\le j\le2(L+1)/\sigma\) and \((j\sigma)^2\ge\overline y\), and set \(w(v)\defeq j\sigma\).
The search upper bound is an integer because \(1/\sigma\) is an integer, and is feasible because \(\overline y\le4(L+1)^2\).
Since \(\overline y\ge q(v)^2>0\), we have \(j\ge1\), so \(w(v)>0\).
Minimality of \(j\) gives
\[
    0\le w(v)-q(v)
    <\sigma+\sqrt{\overline y-\underline y}
    \le2\sigma,
    \qquad w(v)\le2(L+1).
\]
Thus the initial weights lie directly on the fixed grid, and \eqref{eq:numerics-upward-perturbation} with \(\tau=2\sigma\) bounds the initialization loss by
\begin{equation}\label{eq:numerics-initialization-loss}
    8n^2W_{\max}\sigma
    \le\frac{\gamma}{125}<\frac{1}{100n}.
\end{equation}
The integer search uses \(O(\log(L+1)+\log(1/\sigma))\) exact rational square comparisons.
Since \(d\) has encoding length polynomial in \(L\) and the tolerance \(\sigma^2\) has \(O(1+\log(1/\sigma))\) bits, the interval evaluation and the integer search have total bit cost polynomial in \(L+\log(1/\sigma)\).

\subsection{Uniform running time bound}\label{subsec:numerics-complexity}
We prove the uniform running time bounds in \Cref{thm:algorithm,cor:algorithm-local} by controlling rational encoding lengths and operation counts.
We first treat the independent-set algorithm, then apply its bound to the distribution construction.

Put
\begin{equation}\label{eq:numerics-runtime-factor}
    \rtfactor\defeq2^{(12r-14)P_n}.
\end{equation}

\step{Parameter magnitudes and stored weights.}
Set \(\rtcontrol\defeq L+r+\rtfactor+1\ge2\).
There is an absolute integer \(\rtpower\ge1\), independent of \(n,r,L\), such that
\begin{equation}\label{eq:numerics-uniform-parameters}
    \max\bigl\{\lambda^{-1},J,\eta^{-1},\ell^{-1},
    \gamma^{-1},W_{\max},\sigma^{-1}\bigr\}
    \le\rtcontrol^{\rtpower}.
\end{equation}
This follows from the definition of \(\lambda\) and the parameter definitions in \eqref{alg:coordinate-parameters}, \eqref{alg:tolerance-parameters}, and \eqref{eq:numerics-rounding-grid}, using \(n\le L\).
All quantities in this maximum are positive integers by their definitions, so they and their reciprocals have \(O(\log\rtcontrol)\)-bit encodings.
By \Cref{subsec:numerics-rounding}, every stored grid coordinate also has \(O(\log\rtcontrol)\) bits.
The normalized objective weights \(u(v)\) are computed with polynomial bit cost in \(L\), have encoding lengths polynomial in \(L\), and remain fixed throughout the independent-set loop.

\step{Scalar evaluation and initialization.}
With the selection-gain accuracy fixed at \(\lambda^3/(100n)\), every error tolerance prescribed in the table of \Cref{subsec:numerics-evaluation} has reciprocal bounded by \(\rtcontrol^{O(1)}\), by \eqref{eq:numerics-uniform-parameters}.
For either derivative tolerance \(\varepsilon\in\{\eta/8,\lambda^4/12\}\), the forward-difference procedure uses function-value tolerance \(\varepsilon_{\mathrm{val}}=\varepsilon^2\lambda^2/128\), whose reciprocal satisfies the same bound.
The choices of \(N,\delta\) in \Cref{subsec:numerics-evaluation}, together with \(x\le W_{\max}+1\), therefore bound \(N\) and the number of grid points in every function-value evaluation by \(\rtcontrol^{O(1)}\).

Each grid point and individual summand has encoding length polynomial in \(L+\log\rtcontrol\), by the rational formulas in \Cref{subsec:numerics-evaluation} and the encoding bounds above for the parameters, stored weights, and scalar parameter \(a=u(v)\).
Adding \(N-1\) rational summands increases the encoding length by at most a factor proportional to \(N\), by using the product of their denominators as a common denominator.
At each grid point, we evaluate \(f_N\) by summing \(N-1\) rational terms and compare the resulting rational objective values to find the grid minimum.
Each iteration makes \(O(n^2)\) scalar evaluations, including the function-value calls within forward differences: the derivative and \(\xi_v\) computations use \(O(n)\) calls, and evaluating all selection gains separately uses \(O(n^2)\) calls.
Together with the grid counts and encoding-length bounds above, this gives total bit cost \(\rtcontrol^{O(1)}\) for all scalar evaluations in each iteration, using exact integer arithmetic and comparisons.

By \Cref{subsec:numerics-initialization}, initialization has cost polynomial in \(L+\log(1/\sigma)\), also bounded by \(\rtcontrol^{O(1)}\).

\step{Structural calls and the outer loop.}
Each dichotomy input \((H,w,\theta)\) has total encoding length \(\rtcontrol^{O(1)}\): \(H\) has at most \(n\le L\) vertices, and the scalar bounds above give \(\rtcontrol^{O(1)}\)-bit encodings for \(\theta_v=2r\bar\xi_v/w(v)\).
The uniform bit bound in \Cref{subsec:split-complexity} therefore gives cost \(\rtcontrol^{O(1)}\) for each call and also bounds the encoding length of its returned witness by \(\rtcontrol^{O(1)}\).
The graph \(H\), weights \(u,w,\theta\), and returned witness thus give an input of \(\rtcontrol^{O(1)}\) bits to the reweighting construction in \Cref{lem:improvement}.
That construction has a uniform polynomial bit bound in its input length and the numerical value of \(p\le P_n=O(\log\rtfactor)\), so its bit cost is also \(\rtcontrol^{O(1)}\).
Rounding then returns every surviving coordinate to the fixed grid, so intermediate encoding lengths do not accumulate across iterations.
Since \(c^{-1}=O(r)\), the iteration bound in \Cref{subsec:algorithm-analysis} and \eqref{eq:numerics-uniform-parameters} give \(O(r\gamma^{-1}+n)=\rtcontrol^{O(1)}\) iterations.
The total bit cost of the independent-set algorithm is therefore \(\rtcontrol^{O(1)}\).

\step{Absorbing the witness factor.}
For \(b_n=\lceil\log_2(2n)\rceil\ge1\), the definition of \(q_n\) gives \(q_n=O(1+\log b_n)=O(b_n^{1/6})\).
Since \(P_n=30q_n^3\) by \eqref{alg:witness-exponent}, it follows that \(P_n=O(\sqrt{b_n})\) and \(\log_2\rtfactor=O(r\sqrt{b_n})\).
Applying \(2xy\le x^2+y^2\) to an absolute constant times \(r\sqrt{b_n}\) gives
\[
    \log_2\rtfactor\le b_n+O(r^2),\qquad
    \rtfactor\le4n\,2^{O(r^2)}.
\]
Here \(2^{b_n}\le4n\).
Substitution in the bound \(\rtcontrol^{O(1)}\), using \(L\ge1\), \(n\le L\), and absorbing polynomial factors in \(r\) into \(2^{O(r^2)}\), gives \(2^{O(r^2)}L^{O(1)}\), with absolute constants in both exponents.

\step{The distribution algorithm.}
For this part let \(n\ge1\) and use the construction in \Cref{subsec:construct-distribution}.
It makes exactly \(m=32\lceil8n/c_r\rceil b_n=O(rn\log(2n))\) calls to \Cref{thm:algorithm}.
The inclusion-frequency guarantee was proved there; here we bound the bit cost.
Let \(\beta_a\) be the maximum binary encoding length of a target \(a_v\).
Each value \(f_{4n}(d_G(v))\) is a sum of \(O(n)\) rationals with \(O(\log(2n))\)-bit encodings, and \(c_r\) has \(O(\log(r+1))\) bits.
Thus the targets \(a_v=(c_r/2)f_{4n}(d_G(v))\) are exactly computable in time polynomial in \(n,r\), with absolute degree, and satisfy
\[
    \beta_a=O(n\log(2n)+\log(r+1)).
\]
Each update factor in \eqref{eq:rational-mw-update} has \(O(\beta_a+1)\) bits.
Consequently, after \(t\) updates, the unnormalized objective weight \(z_v^{(t)}\) has \(O(1+t(\beta_a+1))\) bits.
Let \(L_t\) denote the graph-and-weight input length of the call at round \(t\), which uses \(u(v)=z_v^{(t)}\) and \(s\equiv1\).
Since the graph has an encoding of length polynomial in \(n\), the preceding weight bound and \(m=O(rn\log(2n))\) show that \(L_t\) is bounded by a polynomial of absolute degree in \(n,r\) for every \(0\le t<m\).
There are \(nm\) rational weight updates, each on numbers of polynomial encoding length in \(n,r\), so all updates together also have polynomial bit cost with absolute degree.
Applying the independent-set running time bound with \(L=L_t\) and summing over the \(m\) calls gives total bit cost \(2^{O(r^2)}n^{O(1)}\); the polynomial factors in \(r\) are absorbed into \(2^{O(r^2)}\).
For \(n\le1\), the independent-set algorithm returns after the initial checks in \Cref{sec:outer}, with bit cost polynomial in \(L\) of absolute degree, which is within the same bound.
Choosing common absolute constants \(\rtconstant,\rtdegree>0\) for both algorithms yields the bounds in \Cref{thm:algorithm,cor:algorithm-local}.

\printbibliography

\end{document}